\documentclass[graybox,envcountchap]{svmono}

\usepackage{type1cm}         

\usepackage{makeidx}         
\usepackage{graphicx}        
\usepackage{multicol}        
\usepackage[bottom]{footmisc}

\usepackage{newtxtext}       %
\usepackage[varvw]{newtxmath}       

\usepackage[psdextra, pdfpagelabels, colorlinks,linkcolor= blue, citecolor= blue]{hyperref}

\usepackage[titletoc]{appendix}

\usepackage{cancel}
\usepackage[normalem]{ulem}

\usepackage{lipsum}
\usepackage[most]{tcolorbox}

\definecolor{blue}{RGB}{0, 0, 0}

\newcommand{\bleu}[1]{{\color{blue}#1}}

\newcommand\magentaout{\bgroup\markoverwith
	{\textcolor{magenta}{\rule[.5ex]{2pt}{.7pt}}}\ULon}

\newcommand\redout{\bgroup\markoverwith
	{\textcolor{violet}{\rule[.5ex]{2pt}{.7pt}}}\ULon}

\spnewtheorem{exemple}[example]{Example}{\itshape}{\rmfamily}

\newcommand{\R}{\mathbb{R}}

\newcommand{\N}{\mathbb{N}}

\newcommand{\n}{\noindent}
\newcommand{\m}{\medbreak}

\graphicspath{{IMAGES/}{IMAGES/OBS/}}

\makeindex             

\begin{document}

\author{N.~Cunniffe, F.~Hamelin, A.~Iggidr, A.~Rapaport, G.~Sallet}
\title{Identifiability and Observability in Epidemiological Models}
\subtitle{-- a survey --}
\maketitle

\frontmatter

%
%
%

\begin{dedication}
The World requires at least ten years to understand a new idea, however important or simple it may be.

\medskip

Ronald Ross (1902 Nobel Prize).
\end{dedication}

%
%

\preface

In Mathematical Epidemiology, many papers present the following structure:

\begin{itemize}
	\item[-]  a model is proposed,
	\item[-]  some parameters are given, extracted from literature,
	\item[-]  remaining unknown parameters are estimated by fitting the model to some observed data. 
\end{itemize}
Fitting is done usually by using an optimization algorithm with the use for example of a least square method or  a maximum likelihood estimation.
To validate the parameters estimation, one can use noisy synthetic simulated data obtained from the model for given values of the parameters, to check that the algorithm is able to reconstruct from the data the values of these parameters with accuracy.\m

One objective of this book is to show that this procedure is not always safe 
and that an examination of the identifiability of parameters is a prerequisite before a numerical determination of parameters. We will review different methods to study identifiability and observability and then consider the problem of numerical identifiability. Our touchstone will be the most famous, however simple, model in Mathematical Epidemiology,  the SIR model of Kermack and Mckendrick  \cite{KmK1927}. This model gets a renewed attention with the 
COVID-19 pandemic \cite{MR4108968,Roda:2020aa}. 
Parameter identifiability analysis addresses the problem of which unknown parameters of an ODE model can  uniquely be recovered from observed data. We will show that, even for very simple models, identifiability is far from being guaranteed.\m

The problem of identifiability for  epidemiological models is relatively rarely addressed. For instance, a research in the Mathematical Reviews of the American Mathematical Society\footnote{\tt https://mathscinet.ams.org/mathscinet} en 2020 with \texttt{epid*  AND identifiability} gives only 4 papers,  while \texttt{epidem* AND  parameter} returns  68 publications. Only a small part of the later publications address the problem of identifiability. In particular, the following publications consider the problem of identifiability in epidemiological models: \cite{Anguelova2007,Audoly:2001aa,Bellu:2007aa,Chis:2011aa,Chis:2011ab,Eisenberg:2013aa,MR2142487,Jacquez-Greif-MBS85,Lintusaari:2016aa,MR2785878,Nguyen:2015aa,MR2726194,MR2817815,MR2019304,MR3553947,MR3784444,MR3825986,Villaverde:2016aa,MR2393024,MR1957979}. However the majority of these papers is published  elsewhere than in Biomathematics journals.

\m

The question of observability, i.e.~the ability to reconstruct state variables of the model from measurements, is often considered separately from the problem of identifiability. Either model parameters are known, or the identifiability analysis is performed prior to the study of observability. Indeed, the concepts of identifiability and observability are closely related, as we show in this book. However, for certain models, it is possible to reconstruct state variables with observers, while the model is not identifiable. In other situations, we show that considering jointly identifiability and observability with observers can be a way to solve the identifiability problem. This is another illustration of the interest of observers. This is why we shall dedicate a fair part of this monograph reviewing the concept of observers and their practical constructions in epidemiology.

\m

This book is aimed at scientists, researchers, graduate students, who use or develop mathematical models for epidemiology, and who are not yet familiar with the concepts of control science (detectability, observability, observers) applied to this field.

\vspace{\baselineskip}
\begin{flushright}\noindent
Cambridge,\hfill {\it Nik Cunniffe}\\
Rennes,\hfill {\it Frédéric Hamelin}\\
Metz,\hfill {\it Abderrahman Iggidr, Gauthier Sallet}\\
Montpellier,\hfill {\it Alain Rapaport}\\
May 2023 \hfill 
\end{flushright}

%
%

\extrachap{Acknowledgements}

The authors are deeply grateful to P.-A.~Bliman, C.~Lobry, J.~Harmand, T.~Sari, M.~Sofonea, M.~Souza and many other colleagues or students, for exchanges and fruitful discussions that gave them the willingness to write this monograph.

\tableofcontents

\mainmatter

\chapter{Introduction}
\label{chap:intro} 

\abstract*{We introduce the basic concepts of identifiability and observability, and present the structure of the book.}

\abstract{We introduce the basic concepts and present the structure of the book.}

\section{Definitions}
\label{secdef}
The question of parameter identifiability originates from control theory and is related to observability and controllability  concepts \cite{Son90}.   The first appearance is in Kalman \cite{MR0152167} and is now sixty years old.
Identifiability is related to observability: the observability of a model is the ability to reconstruct the state of a system from the observation (we give a precise definition below). In the language of dynamical systems with inputs and outputs,
which is the standard paradigm in control systems theory, an input-output relation is defined. 
\bleu{Typically, models take the form
\[
\left\{\begin{array}{l}
\dfrac{dx(t)}{dt} = f(x(t),u(t)),\\
\\
y(t)=h(x(t))
\end{array}\right.
\]
where $u(t)$ and $y(t)$ are respectively "input" and "output" vectors at time $t$, and the vector $x(t)$ represents the internal variables of the model at ime $t$.
The inputs, represented by the function $u(\cdot)$, also called "controls" or "control variables" are considered as known. The outputs, represented by the function $y(\cdot)$, are called "observations" or "measurements" and are also considered as known. For simplicity, we will only consider systems without control (which is a peculiar case where the controls take values  in a singleton...).} When controls are known, with more information, observability/identifiability is sometimes easier.  These problems have rarely been considered for uncontrolled systems
whereas many methods have been developed for controlled systems.
To be more precise, let us consider a dynamical system defined in a \bleu{open domain
${\cal D} \subset \R^n$}
\begin{equation}\label{sys}
	\Sigma: \,\left \{
	\begin{array}{rl}
		\dot x (t)&=  f(x(t)), \; x(0)=x_0 \in {\cal D},\\
		\\
		y (t) &=h(x(t)), \\
	\end{array}
	\right.
\end{equation}
where we have denoted  $\dot x(t)=\dfrac{dx(t)}{dt}$.
The  ordinary differential equation (ODE)  $\dot x =f(x)$ is  the dynamics  and  $x$ is called the state of the system. To  avoid technical details we will assume in all the book that for any initial condition  $x_0$ in ${\cal D}$, there exists an unique solution denoted  $x(t,x_0)$ such that $x(0,x_0)=x_0$  and 
\[
\frac{dx(t,x_0)}{dt}\, = f(x(t,x_0)). 
\] 
Moreover, we will assume that this solution  $x(t,x_0)$ is defined for any time  $t\geq 0$, and \bleu{we shall consider a connected subset $\Omega$ of ${\cal D}$ of non-empty interior} that is a positively invariant, which means that for any initial condition $x_0 \in \Omega$, the solution $x(t,x_0)$ belongs to this set for any $t\geq 0$. This is often the case with epidemiological models for which the state vector $x$ naturally evolves in a compact and connected invariant set with non-negative vectors. This situation is also often encountered  in biological systems.
\bleu{Throughout the manuscript, $\Omega$ will then denote the state space.}

\medskip

The output (or ``observation'') of the system is given by  $h(x)$ where $h$ is a differentiable function  $h :  x \in  \Omega \subset \R^n \mapsto h(x) \in Y \subset \R^m$ . The set  $Y$ is the \bleu{output} space. We will denote by 
$h(t,x_0)$ or  $y(t,x_0)$ the observation at time $t$ for an initial condition $x_0$.

\begin{definition}[Observability]\label{def:obs-intro}
	The system  (\ref{sys}) is  observable on $\Omega$ if for two distinct initial  states $x_0$, $x_0^\prime$ in $\Omega$, there exists a time  $t\geq 0$ such that
	\[
	h(x(t,x_0)) \neq h(x(t,x_0^\prime))
	\]
\bleu{This is equivalent to state
\[
\Big\{ h(x(t,x_0)) = h(x(t,x_0^\prime)),\ \forall t\geq 0 \Big\} \Rightarrow x_0=x_0^\prime.
\]
}

Two states $x_0$, $x_0^\prime$ in $\Omega$ are called indistinguishable if we have
	\[
	h(x(t,x_0)) = h(x(t,x_0^\prime)), \quad t \geq 0
	\]
\end{definition}
Indistinguishability means that it is impossible to distinguish the  evolution of the system,   from two distinct initial conditions, by considering only  the observation $y(\cdot)$.\m
\bleu{
Note that the observability property may depend on the choice of the (positively) invariant set $\Omega$. Consider for instance the system
\[
\left\{\begin{array}{l}
\dot x_a = -x_ax_b\\
\dot x_b = x_ax_b\\
\\
y=x_b
\end{array}\right.
\]
It is not observable for $\Omega=\{ x_a \geq 0, \; x_b \geq 0, \; x_a+x_b \leq 1\}$ because all initial conditions of the form $(x_a,0)$ cannot be distinguished, while it is observable for $\Omega=\{ x_a > 0, \; x_b > 0, \; x_a+x_b \leq 1 \}$. If two solutions $(x_a(\cdot),x_b(\cdot))$, $(x_a^\prime(\cdot),x_b^\prime(\cdot))$ generate the same output $y(\cdot)$, one should have
\[
\left\{\begin{array}{l}
y(t)=x_b(t)=x_b^\prime(t),\\
\dot y(t)=x_a(t)y(t)=x_a^\prime(t) y(t),
\end{array}\right. \qquad t \geq 0
\]
but solutions in $\Omega$ are such that $y(t)\neq 0$ at any $t\geq 0$, which implies $x_a(t)=x_a^\prime(t)$ and thus $(x_a(0),x_b(0))=(x_a^\prime(0),x_b^\prime(0))$.

\medskip

Sometimes, we shall simply say that "the system is observable", when $\Omega={\cal D}$ or when there is no ambiguity about the choice of $\Omega$.
}

\bigskip

Now we consider a system depending on a parameters vector  $\theta \in \Theta \subset \R^p$ (sometimes simply called the "parameter $\theta$")
\begin{equation}\label{sys2}
	\left \{
	\begin{array}{rl}
		\dot x (t)&=  f(x(t), \theta), \; x(0)=x_0,\\
		\\
		y (t) &=h(x(t, \theta)). \\
	\end{array}
	\right.
\end{equation}
\color{blue}
We denote by 
$x(t,x_0, \theta)$ the  solution of  (\ref{sys2}) for an initial condition  $x_0$.
Then, the definition of observability needs to be adapted as follows
\begin{definition}[Observability for parameterized systems]\label{def:obs-param-intro}
	The system  (\ref{sys2}) is  observable on $\Omega$ if whatever are $ \theta \in \Theta$ and two distinct initial conditions $x_0$, $x_0^\prime$ in $\Omega$, there exists a time  $t\geq 0$ such that
	\[
	h(x(t,x_0,\theta)) \neq h(x(t,x_0^\prime,\theta ))
	\]
or equivalently 
\[
\Big\{ h(x(t,x_0,\theta)) = h(x(t,x_0^\prime,\theta)),\ \forall t\geq 0 \Big\} \Rightarrow x_0=x_0^\prime.
\]
\end{definition}
\color{black}

\medskip

Identifiability is the ability to recover the unknown parameter from the observation, {\color{blue} when the initial condition is know.}
\begin{definition}[Identifiability]
\bleu{Given an initial state $x_0$ in $\Omega$,}
	system   (\ref{sys2}) is said to be identifiable if for any distinct $\theta_1$, $\theta_2$ in $\Theta$, there exists  $t \geq 0$ such that 
	\[
	h( x(t,x_0, \theta_1) )\neq h(x(t,x_0, \theta_2)).
	\]
\end{definition}

\m

The concepts of observability and identifiability are similar. Consider the augmented system which consists in adding the parameter $\theta$ as part of the (augmented) state vector with a null dynamics:
\begin{equation}\label{sysAugment}
	\left \{
	\begin{array}{l}
		\dot x (t)=  f(x(t), \theta),\\
		\dot \theta = 0,\\
		\\
		y (t) = h(x(t,\theta)). \\
	\end{array}
	\right.
\end{equation}
\bleu{
\begin{proposition}\label{prop:augm}
If the system \eqref{sysAugment} is observable on $\Omega\times\Theta$, then system \eqref{sys2} is observable on $\Omega$ and identifiable for any $x_0 \in \Omega$.
\end{proposition}
\begin{proof}
If system \eqref{sysAugment} is observable then we have the property:
$h( x(t,(x_0, \theta_1)) )=h( x(t,(x_0^\prime, \theta_2) ))$ for all $t\geq 0$ implies $(x_0, \theta_1)=(x_0^\prime, \theta_2)$, which implies $x_0=x_0^\prime$ and  $\theta_1=\theta_2$.
\end{proof}
\n
The converse is not true as it is illustrated by the following example:
\begin{equation}\label{ex-aug}
\left\{\begin{array}{l}
\dot x_1 = x_1 + \theta\, x_2 \\
\dot x_2 = 0 \\
\\
y= x_1
\end{array}
\right.
\end{equation}
System~\eqref{ex-aug} is observable and identifiable.
\begin{equation}\label{ex-aug1}
\left\{\begin{array}{l}
\dot x_1 = x_1 + \theta\, x_2 \\
\dot x_2 = 0 \\
\dot \theta = 0\\
\\
y= x_1
\end{array}
\right.
\end{equation}
System~\eqref{ex-aug1}, whose state is $(x_1,x_2,\theta)$, is not observable since the output is the same for initial conditions $(x_1,x_2,\theta)$ and $(\bar x_1,\bar x_2,\bar \theta)$ satisfying $x_1=\bar x_1$ and $x_2 \theta= \bar x_2 \bar \theta$.
}

\smallskip

At several places in the following, we shall consider this augmented dynamics.

\m

Actually, for an  epidemiological model it is unlikely  to know   the initial condition and it has long been recognized that initial conditions play a role in identifying
the parameters \cite{DiopFliess1991,MR1261705,MR872589,MR1957979,MR2138148}.

\m

What we have called identifiability is also known as \textit{structural identifiability}. This expression has been coined by R. Bellman and K.J. \AA str\"{o}m \cite{Bellmann1970} in 1970. This is to stress that identifiability depends only on the dynamics and the observation, under ideal conditions of noise-free observations and error-free model. This is a mathematical and \textit{a priori} problem \cite{Jacquez-Greif-MBS85}.

\bigskip

\bleu{Sometimes, it happens that the concept of observability is too strong in the sense that a system in the form \eqref{sys} does not satisfies the observability property but it is still possible to know part of the initial condition. A relaxed definition of observability is as follows.

\begin{definition}[Partial observability]\label{def:obs-partial}
Let $\varphi$ be a smooth map from $\Omega$ to $\varphi(\Omega) \subset \R^{n'}$.
We say that the system \eqref{sys} is partially observable on $\Omega$ with respect to $\varphi$ if for two initial states $x_0$, $x_0^\prime$ in $\Omega$ with $\varphi(x_0) \neq \varphi(x_0^\prime)$, there exists a time  $t\geq 0$ such that
	\[
	h(x(t,x_0)) \neq h(x(t,x_0^\prime)) .
	\]
	
	Two states $x_0$, $x_0^\prime$ in $\Omega$ are called indistinguishable with respect to $\varphi$ if one has $\varphi(x_0) \neq \varphi(x_0^\prime)$ and
	\[
	h(x(t,x_0)) = h(x(t,x_0^\prime)), \quad t \geq 0
	\]
Alternatively, we shall say the function $\varphi$ is observable on $\Omega$.
\end{definition}

\bigskip

In a similar way, one can relax the definition of identifiability as follows.

\begin{definition}[Partial Identifiability]\label{def:iden-partial}
Let $\vartheta$ be a smooth map from $\Theta$ to $\vartheta(\Theta) \subset \R^{p'}$. Given  an initial state $x_0$ in $\Omega$, 
	system   (\ref{sys2}) is said to be identifiable with respect to $\vartheta$ if for any $\theta_1$, $\theta_2$ in $\Theta$ with $\vartheta(\theta_1) \neq \vartheta(\theta_2)$, there exists  $t \geq 0$ such that 
	\[
	h( x(t,x_0, \theta_1) )\neq h(x(t,x_0, \theta_2)).
	\]
Alternatively, we shall say the parameter function $\vartheta$ is identifiable.
\end{definition}

\begin{exemple}
Consider the system in $\R^3$
\[
\left\{
\begin{array}{l}
\dot x_a = -\alpha x_b\\
\dot x_b = -\beta x_c\\
\dot x_c =\gamma
\end{array}
\right.
\]
where $\alpha$, $\beta$, $\gamma$ are parameters. Solutions of this system can be made explicit:
\begin{align*}
& \displaystyle x_a(t)=  x_a(0) - \alpha x_b(0) t + \frac{\alpha\beta}{2} x_c(0)  t^2  + \frac{\alpha\beta\gamma}{6}t^3,\\
& \displaystyle x_b(t) = x_b(0) - \beta x_c(0) t - \frac{\beta\gamma}{2} t^2,\\
& \displaystyle x_c(t) = x_c(0) +\gamma t .
\end{align*}
Here, we consider $\Omega={\cal D}=\R^3$ and $\Theta=\{ (\alpha,\beta,\gamma) ; \; \alpha >0 , \; \beta > 0, \; \gamma > 0 \}$.
Let us show that the system is observable for the observation map $h(x_a,x_b,x_c)=x_a$. We consider two initial conditions $(x_a(0),x_b(0),x_c(0))$, $(x_a^\prime(0),x_b^\prime(0),x_c^\prime(0))$ such that their respective solutions verify $x_a(t)=x_a^\prime(t)$ for any $t \geq 0$. From the expression of the solution $x_a(\cdot)$, this implies that one has 
\[
- \alpha x_b(0) t + \frac{\alpha\beta}{2} x_c(0) t^2  = - \alpha x_b^\prime(0) t + \frac{\alpha\beta}{2} x_c^\prime(0)  t^2
\]
for any $t \geq 0$, and for $t>0$ one gets by differentiating with respect to $t$
\[
- \alpha x_b(0) + \frac{\alpha\beta}{2} x_c(0) t = - \alpha x_b^\prime(0) + \frac{\alpha\beta}{2} x_c^\prime(0) t
\]
When $t>0$ tends to $0$, one obtains $x_b(0)=x_b^\prime(0)$ and then $x_c(0)=x_c^\prime(0)$. Therefore, the two initial conditions coincide and there is no indistinguishable distinct initial condition. The system is thus observable.

Consider now the observation map $h(x_a,x_b,x_c)=x_b$. When two initial conditions are such that their solutions verify $x_b(t)=x_b^\prime(t)$, one gets from the expression of the solution $x_b(\cdot)$
\[
-\beta x_c(0)t = -\beta x_c^\prime(0)t
\]
for any $t>0$ and thus $x_c(t)=x_c^\prime(t)$. However, having $x_a(0)\neq x_b^\prime(0)$ with $x_b(t)=x_b^\prime(t)$ and $x_c(t)=x_c^\prime(t)$ gives exactly the same solution $x_b(\cdot)$, $x_c(\cdot)$. The system is not observable but it is partially observable  with respect to the map
\[
\varphi(x_a,x_b,x_c)=(x_b,x_c).
\]

\medskip

Let us now study the identifiability of the system, first for the observation map $h(x_a,x_b,x_c)=x_a$. For a given initial condition with $x_b(0)\neq 0$ and $x_c(0)\neq 0$, if two sets of parameters $(\alpha_1,\beta_1,\gamma_1)$, $(\alpha_2,\beta_2,\gamma_2)$ give the same solution $x_a(\cdot)$, one has
\[
- \alpha_1 x_b(0) t + \frac{\alpha_1\beta_1}{2} x_c(0)  t^2  + \frac{\alpha_1\beta_1\gamma_1}{6}t^3 =
- \alpha_2 x_b(0) t + \frac{\alpha_2\beta_2}{2} x_c(0)  t^2  + \frac{\alpha_2\beta_2\gamma_2}{6}t^3 
\]
for any $t \geq 0$, and for $t>0$ one gets
\[
- \alpha_1 x_b(0)  + \frac{\alpha_1\beta_1}{2} x_c(0)  t  + \frac{\alpha_1\beta_1\gamma_1}{6}t^2 =
- \alpha_2 x_b(0)  + \frac{\alpha_2\beta_2}{2} x_c(0)  t  + \frac{\alpha_2\beta_2\gamma_2}{6}t^2 
\]
When $t>0$ tends to $0$, one obtains $\alpha_1=\alpha_2$ and thus one has
\[
\frac{\alpha_1\beta_1}{2} x_c(0)  t  + \frac{\alpha_1\beta_1\gamma_1}{6}t^2 =
\frac{\alpha_2\beta_2}{2} x_c(0)  t  + \frac{\alpha_2\beta_2\gamma_2}{6}t^2 
\]
for any $t>0$. In a similar way, one can show that one has $\alpha_1\beta_1=\alpha_2\beta_2$ and then $\alpha_1\beta_1\gamma_1=\alpha_2\beta_2\gamma_2$, which proves that one has $\beta_1=\beta_2$ and $\gamma_1=\gamma_2$ (remind that parameters are assumed to be non null). The system is thus identifiable.

When the observation map is $h(x_a,x_b,x_c)=x_b$, having the same output $x_b(\cdot)$ implies
\[
-\beta_1 x_c(0) t -\frac{\beta_1\gamma_1}{2}t^2 = -\beta_2 x_c(0) t -\frac{\beta_2\gamma_2}{2}t^2
\]
for any $t \geq 0$
As before, one obtains $\beta_1=\beta_2$ and $\gamma_1=\gamma_2$ but having $\alpha_1\neq \alpha_2$ with $\beta_1=\beta_2$ and $\gamma_1=\gamma_2$ provides the same output $x_b(\cdot)$. Indeed the dynamics of $x_b$ and $x_c$ are decoupled from the one of $x_a$, where parameter $\alpha$ appears only. The system is not identifiable but is partially identifiable with respect to the map
$\vartheta(\alpha,\beta,\gamma)=(\beta,\gamma)$.
\end{exemple}
}

\section{Historical notes}
The observability concept has been introduced by   Kalman  \cite{MR0152167}  in the sixties for linear systems.
For nonlinear systems, observability has been characterized circa the seventies \cite{MR389284,zbMATH03616280}. The definition is  given by Hermann and Krener in the framework of differential geometry.
Identifiability and structural identifiability  has been introduced in compartmental analysis  in 1970 by Bellman and  \AA str\"{o}m \cite{Bellmann1970} in a paper that appeared in a bio-mathematics journal.
The problem of identifiability is now addressed in text-books  \cite{Ljung99,MR1482525,WalPro,MR710757}. Numerical identifiability of linear control system is implemented in softwares such as {\tt Matlab} and {\tt Scilab}.

\m

Identifiability of nonlinear systems has been addressed in different context and the first systematic approach is by Tunali and Tarn in 1987 
\cite{MR872589} in the differential geometry framework.
The introduction of the concepts of differential  algebra in control theory is due to Fliess  around  1990 \cite{DiopFliess1991,diopfliessECC91,MR958700}  followed by Glad \cite{FliessGlad93,MR1261705}.
Identifiability is a general problem which has received different names depending on the  community:
\begin{itemize}
	\item[-] observation, identification,
	\item[-]  data assimilation,
	\item[-]  inverse problem,
	\item[-]  parameters estimation.
\end{itemize}

\m

``Data assimilation'' is mainly used in meteorology and oceanography
\cite{LeDimetTalagrand86,MR628731}.
A direct (as opposed to inverse) problem is considering a model which, when introducing an input, gives an observed output. The parameters are considered as known. Conversely the ``inverse problem'' is to reconstruct the parameters from the knowledge of the output \cite{MR2130010}. Finally, ``parameters estimation'' is used in the probability and statistics domains  \cite{Akaike74,bolker2008ecological,MR3197254,MR3807914,MR3807916,ONEILL2002103,Roda:2020ab}.

\section{Identifiability in mathematical epidemiology}

\n
Identifiability is well known in bio-mathematics from the seventies, as already mentioned with  the paper of Bellman and  \AA str\"{o}m \cite{Bellmann1970}. However, considering identifiability in mathematical epidemiology is relatively recent 
\cite{MR3553947,MR2726194,MR1957979,MR2817815,MR2785878,Eisenberg:2013aa,MR2142487}. The first paper, to our knowledge, considering identifiability of an intra-host model of HIV is by Xia and Moog \cite{MR1957979}, and has been published in 2003  in a journal of automatic control.

\section{The concept of observers}

The construction of an {\em observer} is based on an estimation approach different from statistical methods: it consists of determining a dynamical system (called an ``observer'') whose input is the vector $y(\cdot)$ of measures acquired over time, and whose state is an estimate $\hat x(t)$ of the (unknown) true state $x(t)$ of the system at time $t$. 

An observer estimates $x(t)$ continuously over time and without anticipation, in the sense that the estimate $\hat x(t)$ is updated at each instant $t$ through its dynamics as measurement $y(t)$ is available, without requiring the knowledge of any future measurement. This is why an observer is sometimes also called a ``software sensor''. Since the estimate $\hat x(t)$ is given by the solution of a system of differential equations, the main idea behind an observer is the paradigm of {\em integrating instead of differentiating} the signal $y(\cdot)$. Note that although an observer is primarily devoted to state estimation, an observer can also aim to reconstructing simultaneously state and parameters, when some parameters are unknown (in this case a parameter vector $p$ is simply considered a part of the system dynamics with $\dot p=0$).

The most well-known observer is the so-called {\em Luenberger observer} \cite{Luenberger71} that is recalled in Chapter \ref{secobservers}, and that has inspired most of the existing observers (several ones are discussed in Chapter \ref{secobservers}).  However, observers are yet relatively unpopular in Mathematical Epidemiology, comparatively to other application domains (such as mechanics, aeronautics, automobile, \textit{etc}). The aim of the present  review is also to promote the development and use of observers for epidemiological models. 

Chapter \ref{secobservers} presents the theoretical background of observers construction and their convergence as estimators based on the model equations, independently of the quality of real data. In a complementary way, Chapter \ref{sectionPractical} discusses some implementation issues when observers are used with real world data that could be corrupted with noise.
\chapter{Mathematical foundations}

\abstract*{We give mathematical characterizations of observability and identifiability.}

\abstract{We give mathematical characterizations \bleu{and properties} of observability and identifiability.}

\bleu{
\section{Preliminaries}

Here and in all the following chapters, we shall consider that the maps $f$ and $g$ that define the system
\begin{equation}\label{sys4}
	\left \{
	\begin{array}{l}
		\dot x =  f(x), \quad x\in \Omega \subset \R^n, \\[2mm]
		y  =h(x) \in Y \subset \R^m
	\end{array}
	\right.
\end{equation}
are analytic at any point $x \in \Omega$.
We shall denote by $x(t,x_0)$ the solution of $\dot x=f(x)$ for the initial condition $x(0)=x_0$.

\medskip

We recall that a function $\varphi: {\cal D} \mapsto \R$ is analytic on a open domain ${\cal D}$ of $\R^n$ 
\bleu{(we will write $\varphi \in C^{\omega}(\cal D,\R)$)}
if it is $C^\infty$ (i.e.~infinitely differentiable) and its Taylor series locally converge, that is for any $x_0$ in ${\cal D}$ there exists a neighborhood ${\cal V}$ of $x_0$ in ${\cal D}$ such that
\[
\varphi(x)=\lim_{n \to +\infty} \sum_{k=0}^n \frac{\varphi^{(k)}(x_0)}{k!}(x-x_0)^k, \quad x \in {\cal V} .
\]
Indeed, up to our knowledge, the great majority of epidemiological model in the literature are  analytic.

\bigskip

The map $f : \Omega \longrightarrow \R^n$ is also called a {\em vector field} on $\Omega$, as its image for each $x \in \Omega$ is a (velocity) vector in $\R^n$. 

Given a $C^\infty$ function $g: \R^n \mapsto \R$, the classical definition of Lie derivative  of $g$ with respect to the vector field $f$  is given by 
\begin{equation}
\label{defLiederivative}
\mathcal L_f (g)(x) = \dfrac{d}{dt} \,  g(x(t,x) ) {\bigg |}_{t=0}=  \langle \nabla g (x)  | f (x)\rangle , \quad x \in \Omega
\end{equation}
where $\nabla g$ is the gradient of $g$ and $\langle \; |  \; \rangle $ the inner product of $\R^n$. For $k>1$, we define by induction
\[
\mathcal L_f^k (g)(x)= \mathcal L_f ( \mathcal L_f^{k-1} (g))(x)
\]
\bleu{Note that $\mathcal L_f^k (h) (x)$ is the value of the $k$-th time derivative at time $0$ of the output $y(t)=h(x(t,x))$ along the solution of system~\eqref{sys4} with initial state $x$, but considered as function of $x$:
\begin{align*}
    & \mathcal L_f (h) (x)=\dot y(x),\\ 
    & \mathcal L_f^2 (h) (x)=\ddot y(x),\\
    & \qquad \vdots\\
    & \mathcal L_f^k (h) (x)=y^{(k)}(x)
\end{align*}
}
}

\bleu{
For a $C^\infty$ vector-valued function $g: \R^n \mapsto \R^q$, the Lie derivative $\mathcal L_f (g)(x)$ at point $x$ is the vector of Lie derivatives of each component $g_j$ ($j=1\cdots q$)
\[
\mathcal L_f (g)(x) = \left[\begin{array}{c}
\mathcal L_f (g_1)(x)\\
\vdots \\
\mathcal L_f (g_q)(x)
\end{array}\right]
\]
When $f$ and $g$ are linear, that is $f(x)=Ax$ and $g(x)=Cx$, then $\mathcal L_f^k (g) (x)$ is equal to $CA^kx$.

\bigskip

We shall denote ${}^\top$ the transposition operator for a vector or a matrix.
}

\section{Observability}

\bleu{In this section, we introduce the {\em observation space} and the concept of {\em local observability} that can be checked with the help of differential calculus in terms of a (local) {\em rank condition} to be checked. 
For linear systems, this rank condition depends only on the matrices defining the system and the output and ensures the observability property as defined in Section \ref{secdef}. For analytic nonlinear systems, we shall see that this rank condition depends on the vector field and the output function as well as on the initial condition $x_0$ and only ensures a "local" observability of the system on a neighborhood of $x_0$.}

\bigskip

\bleu{The observability definition~\ref{def:obs-intro} given in Chapter~\ref{chap:intro} states a global observability property. It also can be formulated as follows. 
\begin{definition}[Observability]\label{def:obs-2}
The analytic system~\eqref{sys4} is observable if for any initial conditions $x_1$ and $x_2$, 
for any $T>0$, one has:
$h(x(t,x_1)) = h(x(t,x_2))$ for all $t\in [0, T]$ implies $x_1=x_2$.
This is equivalent to say that the map:
\[\begin{array}{ll}
\Omega \rightarrow & {\mathcal C}^\omega([0,T],\R^m)\\[2mm]
x_0\mapsto &
\left\{\begin{array}{ccc}
[0,T]&\rightarrow& \R^m \\
t&\mapsto &h(x(t,x_0))
\end{array}\right\}
\end{array}
\]
is injective.
\end{definition}
Observability means that the initial state $x_0$ (and therefore the trajectory starting from this initial state) can be uniquely determined by the knowledge of the data of the output $y(\cdot)$  on any nontrivial time interval.
}


The components of the observation map $h$ are denoted by $h=(h_1, \cdots, h_m)$.
Each $h_i$ is a  $\mathcal C^{\infty}$ function  from the state space  $\R^n$ to  $\R$.

\begin{definition}[\cite{zbMATH03616280}]\label{def:obs-space} 
	The observation space  $\mathcal O$  of  (\ref{sys4}) is the subspace of the vector space  $\mathcal C^{\infty} (\R^n,\R)$ containing  $h_i$ and  invariant by the Lie derivative $\mathcal L_f$. 
\end{definition}
\medskip
The observation space is generated by the different Lie derivatives of the  $h_i$:
\[\bleu{
\begin{array}{ll}
    \mathcal O  &  = \text{span}_{\R}
\left\{ \mathcal L^k_f h_i\; : i=1, \dots,m,\; 
k \in \N \right\}  \\
      &  =\text{span}_{\R}
\left\{ h_i, \,  \mathcal L_f \, h_i , \mathcal L^2_f h_i, \dots\; : i=1, \dots,m \right\}. 
\end{array}
}
\]
\bleu{The observation space $\mathcal O$ contains the observation function (also called the output function) and all derivatives of the
output function along the system trajectories.
For a linear system: $f(x)=Ax$ and $y=h(x)=Cx$, the observation space is generated by the $n$ functions 
\[Cx,\; CAx,\;\ldots\;,\; CA^{n-1}x.
\]
We recall that thanks to Cayley-Hamilton Theorem $A^k$ for $k\geq n$ is a linear combination of $A^q$ with $q\leq n-1$.
}

\n
We have the following result 
\bleu{that relates the observation space to the observability property.}
\begin{theorem}\label{th:obs-space}
	For an analytic system  (i.e., $f$ and  $h$ are analytic functions) the observability is equivalent to the separation of the points of the state space $\R^n$ by $\mathcal O$ i.e., if  $x_1 \neq x_2$ there exists  $g \in \mathcal O$ such that  $g(x_1) \neq g (x_2)$.
\end{theorem}
\begin{proof}
By analyticity we have
\bleu{
\[
y(t, x_0)= h ( x(t, x_0))= \sum_{k\geq 0} \; 
 \left(\dfrac{d^k}{ds^k} h(x(s,x_0)){ \bigg |}_{s=0}\right) \, \dfrac{t^k}{k!} ,
\]
}
but, by induction we have the following relation 
\bleu{\[ \dfrac{d^k}{ds^k} h (x( s, x_0)){ \bigg |}_{s=0}= \mathcal L_f^k h (x_0).\]
}  
Then a necessary and sufficient condition to distinguish  $x_1\neq x_2$ is that there exists  \bleu{$k$}  such that 
\bleu{
\[
\mathcal L_f^k. h \;(x_1) \neq  \mathcal L_f^k. h \;(x_2).
\]
}
\end{proof}
\bleu{
\begin{exemple} [Application to a virus dynamics model]\label{ex-VIH}
We consider a simple model of an HIV-1 infection \cite{1078.92502}: 
\begin{equation}\label{hiv1}
\left \{
		\begin{array}{l} 
{\dot T}=\Lambda-\mu_T\,  T-\beta\,  V\, T, \\ 
{\dot T^*}=\beta\,  V\, T-\delta\,  T^*, \\ 
{\dot V}=r\, \delta\, T^*-c\,  V, \\ 
\\
y  = V,
\end{array}
\right.
\end{equation} 
where $T$, $T^*$, $V$ denote the concentrations of uninfected (healthy) and infected host cells, and freevirions, respectively. 
The rate of infection is given by $\beta VT$ , with $\beta$ being the infection rate constant. The parameters  $\delta$ and $c$ are the removal rates of the infected cells and virus particles respectively. The healthy cells (T) are produced at a rate $\Lambda$ and $\mu_T$ is the death rate per T cell.
It is assumed that on average each productively infected cell produces $r$ (a positive integer) virions during its lifetime, so the per-capita rate of viral production for an infected cell is given by $r\,  \delta$. All parameters are positive.
We assume that the measurement
of the viral load is available.
%
%

\medskip

Let us show that system~\eqref{hiv1} is observable on the set $\{(T,T^*,V) \in \R^3 : T\geq 0,\, T^*\geq 0,\,  V>0 \}$.
%
%
We shall prove that the observation space associated to System~\eqref{hiv1} separates the points of set $\{(T,T^*,V)^\top \in \R^3 : T\geq 0,\, T^*\geq 0,\,  V>0 \}$.

Here we have
$x=(T,T^*,V)^\top$ and
\[
f(x)=\left[\begin{array}{c}
			\Lambda-\mu_T\,  T-\beta\,  V\, T \\
			\beta\,  V\, T-\delta\,  T^* \\
			r\, \delta\, T^*-c\,  V
		\end{array}\right]
  \]
and the measurable output $y=h(x)=V$.
Computing the Lie-derivatives of the output gives

\[
\mathcal L_f h(x) =   \langle \nabla h (x)  | f (x) \rangle
=\langle\;
\left[\begin{array}{c}
			0 \\
			0 \\	 			
			1
		\end{array}\right]
| \left[\begin{array}{c}
			\Lambda-\mu_T\,  T-\beta\,  V\, T \\
			\beta\,  V\, T-\delta\,  T^* \\
			r\, \delta\, T^*-c\,  V
		\end{array}\right]
 \;\rangle
 =r\, \delta\, T^*-c\, V
 \]
and
\begin{align*}
\mathcal L_{f}^{2} h(x) =   
\langle \nabla \mathcal L_f h(x)  | f (x) \rangle
& =\langle\;
\left[\begin{array}{c}
			0 \\
			r\, \delta \\	 			
			-c
		\end{array}\right]
| \left[\begin{array}{c}
			\Lambda-\mu_T\,  T-\beta\,  V\, T \\
			\beta\,  V\, T-\delta\,  T^* \\
			r\, \delta\, T^*-c\,  V
		\end{array}\right]
 \;\rangle\\
& =r\, \delta\, \beta\, T\,V  -r\, \delta  \left(c +\delta \right) T^* +c^{2}\,V .
\end{align*}

\n
The explicit expression of $\mathcal L_{f}^{3} h(x)$ is quite long but it can be easily  shown that  $\mathcal L_{f}^{3} h(x)$ can be expressed as a function of $h(x)$, $\mathcal L_{f} h(x)$, $\mathcal L_{f}^{2} h(x)$ and the parameters of the system. So to study observability of 
System~\eqref{hiv1}, one does not need to compute Lie-derivative of order higher than 2.

The observation space $\mathcal O$ contains the functions 
\begin{align*}
 & g_1(x)=V,\\
 & g_2(x)=r\, \delta\, T^*-c\,V,\\
& g_3(x)=r\, \delta\, \beta\,T\,V  -r\,\delta  \left(c +\delta \right) T^* +c^{2}\,V .
\end{align*}
It is easy to show 
\[
\Big\{g_1(x)=g_1(\bar x), g_2(x)=g_2(\bar x), g_3(x)=
g_3(\bar x)\Big\}
\Longrightarrow V=\bar V, \, T^*=\bar T^*
\]
and 
$T= \bar T \text{ if } V\neq 0$.
Hence,
these functions separate the points of 
$\R^3_{+}\setminus \{(T,T^*,0)\}$ but they do not separate points of the form $(T,T^*,0)$ and $(\bar T,T^*,0)$ with $\bar T \neq T$.  
Thus System~\eqref{hiv1} is observable on the set $\{(T,T^*,V) \in \R^3 : T\geq 0,\, T^*\geq 0,\,  V>0 \}$.
\end{exemple}

\bigskip

Applying Theorem~\ref{th:obs-space} to linear systems 
\begin{equation}\label{syslin}\tag{$\Sigma_L$}
\left\{
	\begin{array}{l}
		\dot x =  A x, \quad x\in \R^n, \\[2mm]
		y  =C x, \quad y\in  \R^m,
	\end{array}
\right.
\end{equation}
allows to obtain a simple algebraic necessary and sufficient condition for observability.
\begin{proposition}\label{prop-lin}
The linear system~\eqref{syslin} is observable if and only if 
the {\em observability matrix}
	\begin{equation}
		\label{obsmat}
		O_{(C,A)}=\left[\begin{array}{c}
			C\\
			CA\\\
			\vdots\\\\
			CA^{n-1}
		\end{array}\right]
	\end{equation}
is of full rank, i.e., rank $O_{(C,A)}=n$. 
\end{proposition}
Indeed, using Theorem~\ref{th:obs-space}, we have the following successive equivalences:
\[
\begin{array}{c}
\text{The linear system~\eqref{syslin} is not observable} 
\\
\big\Updownarrow
\\
\text{Its observation space }
\mathcal O_{\Sigma_L} \text{ does not separate the points of the state space } \R^n
\\
\big\Updownarrow
\\
\exists x_1 \neq x_2 : \forall g \in \mathcal O_{\Sigma_L}, g(x_1)=g(x_2)
\\
\big\Updownarrow
\\
\exists x_1 \neq x_2 : \forall i\in \{0,1,\ldots n-1\}, 
C A^i x_1=C A^ix_2
\\
\big\Updownarrow
\\
\exists x_1 \neq x_2 :  
x_1 - x_2 \in Ker\,C A^i$, $\forall i\in \{0,1,\ldots n-1\}
\\
\big\Updownarrow
\\
 x_{1}-x_{2} \in ker\;C\ \cap\ ker\;CA\ \cap\ ker\;CA^2\ \cap \ldots \cap 
ker\;CA^{n-1}
\\
\big\Updownarrow
\\
ker\, O_{(C,A)} \neq \{0\}
\\
\big\Updownarrow
\\
rank\, O_{(C,A)} <n.
\end{array}
\]
%
%

%

\bigskip

The result of Theorem~\ref{th:obs-space} can be also reformulated as follows.
For the sake of writing simplicity, we consider real scalar output. Extension to vector output is straightforward.
Let us define, for $k\in \N$, the map: $x\mapsto H_k(x)=
(h(x),\mathcal L_f h (x),\ldots, \mathcal L_f^k h (x))^\top$. Then we have the following characterization of observability :
\begin{proposition}\cite{Inouye_1977}
\label{prop:H_k}
	Suppose System~\eqref{sys4} is analytic. Then it is observable if and only if $H_k(x_1)=H_k(x_2)$, for all $k\in \N$, imply that $x_1=x_2$.
\end{proposition}
This is equivalent to say that the following real analytic mapping (from the state space to an infinite dimensional space)
\[x\mapsto H_{\infty}(x)=
(h(x),\mathcal L_f h (x),\ldots, \mathcal L_f^k h (x), \ldots)^\top
\]
 is injective.
%

\begin{remark}
In general there is no value of $k$ to stop. This is illustrated by the following example \cite{Inouye_1977}: 
\begin{equation}\label{exo}
\left\{\begin{array}{l}
\dot x = -x, \quad x\in \R \\
\\
\displaystyle y=h(x)= x\,\prod_{i=1}^\infty (1-e^{x^2-i^2})^i
\end{array}\right.
\end{equation}
It has been proved in \cite{Inouye_1977} that
\begin{enumerate}
\item[-] $h$ is analytic on $\R$,
\item[-] system \eqref{exo} is observable,
\item[-] $ \forall k\in \N^*$, the equation $H_k(x)=H_k(0)$ has countably infinite solutions $x=0$ and $x=\pm(k+i),\,i=0, 1, 2, ...$ .
\end{enumerate}
\end{remark}
Observability can also be checked using the following proposition.
\begin{proposition}\label{prop-expression}
If any state $x \in \Omega$  can be expressed as a function of the observation $y$ and its time derivatives $y^{(k)}=\mathcal L_f^k h (x)$, that is there exists a map $\phi$ such that 
$x=\phi\Big(h(x),,\mathcal L_f h (x),\ldots, \mathcal L_f^k h (x), \ldots\Big)$, then System \eqref{sys4} is observable on $\Omega$.
\end{proposition}
\begin{proof}	
Let $x_1 \in \Omega$ and $x_2 \in \Omega$ be such that $H_{\infty}(x_1)=H_{\infty}(x_2)$. Then, applying $\phi$, we obtain $x_1=x_2$ which proves that map $H_{\infty}$ is injective on $\Omega$. 
\end{proof}	
\begin{exemple}
	Consider the SIR model of  Kermack-McKendrick \cite{KmK1927}	
\begin{equation} \label{KmcK} \hfill 
		\left \{
		\begin{array}{rl}
			\dot S &=  -\beta \, \dfrac{S}{N} \, I, \\[3mm]
			\dot I  &=\beta \,  \dfrac{S}{N}\,  I - \gamma  \, I, \\[3mm]
			\dot R  &= \gamma  \,I
		\end{array}
		\right.
\end{equation}
for which  the parameters $\beta, \gamma$ 
and $N$ are assumed to be known, where the total population $N=S+I+R$ is constant since $\dot N=\dot S + \dot I +\dot R=0$.
We assume that the recovery $y=\gamma I$ is observed.
Here we have
\[
x= \left(\begin{array}{c}
S\\ I \\ R
\end{array}\right) , \quad 
f(x)=\left(\begin{array}{c}
			  -\beta \, \dfrac{S}{N} \, I \\[2mm]
			\beta \,  \dfrac{S}{N}\,  I - \gamma  \, I \\[2mm]			 			\gamma  \,I
		\end{array}\right), \quad 
h(x)=\gamma I
\]
Clearly, the solutions of the system \eqref{KmcK} evolve in the positively invariant set  (which makes biological sense) 
	
	\[\Omega=\{ (S,I,R) \;  |  \; S > 0, \; I > 0, \; R>0,\; S+I+R=N\} \] 
on which one has $y \neq 0$. Then we have the relations (we can divide by $y$) 
	\[ 
 S = \dfrac{N}{\beta}\, \dfrac{\dot y +\gamma\, y}{y} ,  \; \; 
 I= \frac{y}{\gamma} , \; \; R=N-\dfrac{N}{\beta}\, \dfrac{\dot y +\gamma\, y}{y}-\frac{y}{\gamma}
 \]
Therefore, the state vector $x$ 
can be expressed as a function of $y$ and $\dot y$. The system is thus  observable in $\Omega$.
This model is more thoroughly studied in Chapter \ref{chapSIR}.
\end{exemple}
The observability (equivalent) conditions given by Theorem~\ref{th:obs-space} and Proposition~\ref{prop:H_k}  are hard to test, since they involve an infinite number of analytic equations and therefore checking the condition that $\mathcal O$ separates points can be quite a formidable task.
This is one of the reasons for which  different concepts of local observability have been introduced. One of them is the local weak observability \cite{zbMATH03616280}.
As we will see, the advantage of the local weak observability as compared to other types of observability is that it is easy to indicate for it simple sufficient conditions of an algebraic nature.}

\bigskip
%

\begin{definition}\label{def:obs-local}
$ $
\begin{itemize}
	\item[-] The system  (\ref{sys4}) $\Sigma$ is said locally observable if, for any  $x_0$, for any open set $U$ containing $x_0$, $x_0$ is distinguishable from all the points of $U$ for the restricted system 
$\Sigma | U $.\\
\item[-] The system (\ref{sys4}) is weakly observable at $x$ if there exists an open neighborhood $U$ of $x$ such that the only point in  $U$ which is indistinguishable from $x$ is $x$ itself.
The system (\ref{sys4}) is weakly observable if it is weakly observable at 
every $x\in M$.\\
\item[-] The system (\ref{sys4}) $\Sigma$ is locally weakly observable if for any $x_0$ there exists an open set $U$ containing $x_0$, such that for any neighborhood $V$ with  $x_0\in V \subset U$,  $x_0$ is distinguishable for  $\Sigma | V$ from all the points of  $V$.
\end{itemize}
\end{definition}
Intuitively,  a system   is locally weakly observable if one
can instantaneously distinguish each point from its neighbors.
The local weak observability can be characterized as follows.
\begin{definition}\cite{SontagCDC91}\label{def:dO}
Let \bleu{$\mathcal O$ be the observation space of System~\eqref{sys4}, we define}
	\[d\mathcal O= \{ d\psi  \mid  \psi \in \mathcal O \}\]
	where $d\psi(x)$ is the differential of $\psi$ at $x$.
\end{definition}
\begin{definition}\label{def:rank}
	A system $\Sigma$ is said to satisfy  the observability rank condition \bleu{(ORC)} at  $x$ if the dimension of  $d\,\mathcal O $ at  $x$ satisfies   
\[\text{dim} \left (d \, \mathcal O  (x)\right )=n \]
\end{definition}
where $d\mathcal O$ is generated by the gradients of the  
$\mathcal L_f^kh$.
\begin{theorem}[Hermann-Krener \cite{zbMATH03616280}]\label{th:rank} 
\bleu{If the {\bf analytic} system~\eqref{sys4} satisfies the observability rank condition (ORC)
	at $x_0$  then $\Sigma$  is locally weakly observable at $x_0$.} 
\end{theorem}
\begin{proof}
\bleu{Since  $\text{dim} \left (d \, \mathcal O  (x_0)\right )=n $, there exists  $n$ functions  $\varphi_1, \cdots, \varphi_n \in \mathcal O$  such that the gradients  $d\varphi_1(x_0), \cdots, d\varphi_n(x_0) $ are linearly independent.
	Therefore the  function $\Phi: x \mapsto (\varphi_1(x), \cdots, \varphi_n(x))$ has a non-singular  Jacobian  in  $x_0$. 
	As a consequence, from the Inverse Function Theorem,  there exists an open set  $U$ containing  $x_0$ where 
\bleu{$\Phi$} 
is a bijection.

On any open set $V \subset U$ suppose that we have 
$h(x(t,x_0))=h(x(t,x_1))$ for $t \in [0,T]$. Then, from the fact that $f$ and $h$ are analytic and 
\[
\dfrac{d^k}{dt^k} h (x( t, x_0)){ \bigg |}_{t=0}= \mathcal L_f^k h (x_0),
\]
we have 
\[h(x(t,x_0))-h(x(t,x_1))
=\sum_{k\geq 0} \; \dfrac{t^k}{k!} \, \left(\mathcal L_{f}^{k}h(x_0) 
-\mathcal L_{f}^{k}h(x_1)\right)=0, \text{ for }t \in [0,T].
\]
This implies 
$\mathcal L_{f}^{k}h(x_0) -\mathcal L_{f}^{k}h(x_1)=0$ for all $k\geq 0$ 
which implies that $\varphi_i(x_0) =\varphi_i(x_1)$ for all $i=1,\ldots,n$ since each $\varphi_i$ is a linear combination of the $\mathcal L_{f}^{k}h$. Therefore $\Phi(x_0)=\Phi(x_1)$ and hence $x_1=x_0$ since $\Phi$ is a bijection. This proves that 
$x_0$ is distinguishable from all points of $V$ and hence the analytic 
system~\eqref{sys4} is locally weakly observable at $x_0$.
} 
\end{proof}

\bleu{
If the observability rank condition is satisfied everywhere the system is locally weakly observable.

\medskip
A converse result has been proved in \cite{zbMATH03616280}:
\begin{proposition}\cite{zbMATH03616280,Casti82}\label{prop:equiv0}
If the system is is locally weakly observable then the rank condition is satisfied almost everywhere, i.e., in an open dense subset of the state space.
\end{proposition}
To summarize, we have for analytic systems:
\begin{proposition}\cite{zbMATH03616280,Casti82}\label{prop:equiv}
For the {\bf analytic} system~\eqref{sys4}, the relationships between the various observability concepts are given in the following diagram
%
%
\[
\begin{array}{ccc}
& & \text{ORC satisfied almost everywhere} \\
& & \big\Updownarrow  \\
\text{\eqref{sys4} locally observable} & \Longrightarrow & \text{\eqref{sys4} locally weakly observable}\\
 \big\Downarrow& & \big\Updownarrow  \\
\text{\eqref{sys4} observable} &\Longrightarrow & 
\text{\eqref{sys4} weakly observable}
\end{array}
\]
\end{proposition}
\begin{remark}
For linear systems, the five properties are equivalent.
\end{remark}
\begin{remark}
The analytic system
$\dot x=0$, $y=x^3$ is observable on $\R$ and weakly locally observable but the ORC is not satisfied at $x=0$. However it is satisfied on $\R \setminus \{0\}$ and so it is satisfied almost everywhere. 
\end{remark}
\begin{remark}
If an analytic system is observable then the observabilty rank condition (ORC) is satisfied almost everywhere. It must be noticed that the converse is not true:  the analytic system
\[\left\{\begin{array}{l}
\dot x=1, \; x\in \R \\
\\
 y=(\sin x, \cos x) \in \R^2
\end{array}\right.
\]
satisfies the ORC at any $x\in \R$ but it is not observable because the states $x$ and $x+ 2 k \pi$ are indistinguishable. 
\end{remark}
It must be emphasized that the study of the observability of analytic systems requires   either to deal with an infinite number of analytic equations or  to compute the dimension of the linear space generated by the gradients of all Lie derivatives of the output. In general, there is no bound on the number of  Lie derivatives necessary to conclude, as seen in example \eqref{exo}.
}

\begin{remark} For linear systems $\dot x=Ax$, $y=Cx$, all the definitions of observability are equivalent to having 
	the {\em observability matrix}
	\begin{equation*}
		O_{(C,A)}=\left[\begin{array}{c}
			C\\
			CA\\\
			\vdots\\\\
			CA^{n-1}
		\end{array}\right]
	\end{equation*}
of full rank (see for instance \cite{Kailath}).
\end{remark}

\medskip

The observability analysis can also be a way to choose the right sensor, as illustrated on the following example.
\begin{exemple}
	Consider a population model structured in five age classes, whose population sizes are
	\begin{itemize}
		\item[] $x_1$ for juveniles,
		\item[] $x_2$ for subadults capable of reproduction when adults,
		\item[] $x_3$ for subadults not capable of reproduction when adults,
		\item[] $x_4$ for adults capable of reproduction,
		\item[] $x_5$ for adults not capable of reproduction,	    
	\end{itemize}
	and the dynamics is
	\begin{equation*}
		\left\{\begin{array}{l}
			\dot x_1 = -\alpha x_1 + \beta x_4\\
			\dot x_2 = \frac{\alpha}{2}x_1 - \alpha x_2 -m_1x_2\\
			\dot x_3 = \frac{\alpha}{2}x_1 - \alpha x_3 -m_1x_3\\
			\dot x_4 = \alpha x_2 -m_2 x_4\\
			\dot x_5 = \alpha x_3 -m_2 x_5
		\end{array}\right.
	\end{equation*}
	(where $\alpha$ is an aging rate, $m_1$, $m_2$ are mortality rates, and $\beta$ is a fecundity rate). If only one sub-population $x_i$ can be targeted for measurement, one can easily check that the only possibility for the system to be observable is to measure the variable $y=x_5$.
 \bleu{
This system is of the form $\dot x = A x$, $y=Cx$ with 
\[
A=\left[\begin{array}{cc}
\overbrace{\left[\begin{array}{cccc}
-\alpha &\hspace{3mm} 0 & 0 & \beta  \\
\alpha/2 &\hspace{3mm} -\alpha-m_1 & 0& 0  \\
\alpha/2 &\hspace{3mm} 0 & -\alpha-m_1 & 0  \\
0 &\hspace{3mm} \alpha & 0& -m_2
\end{array}\right]}^{A_4} & \begin{array}{c}
0\\
0\\
0\\
0
\end{array}\\[9mm]
\begin{array}{llll}
0 \hspace{11mm}& 0\hspace{11mm} & \alpha \hspace{9mm} & 0\end{array}
 & -m_2
\end{array}\right]
\]
that is of the form 
\[
\dot x = \left[\begin{array}{c}
 A_4\, 
\left[\begin{array}{c}
x_1 \\
x_2 \\
x_3 \\
x_4 
\end{array}\right]\\
\alpha x_3 -m_2 x_5
\end{array}\right], \quad y=Cx
\]

From this, we see that for any integer $k\geq 1$, $A^k$ is of the form 
\[
\left[\begin{array}{cc}
\begin{array}{c}
A_4^k
\end{array} & \begin{array}{c}
0\\
0\\
0\\
0
\end{array}\\[9mm]
\begin{array}{llll}
* \hspace{11mm}& *\hspace{11mm} & * \hspace{9mm} & *\end{array}
 & (-m_2)^k
\end{array}\right].
\]
Therefore, if $C=[ c_1\; c_2 \; c_3\; c_4\; 0]$ then 
$C A^k =( *\;  *\; *\;  *\; 0)$ and hence the observability matrix will be of rank $\leq 4$ and so the system is not observable if the output does not depend on $x_5$.

Now, if $C=[ 0\; 0 \; 0\; 0\; 1 ]$ then the observability matrix defined by~\eqref{obsmat} is given by:
\[
O_{(C,A)}=
\left[\begin{array}{ccccc}
0 & 0 & 0 & 0 & 1 
\\
 0 & 0 & \alpha  & 0 & -m_{2}  
\\
 \frac{\alpha^{2}}{2} & 0 & -\alpha  \left(m_1 +m_{2} +\alpha \right) & 0 & m_{2}^{2} 
\\
 -\frac{\alpha^{2} \left(2 \alpha +m_1 +m_{2} \right)}{2} & 0 & h_{43} & \frac{\beta  \,\alpha^{2}}{2} & -m_{2}^{3} 
\\
h_{51} & \frac{\alpha^{3} \beta}{2} & 
h_{53} & 
-\frac{\beta  \,\alpha^{2} \left(2 \alpha +m_1 +2 m_{2} \right)}{2} & m_{2}^{4} 
\end{array}\right]
\]
with
\[\begin{array}{l}
h_{43}=\left(\alpha^{2}+\left(2 m_1 +m_{2} \right) \alpha +m_1^{2}+m_1 m_{2} +m_{2}^{2}\right) \alpha, \\[2mm]
h_{51}=\left(\frac{3 \alpha^{2}}{2}+\frac{\left(3 m_1 +2 m_{2} \right) \alpha}{2}+\frac{m_1^{2}}{2}+\frac{m_1 m_{2}}{2}+\frac{m_{2}^{2}}{2}\right) \alpha^{2}, \\[2mm]
h_{53}=-\alpha  \left(m_1 +m_{2} +\alpha \right) \left(\alpha^{2}+2 \alpha  m_1 +m_1^{2}+m_{2}^{2}\right).
\end{array} 
\]
One has
$det \; O_{(C,A)} =-\dfrac{\alpha^8 \beta^2}{8} \neq 0$. Hence 
$O_{(C,A)}$ is of full rank which proves the observability of the system when $y=x_5$.
}
\end{exemple}

\section{\bleu{About identifiability}}

Since very often the initial conditions are not known, or partially known, we will consider in the following the problem of \bleu{joined} identifiability and observability, considering the augmented system (\ref{sysAugment}). Note that identifiability-only problems are a special case in which $y=x$ \bleu{(this is why we consider here the more general case of joined identifiability and observability).}

\medskip

\color{blue}
Consider a parametrized system
\begin{equation}\label{sys6}
 	\left \{
 	\begin{array}{l}
 		\dot x (t)=  f(x(t), \theta), \; x(0)=x_0\\
 		\\
 		y (t) =h(x(t), \theta)\\
 	\end{array}
 	\right.
 \end{equation}
 with $x \in \R^n$, $y \in \R^m$ and $\theta \in \Theta \subset \R^p$ and assume that system \eqref{sys6} is observable for any known value of the parameter $\theta \in \Theta$. 
 Let us denote the observability map parameterized by $\theta$
 \[
H_k(x,\theta)=\left[ 
\begin{array}{c}
h(x,\theta)\\
{\cal L}_f(h)(x,\theta)\\
\vdots\\
{\cal L}_f^{k-1}(h)(x,\theta)
\end{array}
\right]
 \]
 When there exists an integer $k>1$ and a map $\Phi_\theta: \R^n \mapsto H_k(\R^n,\theta)$ parameterized by $\theta$, such that the output differential equation
 \[
 y^{(k)}(t) = \Phi_\theta(y(t),\dot y(t), \cdots, y^{(k-1)}(t)), \quad t \geq 0
 \]
as an unique solution in $H_k(\R^n,\theta)$, for any $\theta \in \Theta$, one may study the dependency of the map $\Phi_\theta$ with respect to $\theta$ to ensure the identifiability of the system. In particular, when $f$ and $h$ are polynomials with coefficients parameterized by $\theta$, the map $\Phi_\theta$ is a polynomial with coefficients $c_i(\theta)$ where $c: \Theta \mapsto \R^\nu$ for some $\nu$. The injectivity of the map $c$ is clearly a necessary condition to have identifiability, but this is not sufficient as one can see on the following example.
\color{black}

\medskip

\color{blue}
\begin{exemple}
	
	\n
	Consider the system in $\R_+^2$
	
	\begin{equation}\left \{
		\begin{array}{rl}
			\dot x_1&=  -\theta_1 \, x_1+\theta_2\, x_2\\
			\dot x_2&=-\theta_2\,x_2\\
			\\
			y  &=x_1. \\
		\end{array}
		\right.\,.
	\end{equation}
	where $\theta \in \Theta=\R \times \R^\star$ is the unknown vector of parameters.
	\n
	One has
\[
H_1(x,\theta)= 
\left[ \begin{array}{c} x_1\\ -\theta_1 \, x_1+\theta_2\, x_2 \end{array} \right] 
\]
which is invertible w.r.t.~$x$ on $\R_+^2$ for any $\theta \in \Theta$
\[
H_1^{-1}(z,\theta) =
\left[ \begin{array}{c} z_1\\[2mm] \displaystyle \frac{z_2+\theta_1z_1}{\theta_2} \end{array} \right]
\]
The system is thus observable for any $\theta \in \Theta$.
Moreover, one has
\[
\ddot y  = -\theta_1\, \dot y- \theta_1^2\, x_2 = -\theta_2\, \dot y-\theta_1 \, (\dot y +\theta_2 \, y)
\]
that is
\[
\ddot y=\Phi_\theta(y,y'):=-(\theta_1+\theta_2) \, \dot y -\theta_1\,\theta_2 \, y
\]
	Clearly the application $\theta \mapsto (\theta_1+\theta_2, \theta_1\,\theta_2)$ is injective on $\Theta$, but the system is not identifiable if $x_2(0)=0$: the solution verifies $x_2(t)=0$ for any $t>0$ and $x_1$ is solution of $\dot x_1=-\theta_1 x_1$ independently of the value of $\theta_2$.
	
		

	
\end{exemple}
\color{black}

\medskip

\color{blue}
A natural (and usual) way to ensure the joined observability-identifiability property is to consider the augmented state:
\[
\tilde x = \left[\begin{array}{c}
x\\
\theta\end{array}\right] \in \R^n \times \Theta
\]
and require the observability of the extended dynamics
\begin{equation*}\left \{
		\begin{array}{l}
			\displaystyle \frac{d}{dt} \tilde x(t)  = \tilde f(\tilde x(t)):=f(x(t),\theta), \;\;
			\tilde x(0) = (x(0),\theta)\\
			\\
			y(t)  =\tilde h(\tilde x(t))):=h(x(t),\theta) \\
		\end{array}
		\right.
	\end{equation*}
 and we fall back on a problem of pure observability.
\color{black}

\bigskip

\color{blue}
\begin{remark}
    A theoretical answer to the problem of parameters reconstruction when measuring the whole state $x$ in $\R^n$ has been given by D.~Aeyels \cite{MR670047,MR626654} and E. Sontag \cite{MR1938330} in different form. For an (analytic) system with $r$ parameters, it is {\em generically} sufficient to choose  $2\,r+1$ measures at different times to distinguish two different states (the term {\em generically} means here that for any system excepted for a non dense subset of systems among all the analytic systems in $\R^n$).

\end{remark}
\color{black}

\n

\section{Identifiability does not {\color{blue} necessarily} require observability}

\bleu{Let us stress that identifiability does not necessarily imply observability. It can happen that the knowledge of the output function $y(\cdot)$ allows to reconstruct uniquely the set of parameters, but not necessarily the state variables of the system. We give below an example of such a situation.}


\begin{exemple}
The following "academic" model is identifiable but not observable.
	\begin{equation}\label{ex-idenonobs}
		\left \{
		\begin{array}{l}
			\dot x_1 =  -\alpha \, (x_1+x_2)   \\
			\dot x_2  =\alpha \,  (x_1-x_2)\\
			\\
			y = \dfrac{1}{2} \,(x_1^2+x_2^2) 
		\end{array}
		\right.
	\end{equation}
 \bleu{One immediately gets 
 \[
 \dot y= -\alpha \, y .
 \]
 For (unknown) initial conditions $(x_1(0),x_2(0)\neq (0,0)$, one has $y(0)>0$ and $y$ is thus a positive function. Then one obtains $\alpha = -\dot y /y$: the system is identifiable on $\R^2 \setminus \{0\}$. 
 
 Compute now further derivatives: $\ddot y=\alpha^2\, y$, $\cdots$ , $y^{(p)} = (-1)^p \, \alpha^p\, y$. 
	Formally, one gets
	\[\text{Jac } [ h, \mathcal L_f\,h,  \mathcal L^2_f\,h, \cdots , \mathcal L^P_f\,h] = 
	\begin{bmatrix}
		x_1 & -\alpha x_1 & \alpha^2 x_1 & \cdots & (-1)^p \alpha^p x_1\\
		x_2 & -\alpha x_2 & \alpha^2 x_2 & \cdots & (-1)^p \alpha^p x_2\\
		0   &  -y         & 2\alpha y    & \cdots & \alpha (-1)^p \alpha^{p-1}\, y
	\end{bmatrix}, \]
	which is of rank $2$ for any $(x_1,x_2) \in \R^2\setminus\{0\}$ and for any positive integer $p$. The parameter $\alpha$ is identifiable, but the system is not observable. Consider another solution $(\xi_1(\cdot),\xi_2(\cdot))$ for the initial condition $\xi_1(0)=-x_1(0)$, $\xi_2(0)=-x_2(0)$. One can straightforwardly check that this solution verifies  $\xi_1(t)=-x_1(t)$, $\xi_2(t)=-x_2(t)$ with the same output $y(t)$ for any $t \geq 0$. Therefore, these two solutions cannot be distinguished, which shows that the system is not observable.}

\bleu{Another way to show the non-observability is to remark that 
 we have 
 \[
 y(t)=e^{-\alpha t} y(0) 
=\dfrac{1}{2} e^{-\alpha t} \left(x_1(0)^2+x_2(0)^2\right)
=\dfrac{1}{2} e^{-\alpha t} \left(\bar x_1(0)^2+\bar x_2(0)^2\right)
\]
for any  $\left(\bar x_1(0),\bar x_2(0)\right) \in S_0$ where $S_0$ is the circle centered at the origin with radius $r_0=\sqrt{x_1(0)^2+x_2(0)^2}$.
Thus the output cannot distinguish the solutions emanating from different points of this circle. Hence System~\eqref{ex-idenonobs} is non observable.

This can also be proved by remarking that one has
\[
\mathcal L^P_f\,h(x)=(-1)^p \, \alpha^p h(x)=
\dfrac{(-1)^p \, \alpha^p}{2} \,(x_1^2+x_2^2) .
\] 
This means that the associated observation space is generated by the function $x_1^2+x_2^2$ that does not separate the points of $\R^2$ and hence, according to Theorem~\eqref{th:obs-space}, system~\eqref{ex-idenonobs} is non observable.
}
\end{exemple}

\n

\section{Identifiability via \bleu{decoupled variables}}
\label{secidentchgvar}

\bleu{Sometimes it is difficult to prove identifiability of a model in the original set of coordinates. Let us show the interest of considering other variables that possess "good" properties.}
We consider a system in $\R^n$ parameterized by $\theta \in \Theta \subset \R^p$ of the following form
\begin{equation}
	\label{sys_param}
	\left\{\begin{array}{lll}
		\dot x & = & f(x,\theta), \quad x(0)=x_0 \in X\\
  \\
		y & = & h(x)
	\end{array}\right.
\end{equation}
where $X \subset \R^n$ is positively invariant for any $\theta \in \Theta$. \m

\begin{proposition}\cite{bichara-ifac23}
	\label{prop_ident_chtvar}
	Assume that the following properties hold.
	\begin{enumerate}
		\item The map $f$ verifies
		\begin{equation}
			\label{hypo_ident_chtvar}
			\forall \theta_1, \theta_2 \in \Theta, \; \forall x \in X, \quad
			\theta_1 \neq \theta_2 \Longrightarrow f(x,\theta_1)\neq f(x,\theta_2) .
		\end{equation}
		\item There exist smooth maps $g$, $\tilde g$ and $l$ such that
		\begin{enumerate}
			\item for any solution $x(\cdot)$ of \eqref{sys_param} in $X$, $w(t):=g(x(t),y(t)) \in W \subset \R^m$ verifies
			\begin{equation}
				\label{dyn_w}
				\dot w(t) = l(w(t),y(t)), \quad t \geq 0
			\end{equation}
			\item for any $x \in X$, one has 
   \begin{equation}
       \label{hypo_inv_chtvar}
   w=g(x,h(x))) \Longleftrightarrow x=\tilde g(w,h(x))
   \end{equation}
		\end{enumerate}
	\end{enumerate}
	Then the system \eqref{sys} is identifiable over $\Theta$ for any initial condition in $X$.
\end{proposition}

\medskip

\begin{proof}
\bleu{Fix $x_0 \in X$ and denote by $x_{\theta}(\cdot)$ the solution of \eqref{sys_param} for the parameter $\theta \in \Theta$. Consider $\theta_1$, $\theta_2$ in $\Theta$ that give the same output function:  $h(x_{\theta_1}(t))=h(x_{\theta_2}(t))=y(t)$ for any $t\geq 0$.

Let $w(\cdot)$ be the (unique) solution of the Cauchy problem
 \[
\dot w = l(w,y(t)), \quad w(0)=g(x_0,h(x_0)).
 \]
Then, one has 
 \[
w(t)=g(x_{\theta_i}(t),y(t)) \quad t \geq 0, \quad i=1,2, 
 \]
and from property \eqref{hypo_inv_chtvar}, one obtains
\[
	\tilde g(w(t),y(t))=x_{\theta_i}(t) \quad t \geq 0, \quad i=1,2, 
\]
that is $x_{\theta_1}(t)=x_{\theta_2}(t)$ for any $t> 0$. Therefore, for any $t>0$ one has also 
$\dot x_{\theta_1}(t)=\dot x_{\theta_2}(t)$ or $f(x(t),\theta_1)=f(x(t),\theta_2)$ where $x(t)=x_{\theta_1}(t)=x_{\theta_2}(t)$.
Finally, from condition \eqref{hypo_ident_chtvar}, we deduce that one has necessarily $\theta_1=\theta_2$, which shows the identifiability of the system.}
\end{proof}

\medskip

\bleu{This result states that when a system is identifiable when measuring the whole state $x$ (condition \eqref{hypo_ident_chtvar}), and there exits a variable $w$ whose dynamics is {\em decoupled} in the sense that it depends on $w$ and $y$ only (condition \eqref{dyn_w}) such that the map $(x,y)\mapsto w$ is invertible with respect to $x$ (condition \eqref{hypo_inv_chtvar}), then the system is identifiable when measuring $y$ only. Up to our knowledge, this approach has  not been deployed in the literature.}
Let us illustrate this result on an intra-host model for malaria infection \cite{MR3181991}.

\begin{exemple}[Malaria model]
	\label{example_malaria}
	The state vector is $x=(S,\ I_1,\ldots\ ,I_5\ ,M)^\top$  in $\R_+^7$, where $S$ is the concentration of uninfected erythrocytes in the blood, $I_i$ are the concentrations of infected erythrocytes in different age classes, and $M$ is the concentration of free merozoites. The dynamics is given by the following system
	\begin{equation} \label{AMGreduit}
		\left\{\begin{array}{l}
			\dot{S}=\Lambda -\mu_S\, S-\beta\, SM,\\[2mm]
			\dot{I}_1=\beta\, SM - (\gamma_1+\mu_1)\,I_1,\\[2mm]
			\dot{I}_2=\gamma_1\, I_1-(\gamma_2+\mu_2)\,I_2,\\[2mm]
			\vdots\\[2mm]
			\dot{I}_5=\gamma_{4}\, I_{4}-(\gamma_5+\mu_5)\,I_5,\\[2mm]
			\dot{M}=r\,\gamma_5\,I_5-\mu_M\, M-\beta\, SM ,
		\end{array}\right.
	\end{equation}
	where the different parameters are
	\begin{itemize}
		\item[]  $\Lambda$: recruitment of the healthy red blood cells (RBC).
		\item[] $\beta$: rate of infection of RBC by merozoites.
		\item[] $\mu_S$: natural death rate of healthy cells.
		\item[] $\mu_i$: natural death rate of $i$-th stage of infected cells.
		\item[] $\gamma_i$: transition rate from $i$-th stage to $(i+1)$-th stage of infected cells.
		\item[] $r$ : number of merozoites released by the late stage of infected cells.
		\item[] $\mu_M$ : natural death rate of merozoites.
	\end{itemize}
	The two first stages  of infected erythrocytes ($I_1$ and $I_2$) correspond to the concentration of free circulating parasitized erythrocytes than can be observed (seen on peripheral blood smears). \bleu{Typically, the quantity 
 \[
 y(t)=h(x(t))=I_1(t)+I_2(t)
 \]
 is measured at any time $t$. Among parameters in (\ref{AMGreduit}), most of them ($\mu_i,\ \gamma_i$, and $r$) are known or at least widely accepted by the community, but
	the infection rate $\beta$, which is playing a crucial role, is unknown and cannot be estimated by biological considerations. Let us then write the dynamics \eqref{AMGreduit} as $\dot x = f(x,\beta)$.} It takes the form
	\begin{equation}\label{AMGx}
		\left\{\begin{array}{ll}
			\dot{x}=& f(x,\beta):= A\, x+\beta\,  S M\ E +\Lambda\, e_1, \\
   \\
			y=&C\, x 
		\end{array}
		\right.
	\end{equation}
	with
	\[
	A=\left[\begin{array}{ccccccc}
		-\mu_S & 0 & 0 & 0 & 0 & 0 & 0\\
		0 & -\gamma_1-\mu_1 & 0 & 0 & 0 & 0 & 0\\
		0 & \gamma_1 & -\gamma_2-\mu_2 & 0 & 0 & 0 & 0 \\
		0 & 0 & \gamma_2 & -\gamma_3-\mu_3 & 0 & 0 & 0 \\
		0 & 0 & 0 & \gamma_3 & -\gamma_4-\mu_4 & 0 & 0 \\
		0 & 0 & 0 & 0 & \gamma_4 & -\gamma_5-\mu_5 & 0 \\
		0 & 0 & 0 & 0 & 0 & r\gamma_5 & -\mu_M
	\end{array}\right],
	\]
	\[
	E=\left[\begin{array}{r}
		-1\\ 1\\ 0\\ 0\\ 0\\ 0\\ -1
	\end{array}\right], \quad
	e_1=\left[\begin{array}{r}
		1\\ 0\\ 0\\ 0\\ 0\\ 0\\ 0
	\end{array}\right], \quad 
	C=\left[\begin{array}{ccccccc}
		0 & 1 & 1 & 0 & 0 & 0 & 0
	\end{array}\right]
	\]
 \bleu{Due to the dimension of the dynamics, it is not easy to check the identifiability of he parameter $\beta$.
However, on the domain $X=(\R_+\setminus\{0\})^7$, one has $SM\neq 0$, which implies the property
\[
f_1(x,\beta_1)=f_1(x,\beta_2) \Rightarrow \beta_1=\beta_2, \quad x \in X .
\]
For the parameter $\theta=\beta$, condition \eqref{hypo_ident_chtvar} of Proposition \ref{prop_ident_chtvar} is thus fulfilled.
 Note that one has $ECE=E$. Therefore one can consider the variable
 \[
 w=g(x,y):=x-Ey=(I-EC)x
 \]
 whose dynamics is independent of the non-linear term $\beta S M$:
	\[
	\dot w = l(w,y):=\bar Aw+\bar AEy+\Lambda\, e_1
	\]
	where we posit $\bar A= A-ECA$. Given $w$ and $y$, the state $x$ is then given by 
 \[
 x=\tilde g(w,y):=w+Ey .
 \]
 Conditions \eqref{dyn_w} and \eqref{hypo_inv_chtvar} of Proposition \ref{prop_ident_chtvar} are thus also satisfied, which allows to conclude without any other calculation that the parameter $\beta$ is identifiable.}
	
\end{exemple}

This example illustrates the \bleu{possible} interest of exploiting conjointly identifiability and observability to solve the identifiability problem.

\chapter{Analysis of the Kermack-McKendrick model}
\label{chapSIR}

\abstract*{This chapter is dedicated to the study of this well-known SIR model.}

\abstract{This chapter is dedicated to the study of the well-known SIR model.}

\section{History}

\n
The SIR model of  Kermack and McKendrick \cite{KmK1927} is certainly  one of the most famous model in Epidemiology. It is given and studied in all of  the  classic books of Mathematical Epidemiology. This model appears in the book of Bailey, which is probably the first book in Mathematical Epidemiology. Some examples can be found in  \cite{AndMay91,MR2002k:92001,MR3969982,Brauer1945,DalGan99,MR3752193,MR3409181,Murr2002}. The figure, in the original paper, fitting the model to plague data in Bombay during the 1906 year, is one of the most famous pictures in Epidemiology. A research with  SIR in  {\tt MathScinet} returns  $11\,106$ articles. 

\m
\n
In the quoted books the SIR model is fitted to data in the following ways:

\begin{itemize}
	\item[-] in  \cite{MR2002k:92001,MR3969982,Brauer1945} the model is fitted to the plague in  Eyam (in the year 1666);
	\item[-]  in \cite{DalGan99} the model is fitted to an influenza epidemic in England and Wales;
	\item[-]  in \cite{MR3752193} a fitting is done with simulated noisy data;
	\item[-] in  \cite{MR3409181,MR2002k:92001}, in a chapter devoted to fitting epidemiological models to data, a SIR model is fitted to an influenza outbreak in an English boarding school. 
\end{itemize}

\m
\n
More recently two publications  \cite{MR2886018,MR3881860} revisit the fit of the Kermack-McKendrick SIR model to the plague in Bombay.

\m
\n
As already mentioned, before attempting to adjust parameters, an  identifiability analysis should be performed.

\section{The different forms of the SIR model}

\n
\n
The original model  \cite{KmK1927}  is 
\begin{equation}\label{KmcKoriginal}
	\left \{
	\begin{array}{rl}
		\dot S &=  -\tilde{\beta} \, S \, I \\
		\dot I  &=\tilde{\beta} \,S\,  I - \gamma  \, I \\
		\dot R &= \gamma\,I  
	\end{array}
	\right.
\end{equation}
where $S,I,R$  represent respectively the numbers of susceptible, infectious and removed individuals. 

\n
This model can also be found in a slightly different form
\begin{equation}
	\left \{
	\begin{array}{rl} \label{KmcKN}
		\dot S &=  -\beta \, \dfrac{S}{N} \, I \\[2mm]
		\dot I  &=\beta \,  \dfrac{S}{N}\,  I - \gamma  \, I \\[2mm]
		\dot R &= \gamma \,I  
	\end{array}
	\right.
\end{equation}
where  $N=S+I+R$ is the total population. Obviously, one can pass from one model to the other though $\tilde{\beta} =\beta/N$. Both models are mathematically equivalent as long as $N$ is a constant. However, we stress that identifying $\tilde\beta$ does not allow one to estimate the parameters $\beta$ and $N$ separately. For instance, estimating $\tilde\beta$ and $\gamma$ only (without knowing $N$ or $\beta$) does not allow one to estimate the basic reproduction number:
\[ 
\mathcal R_0= \dfrac{\tilde{\beta}\, N}{\gamma}= \dfrac{\beta }{\gamma} \,.
\] 



\section{Observability and identifiability of the SIR model} Quite surprisingly, the observability and identifiability of the original Kermack-Mckendrick SIR model has not been studied much, although this model is commonly used to model epidemics. 

Interestingly, the observability and identifiability of the SIR model with births and deaths, constant population, and an observation $y=k\,I$, has been first studied in 2005 \cite{MR2142487}:
\begin{equation}
	\left \{
	\begin{array}{rl}
		\dot S &= \mu \, N -\beta \, \dfrac{S}{N} \, I -\mu\, S \\
		\dot I  &=\beta \,  \dfrac{S}{N}\,  I - (\gamma +\mu)  \, I \\[2mm]
		\dot R &= \gamma \,I   - \mu\, R
	\end{array}
	\right.
\end{equation}
where $\mu$ is the renewal rate of the population. The article \cite{MR2142487} concludes that the system is neither observable nor identifiable.

\m
\n
In \cite{MR3784444} the identifiability of     (\ref{KmcKN}) is addressed assuming (i) that the initial conditions (and therefore $N=S(0)+I(0)+R(0)$) are known, (ii) observing $y=kI$ with $k=1$, and using only the input-output relation to conclude. Under assumptions (i) and (ii), the identifiability is quite immediate, as we shall see, but of limited interest.

\color{blue}
\subsection{The SIR model when observing a ratio of the infected population}
Here we study the observability-identifiablity property of the SIR model
\color{black}
\begin{equation}\label{KmcK}
	\left \{
	\begin{array}{rl}
		\dot S &=  -\beta \, \dfrac{S}{N} \, I \\
		\dot I  &=\beta \,  \dfrac{S}{N}\,  I - \gamma  \, I \\
  \\
		y  &= k \,I  
	\end{array}
	\right.
\end{equation}
The {\color{blue} dynamics of $R$} has been omitted since $R=N-S-I$ . The observation is $y=k\,I$, in other words only a fraction of the infectious individuals are observed. This situation is used for example in \cite{MR3881860,Roda:2020aa}. 
\bleu{In general, the values of $N$ and $k$ are also not known. We have the following result.}
\bleu{\begin{theorem}\label{th:sir}
System~(\ref{KmcK}) is neither observable, nor identifiable. 
\end{theorem}
\begin{remark}\label{rk:sir0b}
Theorem~\ref{th:sir} can be obtained from \cite{MR2142487} by setting $\mu=0$. However we provide a short and elementary proof.
\end{remark}
}

\bleu{
\begin{proof}
System \eqref{KmcK} is obviously not observable on the invariant set $I=0$. Therefore we study the observability-identifiability properties of System \eqref{KmcK} on the following positively invariant open set  
	\[ \mathcal D =\{ (S,I)  \mid  S > 0, \; I > 0, \; S+I <N \}. \]
We will show that there exist a couple of distinct initial condition and distinct parameters that generate the same output.
The computation of successive time derivatives of the output $y$ gives:
\begin{equation}
\label{deriv-y}
\left\{\begin{array}{l}
y = k I, \\[3mm]
\dot y = y \left(\dfrac{\beta  S}{N}-\gamma \right) ,   \\[3mm]
\ddot y = -\dfrac{\beta^{2} S}{k N^{2}} y^2 \,+\dfrac{\dot y^{2}}{y} ,
\\[3mm]
\dddot y = \left(-y \ddot y +\dot y^{2}\right)\dfrac{ \beta}{k N}-\dfrac{3 \dot y^{3}}{y^{2}}+\dfrac{4 \dot y \ddot y}{y}. 
\end{array}
\right.
\end{equation}
\color{blue}
One can observe that the third derivative of $y$ is expressed as a function of lower derivatives of $y$ of and parameter $\frac{\beta}{kN}$ only. This property remains true for any further derivative of $y$. This implies that if one considers two sets of initial condition and parameters such that the expressions of $y$, $\dot y$, $\ddot y$ and $\frac{\beta}{kN}$ coincide, then any further derivative also coincides. By analyticity of the solutions of the system 
\bleu{as well as of the corresponding outputs}, we deduce that \bleu{their outputs} are the same for any time $t>0$.
\color{black}
More precisely, take two different initial conditions $(S_0,I_0) \ne 
(\bar S_0,\bar I_0)$ and sets of parameters 
$(N,k,\gamma,\beta) \neq (\bar N,\bar k,\bar \gamma,\bar \beta)$ such that
\begin{align*}
    & k I_0 = \bar k \bar I_0 , \\
    & \dfrac{\beta}{k N}=\dfrac{\bar \beta}{\bar k \bar N}, \\
    & \dfrac{\beta  S_0}{N}-\gamma = \dfrac{\bar \beta  \bar S_0}{\bar N}-\bar \gamma, \\
    & \dfrac{\beta^{2} S_0}{k N^{2}} = \dfrac{\bar \beta^{2} \bar S_0}{\bar k \bar N^{2}}
\end{align*}
then $x_0=(S_0,I_0,N,k,\gamma,\beta)$ and 
$\bar x_0=(\bar S_0,\bar I_0,\bar N,\bar k,\bar \gamma,\bar \beta)$ generate the same output function $y(\cdot)$.
\end{proof}
Therefore, system~(\ref{KmcK}) is not observable nor identifiable. However, it is {\it partially} observable (respectively {\it partially} identifiable) in the sense of Definitions~\ref{def:obs-partial} and \ref{def:iden-partial} that is some functions of the state variables (respectively of the parameters) are observable (respectively identifiable). This is made precise in the following proposition.
\begin{proposition}\label{prop_obs-partialSIR} System \eqref{KmcK}
 satisfies the following properties.
 \begin{itemize}
    \item[i.] The state functions $kI$, $\dfrac{\beta S}{N}$ and $kS$ are observable.\\[2mm]
    \item[ii.] The parameters $\gamma$ and $\dfrac{\beta}{k \, N}$ are identifiable.
\end{itemize}
\end{proposition}
\begin{proof}
From \eqref{deriv-y}, and since one has $-y \ddot y +\dot y^{2} = 
\dfrac{\beta^{2} k^{2} I^{3} S}{N^{2}} >0$, we get
\begin{equation*}
\left\{\begin{array}{l}
k I= y, \\[3mm]
\dfrac{\beta  S}{N}-\gamma = \dfrac{\dot y}{y},   \\[3mm]
\dfrac{\beta^{2} S}{k \,N^{2}} = 
\dfrac{\dot y^{2}}{y^{3}}-\dfrac{\ddot y}{y^{2}}, \\[3mm]
\dfrac{\beta}{k N} = 
\dfrac{\dddot y -\dfrac{4 \dot y \ddot y}{y}+\dfrac{3 \dot y^{3}}{y^{2}}}{-y \ddot y +\dot y^{2}}. 
\end{array}
\right.
\end{equation*}
From these expressions,  we deduce
\begin{align*}
 & \dfrac{\beta S}{N}=\dfrac{\left(y \ddot y -\dot y^{2}\right)^{2}}{y \left(\dddot y \,y^{2}-4 \dot y \ddot y y +3 \dot y^{3}\right)} , \\
& \gamma = \dfrac{\left(y \ddot y -\dot y^{2}\right)^{2}}{y \left(\dddot y \,y^{2}-4 \dot y \ddot y y +3 \dot y^{3}\right)}
+\dfrac{\dot y}{y} , \\
& k S = \dfrac{\dfrac{\beta S}{N}}{\dfrac{\beta}{k N}}
=\dfrac{\left(-y \ddot y +\dot y^{2}\right)^{3} y}{\left(y^{2} \dddot y -4 \dot y \ddot y y +3 \dot y^{3}\right)^{2}}.
\end{align*}
Thus the parameters $\gamma$ and $\dfrac{\beta}{k \, N}$ are identifiable and
the state functions $kI$, $\dfrac{\beta S}{N}$ and $kS$ are observable.
\end{proof}
\medskip
\begin{remark}
    Since the total population $N$ (constant for the model considered)  is often considered known in epidemiology or epidemic modeling (e.g. from the Census Bureau), Proposition~\ref{prop_obs-partialSIR} implies that 
    the quantity $\beta/k$ is identifiable. However $\beta$ and $k$ are not identifiable independently which means that there are infinitely many combinations of $\beta$ and $k$ values for which the model produce the same observable output.
\end{remark}

As a consequence of Proposition~\ref{prop_obs-partialSIR}, we have also the following properties.
\begin{corollary}\label{corol_obs-partialSIR}
For system \eqref{KmcK}, one has 
 \begin{itemize}
    \item[i.] If  $k=\gamma$, then the state variables  $S$ and $I$ are observable, and the parameters  $\gamma, \dfrac{\beta}{ N}$ are identifiable.\\[2mm]
\item[ii.] If  $N$ is known and  if $k=1$ or  $k=\gamma$ then the system is identifiable and observable.
\end{itemize}
\end{corollary}
\medskip
\begin{remark}\label{rk:sir}
	One could believe at the first look that if $k=\gamma$ but with $N$ unknown,  then  (\ref{KmcKoriginal}) is observable, but this is wrong. Certainly $S$ and $I$ are observable, but  $R$ is not observable. Indeed the output 
$y=k I$ and state $I$, solution of system~\eqref{KmcKoriginal}, do not depend on the variable $R$ since $R$ has no influence on the two first equations of \eqref{KmcKoriginal}. So, the output generated by system~\eqref{KmcKoriginal} will be the same for any initial conditions $(S_0,I_0,R_0)$ and $(S_0,I_0,\bar R_0)$ even if 
$R_0 \neq \bar R_0$.
Therefore the value of $N=S+I+R$ is inaccessible. As a consequence, $\mathcal R_0=\tilde\beta N/\gamma$ is not identifiable.
\end{remark}
\begin{remark}\label{rk:sir1}
The second point of Corollary~\ref{corol_obs-partialSIR} can be also proved by considering the "augmented" system (defined in Chapter 1 by Equation~\eqref{sysAugment})  and using Proposition~\ref{prop:augm}. 
Indeed, when $N$ is known and $k=\gamma$, consider the "augmented" SIR model  as follows
\begin{equation}\label{KmcK-augm1}
	\left \{
	\begin{array}{rl}
		\dot S &=  -\beta \, \dfrac{S}{N} \, I \\
		\dot I  &=\beta \,  \dfrac{S}{N}\,  I - \gamma  \, I \\
		\dot \beta &= 0 \\
		\dot \gamma &= 0 \\  
\\
y&= \gamma I
	\end{array}
	\right.
\end{equation}
Posit $x=\left[
\begin{array}{c}
	 S \\
		 I \\
		 \beta \\ 
		 \gamma 
\end{array}\right]$, 
$f(x)=\left[
\begin{array}{c}
	 -\beta \, \dfrac{S}{N} \, I \\
		\beta \,  \dfrac{S}{N}\,  I - \gamma  \, I 
		\\
		 0 \\ 
		 0 
\end{array}\right]$
and $h(x)= \gamma I$.
We have 
\begin{align*}
  \mathcal L_f h(x)  & =   \langle \nabla h (x)  | f (x) \rangle
= \langle\,
\left[
\begin{array}{c}
	 0 \\[3mm]
		\gamma \\[3mm]
		 0 \\  
		 I 
\end{array}\right]
|\left[
\begin{array}{c}
	 -\beta \, \dfrac{S}{N} \, I \\
		\beta \,  \dfrac{S}{N}\,  I - \gamma  \, I \\
		 0 \\
		 0 
\end{array}\right]
\,\rangle\\
 & =\gamma\,\left(\beta \,  \dfrac{S}{N}\,  I - \gamma  \, I\right)
=\gamma\,\left(\beta \,  \dfrac{S}{N} - \gamma \right)\,  I
=\left(\beta \,  \dfrac{S}{N} - \gamma \right) \, h(x)   .
\end{align*}

\begin{align*}
\mathcal L_f^{2} h(x)  & =   \langle \nabla \mathcal L_f h(x)  | f (x) \rangle
= \langle\,
\left[
\begin{array}{c}
	  \gamma\, \dfrac{\beta}{N} \, I \\
		\gamma\,\left(\beta \,  \dfrac{S}{N} - \gamma \right) \\
		 \gamma\, \dfrac{S}{N} \,I\\
		 \beta\,\dfrac{S}{N}\,  I 
\end{array}\right]
|\left[
\begin{array}{c}
	 -\beta \, \dfrac{S}{N} \, I \\
		\beta \,  \dfrac{S}{N}\,  I - \gamma  \, I \\
		 0 \\ 
		 0 
\end{array}\right]
\,\rangle\\
& = -\gamma\, \dfrac{\beta^2}{N^2} \,S I^2
+\gamma\,\left(\beta \,  \dfrac{S}{N} - \gamma \right)^2\,I
=-\gamma\, \dfrac{\beta^2}{N^2} \,S I^2 
+ \left(\beta \,  \dfrac{S}{N} - \gamma \right)\,\mathcal L_f h(x)\\
& =- \dfrac{\beta^2}{N^2} \,S I \, h(x) 
+ \left(\beta \,  \dfrac{S}{N} - \gamma \right)\,\mathcal L_f h(x).
\end{align*}
and computing the third  Lie-derivative of $h$ gives the expression
\[\mathcal L_f^{3} h(x)=
-\dfrac{\beta^2}{N^2} \,S I
\left(\dfrac{\beta S}{N} - \gamma-\dfrac{\beta I}{N}\right)
h(x)
-2 \dfrac{\beta^2}{N^2} \,S I\,\mathcal L_f h(x)
+\left(\dfrac{\beta S}{N} - \gamma\right)\,\mathcal L_f^2 h(x)
\]
\[
\mathcal L_f^{3} h(x)
=\dfrac{\left(S \left(I^{2}-4 I S +S^{2}\right) \beta^{3}+\left(4 I -3 S \right) S \gamma  N \,\beta^{2}+3 N^{2} S \beta  \,\gamma^{2}-N^{3} \gamma^{3}\right)  \gamma\,I}{N^{3}}.
\]
We show that the map 
\begin{equation}\label{H3sir}
 x\mapsto H_{3}(x)=
(h(x),\mathcal L_f h (x),\mathcal L_f^2 h (x), \mathcal L_f^3 h (x))^\top  
\end{equation}
is injective.
Let $x$ and $\bar x$ two elements of the state space
$\Omega=\{(S,I,\beta,\gamma) \in \R^4 : S>0, I>0, \beta>0, \gamma>0\}$.
Suppose that $H_{3}(x)=H_{3}(\bar x)$. Then, we have the following successive implications:

${h(x)=h(\bar x),\mathcal L_f(h)(x)=\mathcal L_f(h)(\bar x),\mathcal L_f^2(h)(x)=\mathcal L_f^2(h)(\bar x)  }
\Longrightarrow$
\[
\left\{
\begin{array}{l}
\gamma I = \bar\gamma \bar I \\
\left(\beta \,  \dfrac{S}{N} - \gamma \right) \, h(x)
=\left(\bar \beta \,  \dfrac{\bar S}{N} - \bar \gamma \right) \, h(\bar x) \\
- \dfrac{\beta^2}{N^2} \,S I \, h(x) 
+ \left(\beta \,  \dfrac{S}{N} - \gamma \right)\,\mathcal L_f h(x)
=- \dfrac{\bar \beta^2}{N^2} \,\bar S \bar I \, h(\bar x) 
+ \left(\bar \beta \,  \dfrac{\bar S}{N} - \bar \gamma \right)\,\mathcal L_f h(\bar x)
\end{array}
\right.
\]
implies
\[
\left\{
\begin{array}{l}
\gamma I = \bar\gamma \bar I \\
\left(\beta \,  \dfrac{S}{N} - \gamma \right) 
=\left(\bar \beta \,  \dfrac{\bar S}{N} - \bar \gamma \right)  \\
\dfrac{\beta^2}{N^2} \,S I  
=\dfrac{\bar \beta^2}{N^2} \,\bar S \bar I  
\end{array}
\right.
\Longrightarrow \; \left\{
\begin{array}{l}
\gamma I = \bar\gamma \bar I \\
\left(\beta \,  \dfrac{S}{N} - \gamma \right) 
=\left(\bar \beta \,  \dfrac{\bar S}{N} - \bar \gamma \right)  \\
\beta^2 \,S I  
=\bar \beta^2 \,\bar S \bar I  
\end{array}
\right.
\]
\[
\Longrightarrow \; \left\{
\begin{array}{l}
\gamma I = \bar\gamma \bar I \\
\beta \,  \dfrac{S}{N} - \gamma 
=\bar \beta \,  \dfrac{\bar S}{N} - \bar \gamma   \\
\beta^2 \,S   
=\bar \beta^2 \,\dfrac{\gamma}{\bar \gamma}\,\bar S   
\end{array}
\right.
\Longrightarrow \;
\left\{
\begin{array}{l}
\gamma I = \bar\gamma \bar I \\
\beta \,  \dfrac{S}{N} - \gamma 
= \dfrac{\beta^2 \bar \gamma}
{\bar \beta \gamma} \dfrac{ S}{N} - \bar \gamma   \\
\beta^2 \,S   
=\bar \beta^2 \,\dfrac{\gamma}{\bar \gamma}\,\bar S   
\end{array}
\right.
\]
Reporting in $\mathcal L_f^3(h)(x)=\mathcal L_f^3(h)(\bar x)$, we obtain
\begin{align*}
& -\dfrac{\beta^2}{N^2} \,S I
\left(\dfrac{\beta S}{N} - \gamma-\dfrac{\beta I}{N}\right)
h(x)
-2 \dfrac{\beta^2}{N^2} \,S I\,L_f h(x)
+\left(\dfrac{\beta S}{N} - \gamma\right)\,L_f^2 h(x)\\
= & -\dfrac{\beta^2}{N^2} \,S I
\left(\dfrac{\beta S}{N} - \gamma-
\dfrac{\bar \beta \bar I}{N}\right)
h(x)
-2 \dfrac{\beta^2}{N^2} \,S I\,L_f h(x)
+\left(\dfrac{\beta S}{N} - \gamma\right)\,L_f^2 h(x)
\end{align*}
This implies
\[
\dfrac{\beta I}{N}
=
\dfrac{\bar \beta \bar I}{N}
\mbox{ and }
\beta I
=
\bar \beta \bar I
\]
With $\gamma I = \bar\gamma \bar I$, we obtain 
$\bar \beta= \dfrac{\bar \gamma}{\gamma} \beta$.
Reporting in 
$\beta \,  \dfrac{S}{N} - \gamma 
= \dfrac{\beta^2 \bar \gamma}
{\bar \beta \gamma} \dfrac{ S}{N} - \bar \gamma$, we obtain
$\bar \gamma =\gamma$. We then deduce $\bar I = I$, $\bar \beta = \beta$, and $\bar S = S$, and thus $x=\bar x$. We have shown that the map $H_{3}$ \eqref{H3sir} is injective which proves that the augmented system~\eqref{KmcK-augm1}  is  observable which implies, thanks to Proposition~\ref{prop:augm}, that System~\eqref{KmcK} is observable and identifiable when the total (constant) population $N$ is known and $k=\gamma$.
 \end{remark}
}
\bleu{
\subsection{The SIR model when observing the  incidence}

Quite often, observations  of new cases per unit time or incidence are available.
We study how this changes the observability and identifiability of the SIR model.

\n
We thus consider the system where the observation is given by 
\begin{equation}\label{out-inci}
y =k \beta\,\dfrac{S \,I}{N}. 
\end{equation}
The system under consideration is then
\begin{equation}\label{KmcK-inci0}
	\left \{
	\begin{array}{rl}
		\dot S &=  -\beta \, \dfrac{S}{N} \, I ,\\[2mm]
		\dot I  &=\beta \,  \dfrac{S}{N}\,  I - \gamma  \, I, \\
		\\
		y  &=  k \beta\,\dfrac{S \,I}{N}.  
	\end{array}
	\right.
\end{equation}
This problem has been addressed for the SIR model with demography for constant population in \cite{MR2142487}. Identifiability with known initial conditions for (\ref{KmcK-inci0}) is also considered in \cite{MR3784444} using input-output relations.
%
\begin{theorem}\label{thm:Incidence}
	The system (\ref{KmcK-inci0}) with incidence observation is neither observable, nor identifiable.
\end{theorem}
\begin{proof}
We proceed as we did  for the proof of Theorem~\ref{th:sir}.
The computation of successive time derivatives of the output $y$ gives:
\begin{equation}
\label{deriv-y-incid00}
\left\{\begin{array}{lll}
%
y & = & k \beta\,\dfrac{S \,I}{N} ,   \\[3mm]
\dot y & = & \left(\dfrac{\beta  S}{N}-\dfrac{\beta  I}{N}-\gamma \right) y ,
\\[3mm]
\ddot y & = & \dfrac{\dot y^{2}}{y}+
 \left(-\dfrac{2 \beta}{k N}\,  y
+\dfrac{\gamma  \beta  I}{N}\right)\,y,
\\[3mm]
y^{(3)} & = & {\color{blue} \dfrac{2\dot y \ddot y}{y}-\dfrac{\dot y^{3}}{y^2}
-4\dfrac{\beta}{kN}y\dot y+\dfrac{\gamma\beta I}{N}\dot y 
+ \dfrac{\gamma\beta}{N}y\left(\dfrac{y}{k}-\gamma I\right)}\\[3mm]
& = & {\color{blue} \dfrac{2\dot y \ddot y}{y}-\dfrac{\dot y^{3}}{y^2}
+\dfrac{\beta}{kN}\left(-4y \dot y +\gamma y ^2\right)+\dfrac{\gamma\beta I}{N}(\dot y-\gamma y)}
\end{array}
\right.
\end{equation}
From the expression of $\ddot y$, one takes
\[
\dfrac{\gamma\beta I}{N}=\dfrac{\ddot y}{y}-\dfrac{\dot y^2}{y^2}
+2\dfrac{\beta}{kN}y
\]
and reporting in the expression of $y^{(3)}$ gives
\[
\begin{array}{lll}
y^{(3)} & = & {\color{blue} \dfrac{2\dot y \ddot y}{y}-\dfrac{\dot y^{3}}{y^2}
+\dfrac{\beta}{kN}\left(-4y \dot y +\gamma y ^2\right)
+(\dfrac{\ddot y}{y}-\dfrac{\dot y^2}{y^2}
+2\dfrac{\beta}{kN}y)(\dot y-\gamma y)} \\
& = & \dfrac{2\dot y \ddot y}{y}-\dfrac{\dot y^{3}}{y^2}
+\left(\dfrac{\ddot y}{y}-\dfrac{\dot y^2}{y^2}\right)(\dot y-\gamma y)
+\dfrac{\beta}{kN}\left(-4y \dot y +\gamma y ^2 +2y(\dot y-\gamma y)\right)\\[3mm]
& = & \dfrac{2\dot y \ddot y}{y}-\dfrac{\dot y^{3}}{y^2}
+\left(\dfrac{\ddot y}{y}-\dfrac{\dot y^2}{y^2}\right)(\dot y-\gamma y)
+\dfrac{\beta}{kN}
\left(-2y \dot y -\gamma y^2 \right)\\[3mm]
& = & \dfrac{\ddot y \dot y}{y}
 +\left(\dfrac{2 \dot y}{y}-\gamma \right)
 \dfrac{ \left(y \ddot y -\dot y^{2}\right)}{y}
-\dfrac{\beta}{k N} \left(\gamma y +2\dot y\right) y.
\end{array}
\]
To summarize, the successive derivatives are
\begin{equation}
\label{deriv-y-incid0}
\left\{\begin{array}{l}
y =k \beta\,\dfrac{S \,I}{N} ,   \\[3mm]
\dot y = \left(\dfrac{\beta  S}{N}-\dfrac{\beta  I}{N}-\gamma \right) y ,
\\[3mm]
\ddot y = \dfrac{\dot y^{2}}{y}+
 \left(-\dfrac{2 \beta}{k N}\,  y
+\dfrac{\gamma  \beta  I}{N}\right)\,y,
\\[3mm]
y^{(3)} = \dfrac{\ddot y \dot y}{y}
+\left(\dfrac{2 \dot y}{y}-\gamma \right)
\dfrac{ \left(y \ddot y -\dot y^{2}\right)}{y}
-\dfrac{\beta}{k N} \left(\gamma  y +2 \dot y \right) y,\\[3mm]
y^{(3)} = \dfrac{3y \dot y \ddot y -2 \dot y^3 }{y^2}
+
\dfrac{\dot y^{2}-y \ddot y }{y}\, \gamma
-\dfrac{\beta}{k N} \left(\gamma  y +2 \dot y \right) y
\end{array}
\right.
\end{equation}
%
This shows that all the derivatives of $y$ of order $q\geq 3$ can be expressed as functions of derivatives $y^{(r)}$, $r=0,\cdots, q-1$ and parameters $\gamma$, $\dfrac{\beta}{k N}$, only. We deduce that for two different $x=(S_0,I_0,N,k,\gamma,\beta)$ and 
$\bar x=(\bar S_0,\bar I_0,\bar N,\bar k,\bar \gamma,\bar \beta)$ such that
%
\[
\gamma =\bar \gamma 
\quad \mbox{and} \quad
\dfrac{\beta}{k N}=\dfrac{\bar \beta}{\bar k \bar N}
\]
and derivatives $y^{(r)}(0)$, $r=0\cdots 2$, coincide, then their output $y(t)$ is the same for any $t>0$ (by analyticity of the system). Notice that one can write
\[
y = k\beta\frac{SI}{N}=\frac{\beta I}{N} \frac{\beta S}{N} \left(\frac{\beta}{kN}\right)^{-1}
\]
Therefore, from expressions \eqref{deriv-y-incid0}, $y^{(r)}(0)$, $r\in\{0,1,2\}$ coincide under the single condition
\[
\dfrac{ \beta  I_0}{N}=
\dfrac{ \bar \beta  \bar I_0}{\bar N}, \quad \mbox{and} \quad 
\dfrac{\beta  S_0}{N}
=\dfrac{\bar \beta  \bar S_0}{\bar N}.
\]
We conclude that the system is neither observable, nor identifiable.
\end{proof}
\begin{proposition}\label{prop-incidence}
If the parameters $N$, $\gamma$ and $k$ are known than System~\eqref{KmcK-inci0} is observable and identifiable.    
\end{proposition}
\begin{proof}
Thanks to Proposition~\ref{prop:augm}, to prove that System~\eqref{KmcK-inci0} is observable and identifiable, it is sufficient to prove that the following (augmented) system whose state is $(S,I,\beta)$: 
\begin{equation}\label{KmcKincid-augm1}
	\left \{
	\begin{array}{rl}
		\dot S &=  -\beta \, \dfrac{S}{N} \, I \\[2mm]
		\dot I  &=\beta \,  \dfrac{S}{N}\,  I - \gamma  \, I \\[2mm]
		\dot \beta &= 0  \\  
\\
y&= k\beta\dfrac{SI}{N}
	\end{array}
	\right.
\end{equation}
is observable.
 From the expression of $y^{(3)}$, we obtain
 \[
\beta 
 = \dfrac{\gamma  (\dot y^2- y \ddot y)y 
+(3 \dot y \ddot y
 - y y^{(3)})y  - 2 \dot y^3}
{\gamma  y^4+2 \dot y  y^3 }\, k N.
\]
Reporting in $\ddot y$, we obtain $I$ as a function of $N$, $\gamma$, $k$, and the successive derivatives of $y$. The state variable $S$ is obtained by reporting the expressions of $\beta$ and $I$ in the expression of $y$.
Thus System~\eqref{KmcKincid-augm1} satisfies the condition of Proposition~\ref{prop-expression} and therefore it is observable on 
$\{(S,I,\beta) \in \R^3 : S>0, I>0, S+I<N, \beta >0\}$.
\end{proof}

\bigskip

\begin{remark}
If it is rather the cumulative incidence that is observed, that is
\begin{equation*}\label{cumul}
 	y(t)= k \;\int_0^t  \, 
 	\beta\,\dfrac{S(\tau) \,I(\tau)}{N} \; d\tau. 
 \end{equation*} 
the system is also not observable or identifiable, and the same results as in Theorem \ref{prop-incidence} are available. Indeed all the $k$-th  derivatives of $y$ (for $k \geq 1$) correspond to the $k-1$-th derivatives of the observation \eqref{out-inci}.
\end{remark}
}

\chapter{Observers synthesis}
\label{secobservers}

\abstract*{We review different techniques to build observers.}

\abstract{We \bleu{expose the concept of observers} and review different techniques to build observers.}

\section{Introduction}

\bleu{As in the previous chapters we consider an observed system 
\begin{equation}
\label{sys_sect_obs}
	\left \{
	\begin{array}{l}
		\dot x =  f(x), \quad x\in \Omega \subset \R^n, \\
  \\
		y  =h(x) \in Y \subset \R^m
	\end{array}
	\right.
\end{equation}
and denote by $x(t,x_0)$ the solution of $\dot x=f(x)$ for the initial condition $x(0)=x_0$.
}

\medskip

\bleu{
So far we have studied observability as a property ensuring that  knowledge of a  measured ``signal'' $y(\cdot)$ results in   the uniqueness of the initial  condition $x_0$. In this chapter we will address the  {\em state estmation problem} which consist in obtaining an estimate $\hat{x}(t)$ of the state $x(t)$ of the system at time $t$, with the knowledge of the output $y(\cdot)$ up to time $t$.
}

\medskip

\bleu{We have for addressing this problem different possibilities}

\medskip

\bleu{
The first one is {\em given the output $y(\cdot)$ up to time $t$, find the possible initial condition $x_0$ which produces the same output $y(\cdot)$ up to time $t$}. Then, the estimate of $x(t)$ is given by the solution $x(t,x_0)$. The uniqueness of this problem is ascertained by the observability property of the system. This approach leads to the resolution of a {\em minimization problem}.
\[
\min_{x_0} \; \int_0^t   \parallel h (x(x_0,s) )-y(s) \parallel^2  \, ds .
\] 
In other words, we look,  by simulating the system for different initial conditions, to the ``best" one. The drawback of this approach is related to the difficulties of the non-linear minimization algorithms (existence of different local minimal, convergence speed...). Once one obtains an estimation $\hat x_0$ of the solution of this minimization problem, the estimation of the state $x(t)$ is given by $\hat x((t)=x(\hat x_0,t)$.
}

\medskip 

\bleu{
A second  approach is to differentiate the available outputs a number of times and then combine these derivatives appropriately to obtain the state vector.  
Formally, when the system is observable and we know (perfectly) enough derivatives of $y(\cdot)$ at a given time $t$, one just has to invert the map $\phi \mapsto (h(x),\mathcal L_fh(x),\cdots)$ at $(y(t),\dot y(t),\ddot y(t),\cdots)$ to reconstruct the state variable $x(t)$. In practice, it is known that numerically calculating derivatives on the raw signal data is imprecise and sensitive to measurement noise, especially if several successive derivatives have to be determined.  It is generally preferable to use a ``filter'' to smooth the data. For instance, polynomial functions or splines can approach with some regularity the measurements obtained over time, on which the derivative calculations can be performed before the inversion operation. This is an approximation method that does not guarantee an exact solution, and whose accuracy can be strongly influenced by the sensitivity of the solutions of the system with respect to the derivatives of the output.
}

\medskip

\bleu{
The  last approach is to look for a dynamical system whose "inputs" are the output of the observed system $y$, and whose "output" is an estimate $\hat{x}$ of the state of the original system, illustrated on the diagram of Figure \ref{fig_obs}. Such dynamical system is called an observer and is classical in the theory of control.}

\bleu{
\begin{figure}[h!]
\begin{center}
{\includegraphics[width=0.9\textwidth]{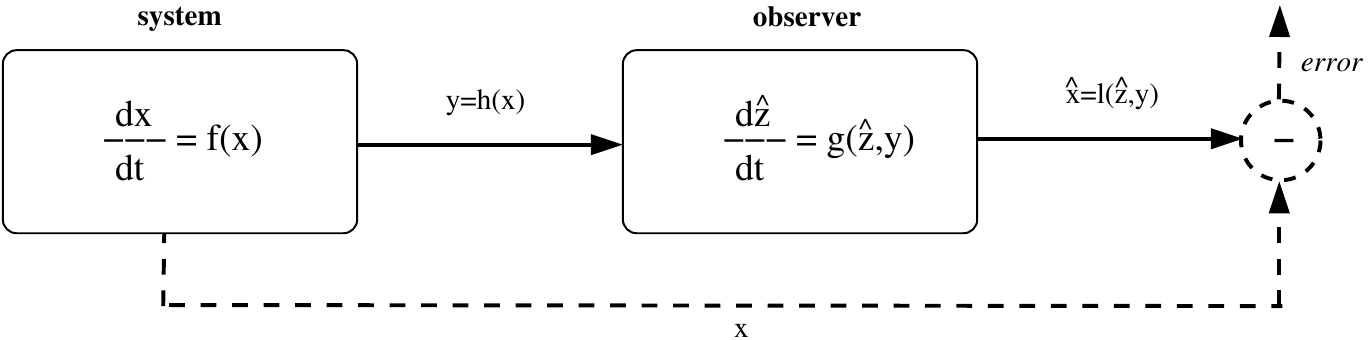}}
\caption{The observer concept}
\label{fig_obs}
\end{center}
\end{figure}
}

\bleu{
\begin{definition}
\label{def_obs}
An observer for system \eqref{sys_sect_obs} is a input-output system of the form
\begin{equation}
    \label{obs_gen}
\left\{\begin{array}{rl}
	\dot{\hat z} &= g(\hat z,y(t)), \; \hat z \in Z \subset \R^{n_z}\\
	\hat x&=l(\hat z,y(t))
\end{array}\right.
\end{equation}
such that $\hat{x}(t)$ is an asymptotic estimate of $x(t)$  satisfying 
\[
\lim_{t \rightarrow  +\infty}  || \hat{x}(t) - x(t) ||=0
\] 
for any initial condition $(x(0),\hat z(0)) \in \Omega\times Z$.
\end{definition}
}

\bleu{
In this definition, note that $n_z$ is not necessarily equal to $n$. It can be less than $n$ for {\em reduced-order} observers (see section \ref{secobsasymp} below) or larger than $n$ (see Remark \ref{rem_embedding} later on).
Let us also underline that this amounts to require the convergence of the solution of the coupled dynamics
\[
\left\{\begin{array}{rl}
	\dot x &= f(x) \\
	\dot{\hat z} &= g(\hat z,h(x))
\end{array}\right.
\; \Rightarrow \; \lim_{t \to +\infty} l(\hat z(t),h(x))-x(t) = 0
\]
for any initial condition in $\Omega\times Z$.
}

\medskip In this Chapter, we study the construction of observers and their theoretical convergence, without considering their practical performances in presence of measurements noise. This point will be addressed with more practical considerations in Chapter \ref{sectionPractical}.

\medskip Many observers are indeed of the form
\[
\dot{\hat x} = f(\hat x) + G(h(\hat x)-y(t))
\]
where $G$ is a constant matrix (i.e.~of the form \eqref{obs_gen} with $\hat z=\hat x$). Such observers are often called {\em Luenberger observers} \cite{Luenberger71}.
Note that this construction consists in a copy of the original dynamics $f$ plus a correcting term which depends on \bleu{the  {\em innovation} term $h(\hat x)-y(t)$, that is defined as the difference between the expected output $\hat y=h(\hat x(t)))$ if the true state was $\hat x(t)$ and the effective measured output $y(t)$}. 
Therefore, if $\hat x(t)$ and $x(t)$ are equal at a certain time $t$, it will remain identical at any future time. The main point is that the matrix $G$, \bleu{often called the {\em gains matrix}}, has to be chosen such that the estimation error $\hat x(t)-x(t)$ does converge to $0$, \bleu{possibly fast}. 
When $f$ and $h$ are linear and the system is observable, the theory of linear automatic control teaches that there always exists $G$ such that the convergence speed of the estimator can be chosen arbitrarily fast (see e.g. \cite{AndreaLara}). Obviously epidemiological models are rarely linear. However, looking for an Luenberger observer is often a first trial before considering more sophisticated estimators. Indeed, we shall see that for certain nonlinear dynamics, such observers do the job and in other cases, observers can be inspired from this form.

\medskip

\color{blue}
In this chapter, we do not pretend to present an exhaustive review of all the possible kinds of observers that exist in the literature (we make some comments at the end of the chapter). We focus on the most general classes of observers, for which the proof of the convergence can be shown, under some assumptions, and that can be derived in a systematic way from the equations of the model, once theses assumptions are fulfilled. 
We begin by some simple cases of observers for particular dynamics exploiting properties of linear systems, and then we consider a more general non-linear framework. The numerical implementation of these observers is illustrated in Chapter \ref{obspratic}.
\color{black}

\section{Observers with linear error dynamics}
\label{secobslinexp}

Consider systems of the form
\begin{equation}
	\label{linear_upto_output}
	\left\{\begin{array}{rl}
		\dot x &= Ax + \Phi(t,y)\\
  \\
		y &= Cx
	\end{array}\right.
\end{equation}
where $x \in \bleu{\Omega}$ and $y \in \R$. \bleu{We are looking for observers in the Luenberger form}
\begin{equation}
	\label{obsex1}
	\dot{\hat x} =A\hat x + \Phi(t,y(t)) + G(C\hat x-y(t)), \quad \bleu{\hat x \in \R^n}
\end{equation}
for which the error vector $e(t)=\hat x(t)-x(t)$ is solution of the linear dynamics
\begin{equation}
	\label{dynerr}
	\dot e =(A+GC)e
\end{equation}
The choice of the {\em gains} vector $G$ providing a convergence of $e(t)$ to $0$ comes directly from {\em the poles placement} technique of the theory of linear systems that we recall below.
\bleu{Let us underline that although $\Omega$ is assumed to be positively invariant by the dynamics
\eqref{linear_upto_output}, there is no reason for $\Omega$ to be invariant by \eqref{obsex1} because of the additive correction term. This is why we consider the dynamics \eqref{obsex1} in whole $\R^n$.}
\medskip

The following lemma gives a key result to construct the matrix $G$, which is well-known and often used in automatic control.

\begin{lemma}
	\label{lem_poles}
	
	Let $A$ be a square matrix of size $n$ and $C$ be a line vector of length $n$. If the observability matrix ${\cal O}$ defined in \eqref{obsmat}
	is of full rank, then for any set $\Lambda=\{\lambda_1,\cdots,\lambda_n\}$ of $n$ real or complex two-by-two conjugate numbers, there is a $G$ vector of size $n$ such that
	\[
	Sp(A+GC)=\Lambda\,.
	\]
	Specifically, if 
	\[
	\pi_{A}(\xi)=\xi^n+a_{1}\xi^{n-1}+\cdots +a_{n-1}\xi+a_{n}
	\]	
	is the characteristic polynomial of $A$, then one has
	\[
	G=P\Big[a_{n+1-i}+(-1)^{n-i}\sigma_{n+1-i}(\Lambda)\Big]_{i=1}^n
	\]
	where
	\[
	P=O^{-1}
	\left[\begin{array}{c}
		0\\ \vdots\\ \vdots\\ 1
	\end{array}\right]\, [I \; A \; \dots \; A^{n-1}]
	\]
	and the $\sigma_k$ designate the symmetric functions of the roots
	\begin{equation}
		\label{sigma}
		\sigma_k(\Lambda)=\sum_{1\leq i_1\leq \cdots\leq i_k\leq n} \!\!\!\! \lambda_{i_1}\lambda_{i_2}\cdots \lambda_{i_k}
	\end{equation}
\end{lemma}

\medskip

A proof of Lemma \ref{lem_poles} is given in Appendix \ref{appendixLemmaPoles}.

\medskip

\n {\em Remark.} This result can be generalized to vectorial observations, i.e.~for matrices $C$ with $m>1$ rows and $n$ columns.

\bigskip

Finally, by choosing  numbers $\lambda_i$ with negative real parts, one can make the convergence of the error given by the exponentially decreasing  dynamics \eqref{dynerr} as fast as desired. The observer \eqref{obsex1} is thus adjustable.
It should be noted that when the difference $\hat z_1(t) -y(t)$ (usually called "innovation") becomes and remains close to $0$, the trajectories of the observer follow those of the system: we can then consider that the observer has practically converged. The innovation is thus very useful in practice because it provides information on the current stage of convergence of the estimate.

\medskip

We illustrate this technique on a population model with age classes.

\medskip

\begin{exemple}
	\label{ex_obs_age}
	Let us consider a population structured in three stages: young, subadult and adult, of stocks $x_1$, $x_2$, $x_3$ respectively. It is assumed that only adults $x_3$ can reproduce, giving birth to young $x_1$:
	\begin{equation}
		\label{model3stages}
		\left\{\begin{array}{l}
			\dot x_1 = -a_1x_1 -m_1x_1 +r(t,x_3)\\
			\dot x_2 = a_1x_1-a_2x_2 -m_2x_2\\
			\dot x_3 = a_2x_2 -m_3x_3\\\
   \\
			y = x_3
		\end{array}\right.
	\end{equation}
	The coefficients $a_i$ are the transition rates between age classes, $m_i$ are the mortality rates of each class, and $r(\cdot)$ is the reproductive function (usually non-linear and seasonally dependent), for example
	\[
	r(t,x_3)=\frac{\bar r(t)x_3}{k+x_3} , \quad \bar r(t) \in [\bar r_{min},\bar r_{max}]\,.
	\]
	Here, it is assumed that only the size of the adult class is measured over time. The aim is to estimate the stocks of larvae and subadults over time. The model \eqref{model3stages} is of the form \eqref{linear_upto_output} \bleu{with $\Omega=\R_+^3$}, where we posed
	\[
	A=\left[\begin{array}{ccc} -(a_1+m_1) & 0 & 0 \\\ a_1 & -(a_2+m_2) & 0\\\ 0 & a_2 & -m_3+0
	\end{array}\right], \quad C=[0\; 0\; 0\; 1]
	\]
	and
	\[
	\Phi(t,x_3)= \left[\begin{array}{c} 0\\ r(t,x_3)\end{array}\right]
	\]
	One can check that the observability matrix ${\cal O}$ si full rank, or alternatively directly check that the system is observable. Indeed, we obtain by using the expression $\dot x_3$:
	\[
	x_2= \frac{\dot y+m_3 y}{a_2} \,,
	\]
	then with the expression $\dot x_2$:
	\[
	x_1= \frac{\dot x_2 + (a_2+m_2)x_2}{m_1}=\frac{\ddot y + (a_2+m_2)m_3 y}{a_2m_2}\,.
	\]
	Therefore, the following system \bleu{in $\R^3$}
	\begin{equation}
		\label{obs3stades}
		\left\{\begin{array}{rl}
			\dot{\hat x}_1 &= -a_1\hat x_1 -m_1\hat x_1 +r(t,y(t)) + G_1(\hat x_3-y(t))\\
			\dot{\hat x}_2 &= a_1\hat x_1-a_2\hat x_2 -m_2\hat x_2+ G_2(\hat x_3-y(t))\\
			\dot{\hat x}_3 &= a_2\hat x_2 -m_3\hat x_3+ G_3(\hat x_3-y(t))\
		\end{array}\right.
	\end{equation}
	with well chosen gains $G_1$, $G_2$, $G_3$ is an observer for the dynamics \eqref{model3stages}, with an exponential convergence.
\end{exemple}

\color{blue}

In the next section, we study a more general kind of non-linearity.

\section{Observers for systems with Lipschitz non-linearity}

Let us consider a system of the form
\begin{equation}
\label{sysLip}
    \left\{\begin{array}{lll}
    \dot x & = & Ax + \phi(x)\\
    \\
    y & = & Cx
    \end{array}\right.
\end{equation}
where
\begin{enumerate}
\item the pair $(A,C)$ is observable i.e. the observability matrix ${\cal O}$ defined in \eqref{obsmat} is of full rank,
\item the map $\phi$ is globally Lipschitz on $\R^n$,
with a Lipschitz constant $\ell$.
\end{enumerate}
For this non-linear system, consider an observer in Luenberger form
\begin{equation}
    \label{obsLip}
    \dot{\hat x} = A\hat x + \phi(\hat x) + G(C\hat x-y(t))
\end{equation}
We recall the Lyapunov Theorem \cite{Parks92}, which is useful in many situations.
\begin{theorem}
\label{LyapunovTh}
Given a symmetric positive definite matrix $Q$, there exists a unique symmetric definite positive matrix $P$ satisfying $M^\top P+PM+Q=0$ if and only if the linear system $\dot x=Mx$ is globally exponentially stable.
\end{theorem}
A proof of this Theorem is given in Appendix \ref{appendixPropLyapunov}.
Then, one has the following result
\begin{proposition}
\label{prop_obsLip}
    Take $G$ such that $A+GC$ is Hurwitz\footnote{A matrix is Hurwitz if all its eigenvalue have negative real parts.}, then if $\ell$ is small enough, \eqref{obsLip} is an exponential observer of \eqref{sysLip}, in the sense that there exist $\alpha>0$ and $\beta >0$ such that
\[
||\hat x(t)-x(t)|| \leq \alpha ||\hat x(0)-x(0)||e^{-\beta t}, \quad t \geq 0
\]
for any initial condition $x(0)$, $\hat x(0)$. 
\end{proposition}

\begin{proof}
Accordingly to Theorem \ref{LyapunovTh}, there exists an unique symmetric definite matrix $P$ satisfying
\begin{equation}
\label{Lyap_A+GC}
(A+GC)^\top P + P(A+GC)+I = 0
\end{equation}
(where $I$ is the identity matrix).
Let us consider the candidate Lyapunov function
\[
V(e)=e^\top P e
\]
and the time function  $v(t)=V(e(t))$. One gets
\begin{eqnarray*}
\dot v & = & e^\top\Big( (A+GC)^\top P + P(A+GC) \Big)e + 2 e^\top P (\phi(\hat x)-\phi(x))\\
& = & - ||e||_2^2 + 2 e^\top P (\phi(\hat x)-\phi(x))\\
& \leq & - ||e||_2^2 + 2  \ell ||P||_2 ||e||_2^2 = -(1-2 \ell ||P||_2)||e||_2^2
\end{eqnarray*}
If the Lipschitz constant $\ell$ is small, one has 
\begin{equation}
    \label{condobsLip}
    1-2  \ell ||P||_2 > 0 .
\end{equation}

On the other hand, as $P$ is a positive definite matrix, one has
\[
\lambda_{min}||e||_2 \leq V(e) \leq \lambda_{max}||e||_2
\]
where $\lambda_{min}$, $\lambda_{max}$ are the smallest and largest eigenvalues of $P$. Then one obtains $\dot v \leq -2\beta v$ where
\[
\beta=\frac{1-2 \ell ||Pe||_2}{2\lambda_{min}} >0 
\]
which implies
\[
v(t) \leq v(0)e^{-2\beta t} \Rightarrow ||e(t)||^2 \leq \alpha^2 ||e(0)||^2e^{-2\beta t}
\]
with $\alpha=\sqrt{\frac{\lambda_{max}}{\lambda_{min}}}$. Finally, one gets
\[
||e(t)|| \leq \alpha ||e(0)||e^{-\beta t}, \quad t \geq 0
\]
which proves that the estimation error $e$ converges exponentially to $0$.
\end{proof}

\bigskip


There exist other techniques to design a gain $G$ based on the Riccati equation rather than the Lyapunov equation \eqref{Lyap_A+GC}, but that are more technical and thus out of the scope of this book (we refer to \cite{Aboky} for interested readers).

\begin{remark}
If the map $\phi$ can be written as $\phi(x)=\varphi(h(x),x)=\varphi(y,x)$, where $\varphi$ is globally Lipschitz w.r.t.~$x$ uniformly in $y$, then one can consider the observer
\begin{equation}
    \label{obsLipgen}
    \dot{\hat x} = A\hat x + \varphi(y(t),\hat x) + G(C\hat x-y(t))
\end{equation}
which generalizes the observers \eqref{obsex1} and \eqref{obsLip}.
\end{remark}

\medskip

\begin{exemple}
\label{exSIRS}
We consider the "SIRS" model \cite{KOROBEINIKOV}, that is the SIR model with loss of immunity, assuming that the size of the recovered population is tracked over time
\[
\left\{\begin{array}{l}
\dot S  =  -\beta SI + \mu R\\
\dot I  =  \beta SI - \gamma I\\
\dot R  =  \gamma I\\
\\
y=R
\end{array}\right.
\]
Here $S$, $I$ and $R$ denote the densities of susceptible, infected and recovered populations, so that one has $S+I+R=1$ at any time. 
We assume that all the parameters are known, and we aim at estimating $S$ and $I$.
Measuring $R$ is equivalent to measuring $S+I$, and we rewrite this model as follows
\[
\left\{\begin{array}{l}
\dot S  =  -\beta I +\varphi(S,I) + \mu (1-S-I)\\
\dot I  =  \beta I -\varphi(S,I)- \gamma I\\
\\y=S+I
\end{array}\right.
\]
where we posit
\[
\varphi(S,I)=\beta (1-S)I
\]
If we consider the initial stage of an epidemics, variables $I$ and $S$ are close to $0$ and $1$ respectively, that is on a domain $\Delta:=\{ 1-\varepsilon < S < 1; \; 0 < I < \varepsilon \}$ for some $\varepsilon>0$. This domain is not invariant by the dynamics but we shall consider the estimation problem on a time windows for which the solution stays in this set $\Delta$. One can check that the function $\varphi$ is Lipschitz with a constant $\ell = 2\varepsilon$ on $\Delta$, and can consider an extension of $\varphi$ outside $\Delta$ with the same Lipschitz constant 
$\ell$ but on all $\R^2$, for instance
\[
\tilde \varphi(I,S)=\varphi(sat_{[0,\varepsilon]}sat_{[0,\varepsilon]}(I))
\]
where $sat_{[\,]}$ denotes the saturation function
	\[
	sat_{[a,b]}(x)= \max(a,\min(b,x))\,.
	\]
On $\Delta$, the system can be written as follows
\begin{equation}
    \label{dyn_ext}
    \begin{array}{l}
\frac{d}{dt}\left[\begin{array}{c}S\\I\end{array}\right]=
\underbrace{\left[\begin{array}{ccr}-\mu && -\beta-\mu\\ 0 && \beta-\gamma\end{array}\right]}_{A}\left[\begin{array}{c}S\\I\end{array}\right]+\underbrace{\left[\begin{array}{r}\tilde\varphi(S,I)\\-\tilde\varphi(S,I)\end{array}\right]}_{\phi(S,I)}\\
y =\underbrace{[1 \; 1]}_C \left[\begin{array}{c}S\\I\end{array}\right]
\end{array}
\end{equation}
One can check that the pair $(A,C)$ is observable, as the matrix
\[
\left[\begin{array}{c}
C\\
CA
\end{array}\right] = 
\left[\begin{array}{ccc}
1 && 1\\
-\mu && -\mu-\gamma
\end{array}\right] 
\]
is full rank.
The  system \eqref{dyn_ext} is well defined on the whole $\R^2$, and $\phi$ is a globally Lipschitz map, with a Lipschitz constant equal to $2\ell$. We are thus in position to apply Proposition \ref{prop_obsLip}: if $\varepsilon$ is small enough, the system
\[
\frac{d}{dt}\left[\begin{array}{c}\hat S\\ \hat I\end{array}\right]=
A \left[\begin{array}{c}\hat S\\ \hat I\end{array}\right] + \phi(\hat S,\hat I) + G(\hat S+\hat I - y(t))
\]
where $G$ is such that $A+GC$ is Hurwitz, is an exponential observer of \eqref{dyn_ext}.

\medskip

In Section \ref{obsLipnum}, we discuss the applicability of this observer with numerical values.
\end{exemple}

\color{black}

\medskip

It may happen that the estimation error of an observer is only partially assignable, as we shall see in the next example.

\section{Observers via \bleu{decoupled variables}}
\label{secobschtvar}

Consider, as in Section \ref{secidentchgvar}, systems
\[
\left\{\begin{array}{lll}
	\dot x & = & f(x) , \quad \bleu{x \in \Omega}\\
 \\
	y &= & h(x)
\end{array}
\right.
\]
for which there exists a change of coordinates $x \bleu{\in \Omega} \mapsto w=g(x,h(x)) \bleu{\in {\cal W}\subset \R^{n^\prime}}$ such that one has
\[
\dot w(t) = l(w(t),h(x(t))), \quad t \geq 0
\]
for any solution $x(\cdot)$ \bleu{in $\Omega$}, with the properties
\begin{enumerate}
	\item $\left\{w=g(x,h(x)) \Longleftrightarrow x = \tilde g(w,h(x))\right\}, \; \bleu{(x, w) \in \Omega\times{\cal W}}$
	\item $\left\{w=g(x,h(x)) \Longrightarrow h(x)=k(w)\right\}, \; \bleu{(x,w) \in \Omega\times{\cal W}}$
\end{enumerate}
where $g$, $\tilde g$, $l$ and $k$ are smooth maps. 

Then, one can look for an observer $\hat w(\cdot)$ of the system
\[
\left\{\begin{array}{lll}
	\dot w & = & l(w,h(x))\\
 \\
	y &= & k(w)
\end{array}
\right.
\]
and take, as an estimator of $x(\cdot)$
\begin{equation}
	\label{estim_org_cooord}
	\hat x(t)=\tilde g(\hat w(t),y(t))\,.
\end{equation}
There is an advantage of considering such a change of coordinates when the maps $g$, $\tilde g$, $l$ and $k$ are independent of a parameter $\theta$ present in the expression of $f$, as in Proposition \ref{prop_ident_chtvar} of Section \ref{secidentchgvar}. However, not that the estimator \eqref{estim_org_cooord} does not filter the measurement $y(\cdot)$ and might be sensitive to noise. 

\medskip

Let us illustrate this approach on the malaria model (Example \ref{example_malaria}). 

\medskip

\begin{exemple}
	\label{example_malaria_obs}
	Consider the model \eqref{AMGreduit} of Example \ref{example_malaria}.
	With the variable $w=x-Ey=(I-EC)x$ in $\R^7$, the dynamics \eqref{AMGx} is independent of the unknown non-linear term $\beta S M$:
	\[
	\left\{\begin{array}{lll}
		\dot w & = & \bar Aw+\bar AEy+\Lambda\, e_1\\
  \\
		y  & = & Cx
	\end{array}\right.
	\]
	and one can consider the following observer for system \eqref{AMGx} in Luenberger form
	\begin{equation}\label{AMGobs}
		\left\{\begin{array}{lcl}
			\dot{\hat w}&=& \bar A\hat w+\bar AEy+\Lambda\, e_1+L(y(t)-C(\hat w+Ey(t))\\\
    & = & (\bar A-LC)\,\hat w+\big(L+(\bar A-LC)E\big)\,y(t) +
			\Lambda\, e_1,\\[3mm]
			\hat x(t) &=& \hat w(t)+E\,y(t) .
		\end{array}
		\right.
	\end{equation}
	where $L$ is a gains vector in $\mathbb{R}^7$ to be chosen.
	The dynamics of the error $e(t)=\hat x(t)-x(t)$ is given by
	\[
	\dot e = (\bar A-LC)e
	\]
	Note that one has $C\bar A=0$. Therefore the rank of the observability matrix of the pair $(\bar A,C)$ is equal to one, and $0$ is an eigenvalue of $\bar A$. The choice of $L$ allows then to assign only one eigenvalue of $\bar A-LC$, equal to $-(L_2+L_3)$, the other eigenvalues remaining negative.
	Therefore \eqref{AMGobs} is  an observer for system \eqref{AMGx} with exponential convergence, that does not use the unknown parameter $\beta$.
\end{exemple}

\begin{remark}
	\label{remobsnonassis}
	Differently to observers of Section \ref{secobslinexp}, one cannot expect a convergence speed of the observer \eqref{AMGobs} faster than the dynamics \eqref{AMGreduit}, because the error dynamics is not completely assignable. However, the convergence is exponential. This is illustrated with numerical simulations in Section \ref{obsmalarianum}.
\end{remark}

\section{Reduced-order observers}
\label{secobsasymp}\m

A typical situation is when one can operate a state decomposition \bleu{when $m<n$}, as follows

\begin{enumerate}
	\item decompose the state vector $x$ (may be at the price of a change of variables) as
	\[
	x = \left[\begin{array}{c}
		y\\
		x_u
	\end{array}\right]
	\]
	where $x_u \in \bleu{\R^{n-m}}$ represent the {\em unmeasured} variables,
	\item look for an auxiliary variable (that we called $z$) whose dynamics is independent of $x_u$ 
	\[
	\dot z = g(z,y(t)), \; \bleu{z \in \R^q}
	\]
    \bleu{(for some $q$)}
	and asymptotically stable \bleu{(that is any solution $z(\cdot)$ converges to $0$ when $t \to +\infty$)}, and such that $x_u$ can be globally expressed as
	\[
	x_u=l(z,y)
	\]
	where $l$ is a smooth map (say $C^1$).
\end{enumerate}\m

Then, the dynamics
\[
\left\{\begin{array}{ll}
	\dot{\hat{z}} & =  g(\hat z,y(t))\\
	\hat x_u & = l(\hat z,y(t))	
\end{array}
\right.
\]
is an {\em asymptotic observer}, whose error convergence $\hat x_u-x_u$ is simply provided by the asymptotic convergence of $\hat z_z$ to $0$, whatever is the initial condition $\hat z(0) \bleu{\in \R^q}$. When the convergence speed of an estimator cannot be adjusted, it is usually called an {\em asymptotic observer}, differently to the previous section for which the error convergence can be made arbitrarily fast. Note that differently to the previous section, these observers have no tuning parameters and are not driven by innovation terms. These estimators are  {\em reduced-order} observers when the variable $z$ is of lower dimension than $x$ \bleu{(i.e.~when $q<n$)}. An interest for such observers is that it can possess robustness features when the maps $g$ and $l$ are independent of some terms or parameters of the dynamics $\dot x=f(x)$.
Let us illustrate this feature on the Kermack-McKendrick model with fluctuating rates.

\medskip

\begin{exemple}
	\label{ex_SIR_fluct}
	We consider the SIR model with birth and death terms
	\begin{equation}
		\label{SIRfluctuating}
		\left\{\begin{array}{rl}
			\dot S &= -\beta(t)SI  +\nu N- \mu S\\
			\dot I & =\beta(t)SI -\rho(t)I-\mu I\\
			\dot R &= \rho(t)I-\mu R
		\end{array}\right.
	\end{equation}
	where parameters $\beta$ and $\rho$ fluctuate unpredictably over time.
	We assume, for simplicity, that the birth rate $\nu$ is equal to the death rate $\mu$, so that the total population remains constant of size $N=S+I+R$ (assumed to be known). Let us suppose that the size of the infected population is monitored over time as well as the number of new cured individuals, which amounts considering that the observation vector at time $t$ is
	\[
	y(t)=\left[\begin{array}{c}y_1(t)\\ y_2(t) \end{array}\right] = \left[\begin{array}{c}I(t)\\ \rho(t)I(t) \end{array}\right]\,.
	\]
	Stocks of classes $S$ and $R$ are not initially known. Then, the system
	\begin{equation}
		\label{obsSIRrobust}
		\left\{\begin{array}{rl}
			\dot Z  &= \nu N -y_2(t)-\mu Z\\
			\hat S &= Z-y_1(t)\\
			\hat R & = N-Z
		\end{array}\right.
	\end{equation}
	is an observer allowing to estimate $S$ and $R$ without knowing $\beta(\cdot)$ and $\rho(\cdot)$. Indeed, the dynamics of the estimators verifies
	\begin{align*}
		&  \frac{d}{dt} (\hat S - S) = \dot Z - \dot y_1 -\dot S = -\mu(Z-S-I)=-\mu(\hat S-S)\\[2mm]
		& \frac{d}{dt} (\hat R - R) = -\dot Z - \dot R = -\nu N+\mu(Z+R)=-\mu(\hat R-R)
	\end{align*}
	which ensures the convergence of the $\hat S$ and $\hat R$ estimates. Note that the internal dynamics of the observer is here of smaller dimension than the system, and that the estimate of the unmeasured state variable $S$ is a function of the internal state $Z$ of the observer and the observation $y_1$. The speed of convergence of this observer is not adjustable, but it has the advantage of being perfectly robust to any (unknown) variations of the terms $\beta(\cdot)$ and $\rho(\cdot)$. This is illustrated with numerical simulations in Section \ref{obssirrobustnum}.
\end{exemple}


\section{\bleu{The high-gain observer for nonlinear systems}}

In the two previous examples, the dynamics of the estimation error was linear. For an observable non-linear system, the existence of an observer whose estimation error is linear is not guaranteed. This is a difficult problem. However, one can consider the (nonlinear) {\em observability canonical form} \bleu{in $\R^n$} \cite{GauthierKupka2001} (given here for a scalar output \bleu{i.e.~for $m=1$})
\begin{align}
	& \label{form_norm1}\dot z =F(z):=\underbrace{\left[\begin{array}{cccccc}
			0 & 1 & 0 & \cdots & &\\
			0   &  0  & 1 & 0 &  \cdots &\\
			& & \ddots & \ddots &  &\\
			& && \ddots & \ddots &\\
			& &  &  & 0 & 1\\
			& &  &  &  & 0 \end{array}\right]}_{A}z+\psi(z)\underbrace{\left[\begin{array}{c} 0 \\ \vdots\\ \\ \\ 0\\ 1\end{array}\right]}_B\\[4mm]
	& \label{form_norm2} y = \underbrace{[ 1 \, 0 \, \cdots \, \cdots \, 0]}_{C}z
\end{align}
where the function $\psi$ is Lipschitz on $\R^n$. Then, one can show that there exists an observer of the Luenberger form
\begin{equation}
	\label{obs_normal}
	\dot{\hat z}=F(\hat z) + G(C\hat z -y(t))
\end{equation}
with exponential convergence when $G$ is a well-chosen gains vector. 
When an observable system
\[
\left\{\begin{array}{rl}
	\dot x &= f(x) , \quad x \in \R^n\\
 \\
	y &= h(x) , \quad y \in \R
\end{array}\right.
\]
is not in normal form, but the application
\[
\phi_n(x)=\left[\begin{array}{c}
	h(x)\\
	\mathcal L_f h(x)\\
	\vdots\\
	\mathcal L_f^{n-1}h(x)
\end{array}\right]
\]
is a diffeomorphism\footnote{a diffeomorphism is an invertible map such that both the map and its inverse are differentiable.} from $\R^n$ into $\R^n$ and the function
\[
\psi(z):=\mathcal L^m_f h \circ \phi_n^{-1}(z)
\]
is Lipschitz on $\R^n$, then the observer \eqref{obs_normal} can be written in the $x$ coordinates as follows
\[
\dot{\hat x} = f(\hat x)+[J\phi_n(\hat x)]^{-1} G(h(\hat x)-y(t))
\]
where $J\phi_n(x)$ denotes the Jacobian matrix of $\phi$ at $x$. The observer preserves the Luenberger structure but with variable gains.	

\bigskip

Let us first note that the pair $(A,C)$ as defined in \eqref{form_norm1}-\eqref{form_norm2} is observable. Indeed, we have $O=Id$. Thus, according to the Lemma \ref{lem_poles}, one can freely assign the spectrum of $A+GC$ by choosing the vector $G$. We show now how to choose the eigenvalues of $A+GC$ to ensure the convergence of the non-linear observer \eqref{obs_normal}.
To do this, we begin by giving some properties of the Vandermonde matrices
\[
V_{\lambda_1,\cdots,\lambda_n}:=\left[\begin{array}{ccccc}
	\lambda_1^{n-1} & \lambda_1^{n-2} & \cdots & \lambda_1 & 1\\
	\lambda_2^{n-1} & \lambda_2^{n-2} & \cdots & \lambda_2 & 1\\
	\vdots & \vdots &  & \vdots & \vdots \\
	\lambda_n^{n-1} & \lambda_n^{n-2} & \cdots & \lambda_n & 1
\end{array}\right]
\]
related to the normal form.
\begin{lemma}
	\label{lemVandermonde}
	Let $\Lambda=\left\{ \lambda_1,\cdots,\lambda_n \right\}$ be a set of $n$ distinct real numbers and $G$ a vector such that 
	\[
	Sp(A+GC)=\Lambda:=\left\{ \lambda_1,\cdots,\lambda_n \right\}
	\]
	Then
	\[
	V_{\lambda_1,\cdots,\lambda_n }(A+GC)V_{ \lambda_1,\cdots,\lambda_n }^{-1}=
	\left[\begin{array}{cccc}
		\lambda_1 & & &\\
		& \lambda _2 & &\\
		& & \ddots &\\
		& & & \lambda_n
	\end{array}\right]
	\]
	Moreover, for any numbers $c>0$ and $\theta>0$, there exist $\lambda_n < \lambda_{n-1}< \cdots < \lambda_1<0$ such that
	\[	
	\lambda_1+c||V_{\lambda_1(\theta),\cdots,\lambda_n(\theta)}^{-1}||_\infty=-\theta		
	\]
\end{lemma}

The proof of Lemma \ref{lemVandermonde} is given in Appendix \ref{appendixLemmaVandermonde}.

\bigskip
We are now ready to show the convergence of the observer \eqref{obs_normal} in coordinates $z$, for a gains vector $G$ such that $A+GC$ has $n$ distinct eigenvalues $\lambda_1$, ..., $\lambda_n$ of negative real values.  Denote the error $e=\hat z-z$. We have
\[
\dot e = (A+G C)e+B(\psi(\hat z)-\psi(z))
\]
Let $\xi=V e$ where $V$ designates the Vandermonde matrix $V_{\lambda_1,\cdots,\lambda_n}$. Thanks to Lemma \ref{lemVandermonde}, we obtain
\[
\dot \xi = \Delta \xi + V B (\psi(\hat z)-\psi(z))
\]
where $\Delta$ is the diagonal matrix $diag(\lambda_1,\cdots,\lambda_n)$.
By multiplying on the left by $\xi^\top$, one obtains
\[
\begin{array}{lll}
	\xi^\top\dot \xi & = & \xi^\top\Delta \xi + \xi^\top V B (\psi(\hat z)-\psi(z))\\
	& \leq  & \lambda_1||\xi||^2 + ||\xi||\sqrt{n}|\psi(\hat z)-\psi(z)|\\
	& \leq  & \lambda_1||\xi||^2 + ||\xi||\sqrt{n}L||e||\\
	& \leq  & \big(\lambda_1+\sqrt{n}L||V^{-1}||_\infty\big)||\xi||^2
\end{array}
\]
where $L$ is the Lipschitz constant of $\psi$.
Thus the norm of $\xi$ verifies
\[
||\xi(t)|| \leq ||\xi(0)|| + \int_0^t \big(\lambda_1+\sqrt{n}L||V^{-1}||_\infty\big)||\xi(\tau)||d\tau
\]
and by Gronwall's Lemma, we obtain
\[
||\xi(t)|| \leq ||\xi(0)||e^{\big(\lambda_1+\sqrt{n}L||V^{-1}||_\infty\big)t}
\]
Finally, for any $\theta>$0, Lemma \ref{lemVandermonde} gives the existence of numbers 
$\lambda_n <\lambda_{n-1}< \cdots < \lambda_1<0$ such that
$||\xi(t)|| \leq ||\xi(0)||e^{-\theta t}$ for any $t>0$, which guarantees the exponential convergence of the error $e$ to $0$.

\bigskip 

The observer \eqref{obs_normal} with the gains vector $G_\theta$ is called {\em high gain observer} \cite{GauthierKupka2001}, because the value of $\theta$ must be "sufficiently" large, and its successive powers might take large values.

\begin{remark}
\label{rem_embedding}
	In practice, the function $\psi$ is not necessarily globally Lipschitz on $\R^n$, \bleu{and even not properly defined outside $\Omega$, while the observer \eqref{obs_normal} needs to defined on whole $\R^n$}. Nevertheless, if there exists
	a compact subset $K$ of \bleu{$\Omega$} that is forwardly invariant by the dynamics  \eqref{form_norm1}, one can consider an extension of $\psi$ outside $K$ that is globally Lipschitz on $\R^n$ and define then the observer on whole $\R^n$ \bleu{(see for instance \cite{HGO92,RapaportMaloum04})}. \bleu{This is illustrated on Example \ref{ex_highgainSIR} below.}
    \bleu{It could also happen that for observable systems the map $\phi_n$ is not injective (and thus cannot be a  diffeomorphism), but $\phi_{\tilde n}$ is injective for $\tilde n > n$. Then, it is theoretically possible to embed the system in $\R^{\tilde n}$, that is to write the dynamics with a state vector in dimension $\tilde n >\tilde n$ and then to construct a high gain observers in $\R^{\tilde n}$. Such construction is beyond the present book (see \cite{RapaportMaloum04} for some techniques to build such extensions).}
\end{remark}

\medskip

Let us now illustrate this construction on the Kermack-McKendrick model.

\medskip

\begin{exemple}
	\label{ex_highgainSIR}
	
	We consider the classical SIR model
	\begin{equation}
		\label{SIRclassical}
		\left\{\begin{array}{rl}
			\dot S &= -\beta SI\\
			\dot I & =\beta SI -\rho I\\
			\dot R &= \rho I
		\end{array}\right.
	\end{equation}
	where the parameters $\beta$ and $\rho$ are known, and suppose that the only observation is the cumulative number of recovered individuals since a time $t_0$
	\[
	y(t) = \int_{t_0}^t \rho I(\tau)d\tau\,.
	\]
	It is also assumed that the size of the total population $N=S+I+R$ is known.
	To put the system in canonical form, we write
	\[
	\begin{array}{lll}
		z_1 & = & y\\
		z_2 & = & \dot z_1 = \rho I\\
		z_3 & = & \dot z_2 = (\beta S -\rho)z_2
	\end{array}
	\]
	and
	\[
	\dot z_3 =(\beta S -\rho)z_3 - \beta(\beta S I)z_2
	= \psi(z):=\frac{z_3^2}{z_2}-\frac{\beta}{\rho}\left(\frac{z_3}{z_2}+\rho\right)z_2^2\,.
	\]
	We can then reconstruct the $I$, $S$ and $R$ stocks from the variables $z$ as follows
	\begin{align*}
		& S = \frac{1}{\beta}\left(\frac{z_3}{z_2}+\rho\right)\\
		& I = \frac{z_2}{\rho}\\
		& R = N-\frac{1}{\beta}\left(\frac{z_3}{z_2}+\rho\right)-\frac{z_2}{\rho}
	\end{align*}
	Note that $\psi$ is not globally Lipschitz on $\R^3$, and has a singularity at $z_2=0$. Nevertheless, we notice that the term $z_3/z_2$ can be framed as follows
	\[
	\frac{z_3}{z_2}=\beta S - \rho \in [-\rho,\beta-\rho]
	\]	
	and that one has
	\[
	\dot z_3 \in [-\rho^3 N -\rho\beta^2N^3,\rho(\beta-\rho)N]\,.
	\]
	We can therefore consider the expression
	\[
	\tilde\psi(z)= sat_{[-\rho^3 N -\rho\beta^2N^3,\rho(\beta-\rho)N]}\left(sat_{[-\rho,\beta-\rho]}\left(\frac{z_3}{z_2}\right)z_3 -\frac{\beta}{\rho} z_3z_2-\beta z_2^2)\right)
	\]
	instead of $\psi(z)$, where $sat_{[\,]}$ denotes the saturation function
	\[
	sat_{[a,b]}(x)= \max(a,\min(b,x))\,.
	\]
	Finally, we choose the gains  $G_i$ of the observer such that $Sp(A+GC)=\{\lambda_1,\cdots,\lambda_n\}$ with 
	$\lambda_n <\lambda_{n-1}< \cdots < \lambda_1<0$ and $\lambda_1$ {\em enough} negative. This amounts to take
	$G_i=\sigma_i(\{\lambda_1,\cdots,\lambda_n\})$. One thus obtains the internal dynamics of the observer
	\begin{equation}
		\label{obsSIRhighgain}
		\left\{\begin{array}{lll}
			\dot{\hat z}_1 & = &  \hat z_2 +(\lambda_1+\lambda_2+\lambda_3)(\hat z_1 - y(t))\\
			\dot{\hat z}_2 & = &   \hat z_3 +(\lambda_1\lambda_2+\lambda_1\lambda_3+\lambda_2\lambda_3)(\hat z_1 - y(t))\\
			\dot{\hat z}_3 & = & \tilde\psi(\hat z) +\lambda_1\lambda_2\lambda_3(\hat z_1 - y(t))
		\end{array}\right.
	\end{equation}
	and the estimators
	\begin{equation}
		\label{estSIRhighgain}
		\left\{\begin{array}{lll}
			\hat S & = &  \displaystyle \frac{1}{\beta}\left(sat_{[-\rho,\beta-\rho]}\left(\frac{\hat z_3}{\hat z_2}\right)+\rho\right)\\[3mm]
			\hat I & = &   \displaystyle \frac{\hat z_2}{\rho}\\[3mm]
			\hat R & = & N-\hat S- \hat I
		\end{array}\right.
	\end{equation}
\end{exemple}

\section{Discussion}

The construction of an observer can avoid in certain situations to study the identifiability. For instance, in Sections \ref{secobschtvar} and \ref{secobsasymp}, the observers do not require the knowledge of all the parameters of the model, and even of some functions involved in the model. This is why such observers are also called {\em unknown-inputs observers}. The theory of unknown-inputs observers has been mainly studied for linear systems \cite{MR1164570,ChenSaif2006,Nazari2015}. Very few general results are available for nonlinear systems (this research field is today largely open).

\medskip

The  existence of observers without the possibility of fixing arbitrarily the speed of convergence, as in Examples \ref{example_malaria_obs} and \ref{ex_SIR_fluct}, is connected to the property of {\em detectability} (see for instance \cite{MD08,ABS13}), which is a weaker than observability: a system
\[
\left\{\begin{array}{l}
	\dot x = f(x), \; \bleu{ x \in \Omega}\\
 \\
	y=h(x)
\end{array}\right.
\]
is {\em detectable} \bleu{(in $\Omega$)} if for any pair of solutions $x^a(\cdot)$, $x^b(\cdot)$ \bleu{in $\Omega$}, one has
\[
\left\{ h(x^a(t))=h(x^b(t)), \; t \geq 0 \right\} \; \Rightarrow \;
\lim_{t \to +\infty} x^a(t)-x^b(t) =0\,.
\]

For linear dynamics $\dot x=Ax$, $y=Cx$ that are not observable, there exists a {\em Kalman decomposition} \cite{MR0152167} i.e.~an invertible matrix $P$ such that 
\[
PAP^{-1}=\left[\begin{array}{cc}A_{11} & 0\\
	A_{21} & A_{22}
\end{array}\right], \quad CP^{-1}= \left[\begin{array}{cc}C_{1} & 0
\end{array}\right]
\]
with $A_{11} \in {\cal M}_{l\times l}$, $C_{1} \in {\cal M}_{1 \times l}$ where $l<n$ is the rank of the observability matrix ${\cal O}$ recalled in \eqref{obsmat}, such that the subsystem  
$\dot z =A_{11}z$, $y=C_{1}z$ is observable.
Then, the system is detectable when the matrix $A_{22}$ is Hurwitz. This is exactly the case of Example \ref{example_malaria_obs} with $l=1$.

\bigskip

\color{blue}
Let us mention the more recent technique proposed by Kazantzis and Kravaris for obtaining Luenberger like observers for non-linear systems. It consists in looking at a (non-linear) change of coordinates as a diffeormorphism $T$: $x \mapsto \zeta=T(x)$ such that the dynamics of the variable $\zeta$ writes
\[
\dot \zeta = A\zeta + B(h(T^{-1}(x)))
\]
where $A$ is a Hurwitz matrix and $B$ a smooth map. Then,
\[
\dot z = AZ+B(y), \quad \hat x=T^{-1}(z)
\]
is a natural observer, whose speed of convergence is given by the spectrum of the matrix $A$. However, the map $T$ has to be found as a solution of the partial derivative equation
\[
\frac{\partial T}{\partial x}(x)f(x)=AT(x)+B(h(x))
\]
and has to be Lipschitz and invertible for $\hat x$ to be well defined and converging. This is a difficult problem to solve (see \cite{KK98,AndrieuPraly06}), still open in general. This is why we have not presented this method.

\medskip

Other approaches consider non-smooth observers. In particular, the following observer proposed by Levant \cite{Levant}
\[
\left\{\begin{array}{lll}
\dot{\hat z}_1 & = & z_2 -K_1L^{\frac{1}{3}}|z_1-y|^{\frac{2}{3}}sign(z_1-y)\\
\dot{\hat z}_2 & = & z_3 -K_2L^{\frac{1}{2}}|z_1-y|^{\frac{1}{3}}sign(z_1-y)\\
\dot{\hat z}_3 & = & -K_3Lsign(z_1-y)
\end{array}\right.
\]
reconstructs theoretically $\dot y$ and $\ddot y$ in finite time from the measurement $y$, provided that the solutions of the original system are bounded. However, such observers and more generally sliding-mode observers \cite{Spurgeon}, are extremely sensitive to measurement noise. Although well employed for mechanical or electro-mechanical systems, we believe that there are not very well suited to epidemiological models.
\color{black}
\chapter{Practical and numerical considerations}
\label{sectionPractical}

\abstract*{The objective of this chapter is to show how to put theory in practice, illustrated in some cases studies.}

\abstract{The objective of this chapter is to show how to put theory in practice, illustrated in some cases studies.}

\section{Practical identifiability}

Till now we have studied structural identifiability/observability . While structural identifiability is a property of the  model structure, given a set of outputs, practical identifiability is related to the actual data. In particular, it depends on the amount of information contained in the data. 

\m

A model can be structurally identifiable, but still be practically unidentifiable due to poor data quality, e.g., bad signal-to-noise ratio, errors in measurement or sparse sampling \cite{Raue:2009aa}. Structural identifiability means that parameters are identifiable with ideal (continuous, noise-free) data. While structural identifiability is a prerequisite for parameter identification, it does not guarantee that parameters are practically identifiable with a finite number of noisy data points. 

\n
Moreover, parameter estimation requires using numerical optimization algorithms. The distance, for the problem considered, to the nearest ill-posed problem, \cite{MR895087,MR1041063}, i.e., the conditioning of the problem, can challenge the convergence of algorithms.

\n
Another source of practical  unidentifiability is the lack of information from the data, i.e., the signal from the data does not satisfy the persistence of excitation \cite{MR1261705}. This is the case when the observation  is  near an equilibrium \cite{MR2313504}.

\m
\n
In this section we use sensitivity analysis and results from asymptotic statistical theory to study practical identifiability. We refer to previous surveys and papers on the topic \cite{MR2547126,MR3024538,Banks2009,davidian1995,MR3203115}. Our purpose here is to give an intuitive account of these techniques.

\subsection{Rationale for using sensitivity analysis} Practical identifiability is often assessed in terms of confidence intervals on parameters \cite{Wieland2021}. Confidence intervals can be derived from the Fisher Information Matrix (FIM) \cite{bolker2008ecological}. More specifically, the covariance matrix ($\Sigma$) of the estimated parameters may be approximated as the inverse of the FIM. The diagonal elements of $\Sigma\approx \text{FIM}^{-1}$ correspond to the variance of the parameter estimates. Their square-roots (the standard deviations) give confidence intervals on the parameters, thus providing information on practical identifiability.

In the least-squares framework, the Fisher Information Matrix can be expressed in terms of sensitivity matrices, that we define below.

\subsection{Observed system} We consider that the initial condition $x_0$ is unknown. Unless otherwise specified, the term ``parameter'' now refers to both the parameter $\theta$ and the initial condition $x_0$, i.e. $\Theta=(\theta,x_0)$. We make explicit the dependence of the state variables $x$ and $y$ on $\Theta$ to clarify the following derivations:

\begin{equation}\label{sys10}
	\left \{
	\begin{array}{l}
		\dot x (t,\Theta)=  f(x(t,\Theta), \theta) , \quad x(0,\Theta)=x_0 , \\
		\dot \theta (t)=0, \quad \theta(0)=\theta , \\
		\\
		y (t,\Theta) =h(x(t,\Theta), \theta). \\
	\end{array}
	\right.
\end{equation}
with $x \in \R^n$, $y \in \R^m$ and $\theta \in \R^p$.  

\subsection{Sensitivity analysis} 

We wish to quantify how the observed variable $y(t,\Theta)$ changes for a small parameter variation $\Delta \Theta$.

We denote the Jacobian of the observation $y(t,\Theta)$ with respect to the parameter $\Theta$ as
\[
\chi(t, \Theta)=\dfrac{\partial y}{\partial \Theta} (t, \Theta)\,.
\]
This $m \times (n+p)$ matrix is called the sensitivity matrix.

\n
By linearization (first-order Taylor approximation), one can write
\[
\Delta y (t, \Theta) =  \chi(t, \Theta)\,\Delta \Theta.
\]

\paragraph{Side remarks} Reid \cite{Reid} defined a parameter vector as ``sensitivity identifiable'' if the above equation can be solved uniquely for $\Delta \Theta$. This linear problem is well known: if $\chi$ has maximal rank then the solution is given by means of the Moore-Penrose pseudo-inverse $\chi^+=(\chi^T\,\chi)^{-1}\chi^T$:
\[  
\Delta \Theta =\chi^+(t, \Theta) \,\Delta y(t, \Theta).
\]
It is also well known \cite{Golub1989} that the sensitivity of this solution is ruled by the condition number $\kappa_2(\chi)= \sigma_{\text{max}}/\sigma_{\text{min}}$ with $\sigma_{\text{max}}$ and $\sigma_{\text{min}}$ respectively the greatest and smallest singular value of $\chi$ (which are the  square roots of the corresponding eigenvalues of $\chi^T\chi$).

\subsection{Ordinary Least Squares}
Now we consider a set of $M$ observations $Y_i$, $i = 1,\ldots, M$, that have been obtained at times $t_i$. We assume that the observation is given by 
\[
Y_i= y(t_i, \Theta)+ \mathcal E_i ,
\]
with the error $\mathcal E_i$ assumed to be a random variable satisfying the following assumptions:

\begin{itemize}
	\item the errors $\mathcal E_i$ have mean zero $E[\mathcal E_i]=0$;
	\item the errors have a constant variance $\text{var}(\mathcal E_i)= \sigma^2$;
	\item the errors are independent and identically distributed. 
\end{itemize}

\n
The Fisher Information Matrix, for the preceding defined observations, is defined as 
\begin{equation}\label{FIM}
	\text{FIM} (\Theta,\sigma^2) = \frac{1}{\sigma^2}\sum_{i=1}^{M} \; \chi(t_i, \Theta)^T \; \chi(t_i, \Theta). 
\end{equation}

\n
Solving the ordinary least square (OLS) equations gives an estimator $\hat \Theta_{\text{OLS}}$ of the parameter $\Theta$:
\begin{equation}\label{estimtheta}
	\hat \Theta_{\text{OLS}} = \arg \min_{ \Theta}  \sum_{i=1}^{M}  \left( Y_i - y(t_i, \Theta) \right)^2 
\end{equation}

\n
Even though the error's distribution is not specified,  asymptotic statistical theory can be used 
to approximate the mean and variance of the estimated $\Theta$ (a random variable) \cite{MR986070,MR3203115}: the bias-adjusted approximation for $\sigma^2$ (with $n+p$ ``parameters'') is
\begin{equation}\label{estimbias}
	\hat \sigma_{\text{OLS}}^2= \dfrac{1}{M-(n+p)} \; \sum_{i=1}^{M} \;  \left( Y_i - y\left(t_i,\hat \Theta_{\text{OLS}}\right) \right)^2. 
\end{equation}

\subsection{Confidence intervals} 

The above approximation of the error variance can be used to further approximate the parameter covariance matrix $\Sigma$:
\begin{equation}\label{estimcovmat}
	\hat\Sigma := \left[\text{FIM}\left(\hat \Theta, \hat \sigma_{\text{OLS}}^2\right)\right]^{-1}.
\end{equation}

\n
The standard error (SE) for $\hat \Theta_{\text{OLS}}$ can be approximated by taking the square roots of the diagonal elements of the covariance matrix $\Sigma$: for all $k=1,\ldots,n+p$,
\begin{equation}\label{SE}
	\text{SE} \left( \hat \Theta_{\text{OLS}}(k)\right) = \sqrt{ \hat\Sigma_{kk}}\,.
\end{equation}

\n
Finally, to compute the 95\% confidence interval for the $k$-th component of the parameter vector $\hat \Theta_{\text{OLS}}$ with $n+p$ ``parameters'', one may use the Student's $\mathsf{t}-$distribution with $M-(n+p)$ degrees of freedom: letting
\[
\zeta(k)=\mathsf{t}_{0.025}^{M-(n+p)} \times \text{SE} \left( \hat \Theta_{\text{OLS}}(k)\right)\,,
\]
the confidence interval is defined as
\[  
\hat \Theta_{\text{OLS}}(k) - \zeta(k) <   \hat \Theta_{\text{OLS}}(k) <  \hat \Theta_{\text{OLS}}(k) + \zeta(k) \,.
\]

From  these formulas it appears  that the conditioning of the Fisher Information Matrix plays an essential role. Huge confidence intervals give indications about the practicality of the identification.

\subsection{Computing the sensitivity matrix}

The sensitivity matrix $\chi(t,\Theta)$, with $\Theta=(\theta,x_0)$, is obtained by integrating an ODE. The components of the ODE to be integrated depend on whether one differentiates with respect to $\theta$ or $x_0$.

\paragraph{Differentiating with respect to $\theta$} The first part of the ODE is given by 

\[\dfrac{\partial y }{\partial \theta} (t, \theta,x_0) =\dfrac{\mathrm{d} }{\mathrm{d} \theta} \;h(x(t,\theta,x_0),\theta)= \dfrac{\partial h  }{\partial x}  \; \dfrac{\partial x }{\partial \theta} (t,x_0,\theta) +\dfrac{\partial h }{\partial \theta}(x(t,\theta,x_0),\theta).\]

\n
The Jacobian $ \dfrac{\partial h  }{\partial x} $ is a $m \times n $ matrix while $ \dfrac{\partial h  }{\partial \theta} $ is a $m \times p $ matrix. 

We then have to compute the $ n \times p$ matrix \[
z(t,\Theta):=\dfrac{\partial x }{\partial \theta} (t,\theta,x_0)=\dfrac{\partial x }{\partial \theta} (t,\Theta)\,.
\]
Let $A(t,\Theta)$ and $B(t,\Theta)$ be the following time-dependent $n \times n$ and $n \times p$ matrices, respectively: 
\[
A(t,\Theta) :=\dfrac{\partial f }{\partial x }( x(t,\Theta), \theta)\,,
\] 
and
\[
B(t,\Theta) :=  \dfrac{\partial f }{\partial \theta }( x(t,\Theta), \theta)\,.
\]
It is well known  \cite{MR1929104} that $z(t)$ is the solution of the linear matrix equation:
\[
\dot z (t,\Theta)= A(t,\Theta) \, z(t,\Theta) +B(t,\Theta), 
\]
with the initial condition $z(0,\Theta)=0_{n\times p}$ (a zero matrix of size $n\times p$). 

\paragraph{Differentiating with respect to $x_0$} The second part of the ODE is given by
\[\dfrac{\partial y }{\partial x_0} (t, \theta,x_0) =\dfrac{\partial  }{\partial x_0} \;h(x(t,\theta,x_0),\theta)= \dfrac{\partial h  }{\partial x}  \; \dfrac{\partial x }{\partial x_0} (t,\theta,x_0). \]

Let
\[
w(t,\Theta) := \dfrac{\partial x }{\partial x_0} (t,\theta, x_0)=\dfrac{\partial x }{\partial x_0} (t,\Theta)\,.
\]
\n
Based on the same reference \cite{MR1929104}, $w(t,\Theta)$ is solution of the linear matrix ODE
\[
\dot w(t,\Theta)=A(t,\Theta)\, w(t,\Theta),
\]
with the initial condition $w(0,\Theta) = \text{Id}_{n\times n}$ (the identity matrix of size $n$).

\paragraph{Full system} To summarize, one has to solve the following system in dimension $n^2+np+n$
\begin{equation}\label{eq:fullSystem}
	\left \{
	\begin{array}{rrl}
		\dot x (t,\Theta)&=&  f(x(t,\Theta), \theta), \quad x(0,\Theta)=x_0 , \\[2mm]
		\dot z(t,\Theta) &=&A(t,\Theta)\,z(t,\Theta)+B(t,\Theta) , \quad z(0,\Theta)=0_{n \times p} ,\\[2mm]
		\dot w (t,\Theta) &=&A(t,\Theta)\,w(t,\Theta) , \quad w(0,\Theta)= \text{Id}_{n \times n}
	\end{array}
	\right.
\end{equation}
with $A(t,\Theta) =\dfrac{\partial f }{\partial x }( x(t,\Theta),\theta ) $   and  $B(t,\Theta)=  \dfrac{\partial f }{\partial \theta }( x(t,\Theta),\theta )$ .

\m
\n
For large systems, the computation of the different Jacobians can be prohibitive, in this case  automatic differentiation software has to be used.

\subsection{Some case studies}

In this section, we consider two classical examples as case studies. These examples have been used in many books of mathematical epidemiology, e.g. \cite{Murr2002}.

\paragraph{Case 1. Influenza in a boarding school}

\smallskip

Our first example is an outbreak of influenza in a United Kingdom boarding school which occurred in 1978 \cite{BMJ-boarding}. In \cite{Murr2002} the parameters $\beta, \gamma$ are identified by an unspecified ``best-fit'' algorithm. A more complete analysis is done in \cite{MR2002k:92001}, 
where the analysis is done using sensitivity analysis and asymptotic statistical theory. In \cite{MR3409181} the same example is considered. Different sources exist for the data  \cite{MR3969982,MR2002k:92001,MR2242784} with small differences.  

Using the figure in \cite{BMJ-boarding} and the \textit{Plot Digitizer} software, we got an approximation of the data. It was reported that $N = 763$, and the conditions at the start of
this outbreak were $S_0 = 762$ and $I_0 = 1$. We used the following data, in which time $t$ is in day and $I(t)$ denotes the number of infectious people at time $t$.

\begin{center}
	\begin{tabular}{|c||c|c|c|c|c|c|c|c|c|c|c|c|c|c|}
		\hline $\vphantom{\displaystyle \int}$
		$t$ & 0 & 1 & 2 & 3 & 4 & 5 & 6 & 7 & 8 & 9 & 10 & 11 & 12 & 13\\
		\hline $\vphantom{\displaystyle \int}$
		$I(t)$ & \, \verb"1" \, & \, \verb"6" \, & \, \verb"26" \, & \, \verb"73" \, & \, \verb"222" \, & \, \verb"293" \, & \, \verb"258" \, & \, \verb"237" \, & \, \verb"191" \, & \, \verb"124" \, & \, \verb"68" \, & \, \verb"26" \, & \, \verb"10" \, & \, \verb"3" \, \\
		\hline
	\end{tabular}	
\end{center}

\medskip

\m

Specifically, we considered model \eqref{KmcK} with $k=1$ (all infectious are assumed to be observed). We obtained the following OLS estimation, as given by the {\tt Scilab} software:
\[
\beta \approx 1.96,  \quad \gamma \approx 0.475
\]
(see the numerical code in Appendix \ref{appendixCode}). The fit is shown in Figure \ref{figBoardingSchool}.
\begin{figure}[!h] \centering 
	\includegraphics[scale=0.4]{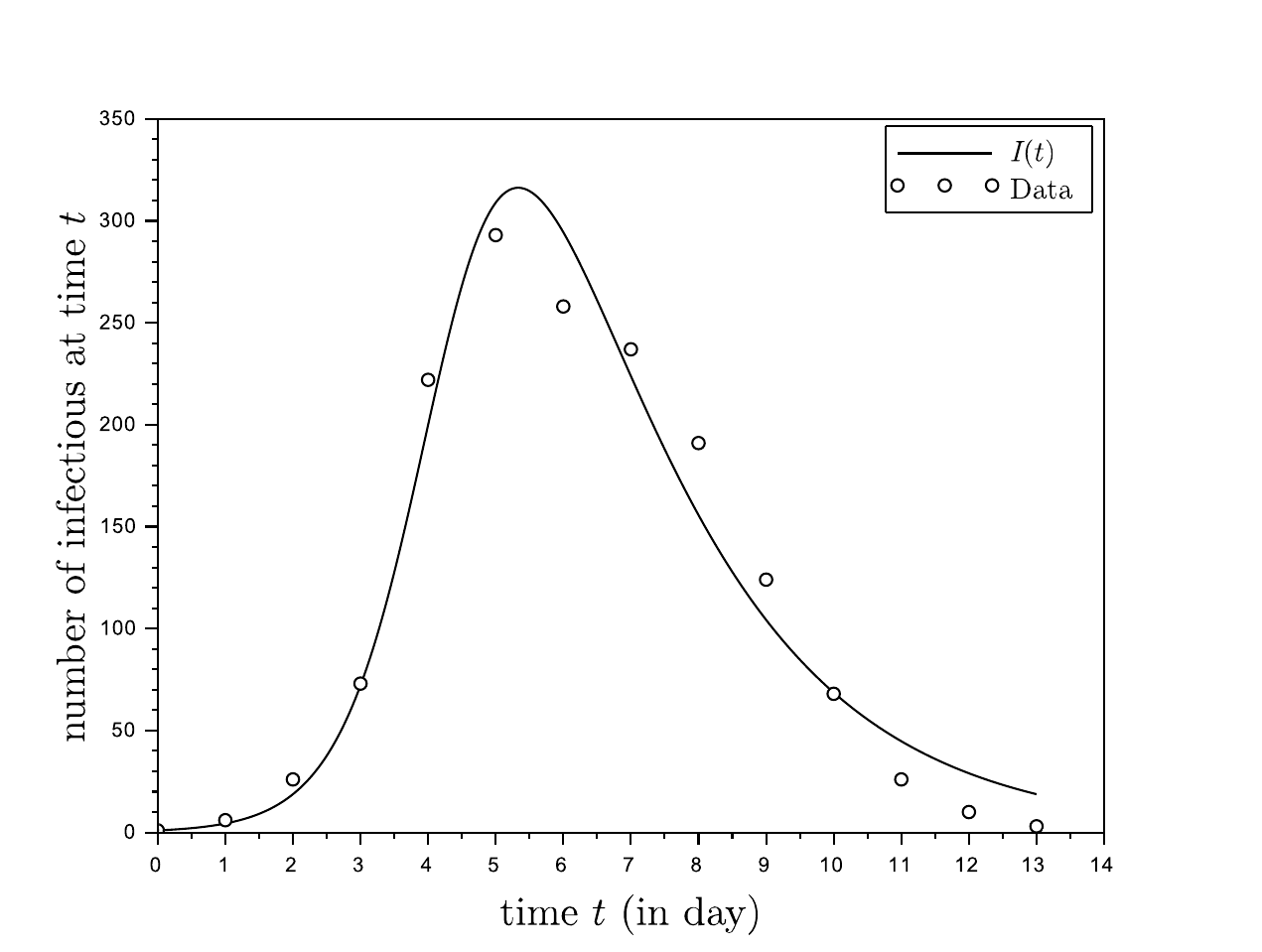}
	\caption{Boarding School example.}
	\label{figBoardingSchool}
\end{figure}

\medskip

\n
We computed the 95\% confidence intervals  using the formulas given in the preceding section: \eqref{FIM}, \eqref{estimtheta}, \eqref{estimbias}, \eqref{estimcovmat}, \eqref{SE}, up to a few changes due to the fact that the initial conditions $x_0=(S_0,I_0)$ are assumed to be known in this example (see Appendix \ref{appendixCode}). We find:
\begin{equation*}
	\begin{array}{l}
		\beta  \approx 1.96 \pm  0.073\,, 
		\\[2mm]
		\gamma  \approx 0.475 \pm  0.04\,. 
	\end{array}
\end{equation*}
\color{black}

One can obtain approximately the same results using the likelihood profile method to compute confidence intervals \cite{bolker2008ecological}, assuming normally distributed errors. However, it is well known that the profile method quickly becomes impractical for model with more than two parameters \cite{bolker2008ecological}, which is the rule rather than the exception, as will be the case in the following example. This is why we stick to the FIM method. Lastly, we note that the condition number of the FIM is approximately equal to 3.80 in this example.

\medskip

\paragraph{Case 2. Plague in Bombay}

Our second example is the Bombay Plague of 1905--1906 \cite{KmK1927}. We collected the data from \cite[Table IX]{PlagueBombay}, over the same period as \cite{KmK1927} (Dec. 17 to Jul. 21). The form of the data is presented in the following table, in which time $t$ is in week, and $\dot R(t)$ denotes the number of death per week at time $t$.

\medskip

\begin{center}
	\begin{tabular}{|c||c|c|c|c|c|c|c|c|c|c|c|c|c|c|c|}
		\hline $\vphantom{\displaystyle \int}$
		$t$ & 0 & 1 & 2 & 3 & 4 & 5 & 6 & $\cdots$ & 24 & 25 & 26 & 27 & 28 & 29 & 30\\
		\hline $\vphantom{\displaystyle \int}$
		$\dot R(t)$ & \, \verb"8" \, & \, \verb"10" \, & \, \verb"12" \, & \, \verb"16" \, & \, \verb"24" \, & \, \verb"48" \, & \, \verb"51" \, & \, $\cdots$ \, & \, \verb"106" \, & \, \verb"64" \,  & \, \verb"46" \, & \, \verb"35" \, & \, \verb"27" \, & \, \verb"28" \, & \, \verb"24" \, \\
		\hline
	\end{tabular}	
\end{center}

\medskip

\begin{figure}[!h] \centering 
	\includegraphics[scale=0.4]{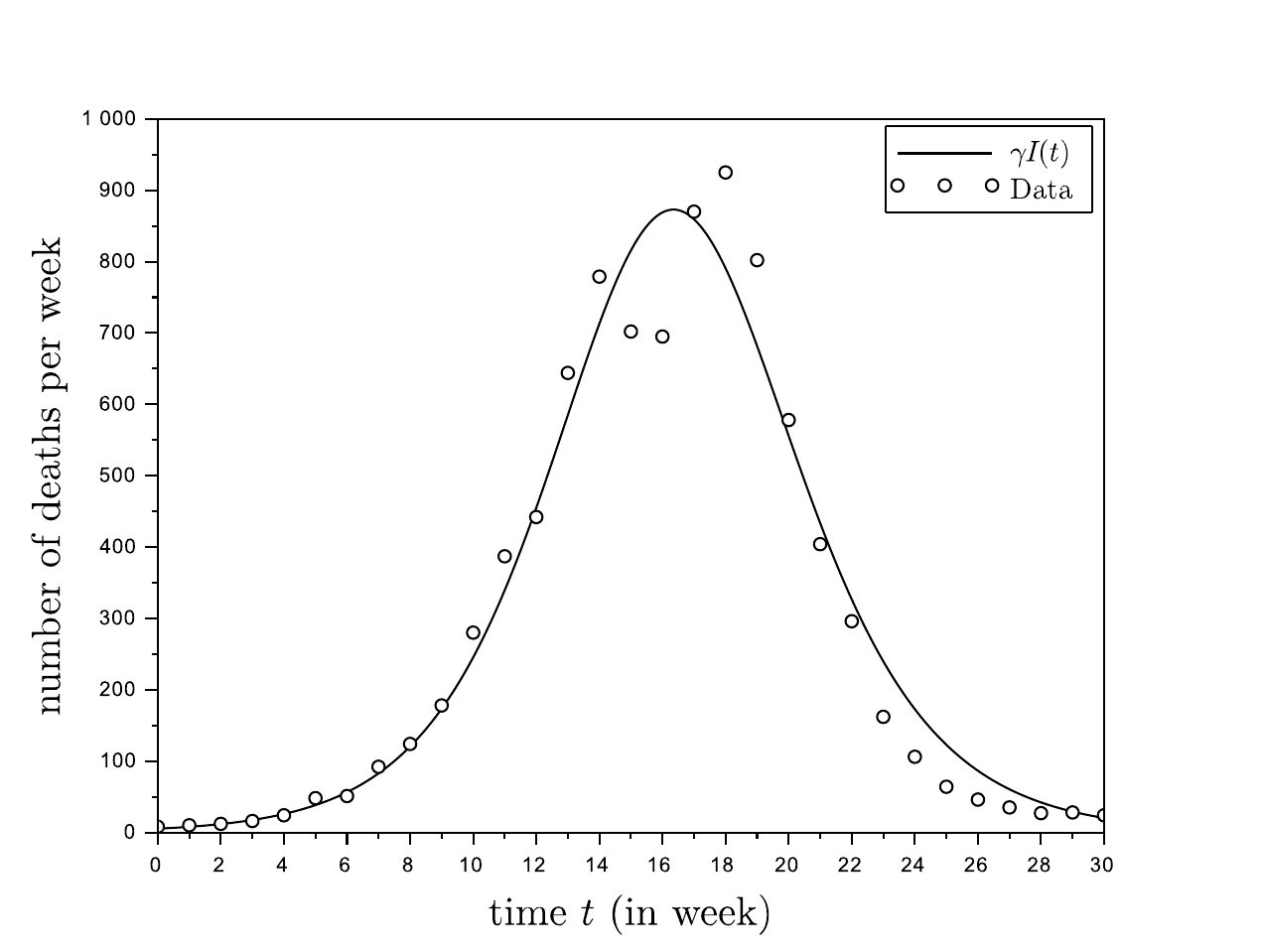}
	\caption{Plague in Bombay example.}
	\label{figBombay}
\end{figure}

We consider that the number of death per week is the same as $\dot R(t)=\gamma I(t)$, meaning that all infections lead to death, which is a reasonable assumption in this context \cite{MR2886018}. Therefore, we consider model \eqref{KmcK} with $k=\gamma$. In this example, not only the parameters $\beta$ and $\gamma$, but also the size of the population, $N$, as well as the initial conditions, $S_0$ and $I_0$, are unknown \cite{MR2886018}. According to Theorem \ref{th:sir} (and Remark \ref{rk:sir}), the model is neither observable nor identifiable.  However, the model is partly identifiable in the sense that $S_0$, $I_0$, $\gamma$ and $\tilde\beta=\beta/N$ are structurally identifiable. Starting from an arbitrary initial guess, we obtained the following OLS estimation, as given by the {\tt Scilab} software:
\[
\tilde\beta \approx 0.0000855\,,   
\quad 
\gamma \approx 3.72\,, 
\quad
S_0 \approx 4.81\,10^4\,, 
\quad
I_0 \approx 1.42\,,
\]
(see the numerical code in Appendix \ref{appendixCodeBombay}). The fit is shown in Figure \ref{figBombay}.

Proceeding as in the previous example, we obtained the following 95\%-confidence intervals:
\begin{equation*}
	\begin{array}{l}
		\tilde\beta  \approx 0.0000855 \pm  0.00157\,,
		\\[2mm]
		\gamma  \approx 3.72 \pm  25\,,
		\\[2mm]
		S_0 \approx 4.81\, 10^4 \pm  6\,10^4\,, 
		\\[2mm]
		I_0 \approx 1.42\pm  34\,.
	\end{array}
\end{equation*}

The confidence intervals are huge, which means that we can have absolutely no confidence in the estimated values of the parameters, even though the fit looks good and these parameters are structurally identifiable in principle. In practice, one can show that many other and very different combinations of the parameters can yield approximately the same fit. Note that if we did as if the initial conditions were known, the confidence intervals on $\tilde\beta$ and $\gamma$ would be reasonable, as in the previous example. In this example, the condition number of the FIM is approximately $9.14\times 10^{24}$, meaning that the problem is ``sloppy'' \cite{Chis:2016}; note however that while it is usual that a model is both sloppy and practically non-identifiable, this is not always the case \cite{Cole:2020}. Altogether, we can conclude that there is a severe practical identifiability issue in this classical example.

\subsection{Discussion}

\n
The SIR model of Kermack-McKendrick has been studied in  a series of papers 
\cite{CRSC-tr08,MR2547126,MR2532016,Capaldi:2012aa} where the problem of observability/identifiability is approached from the statistical point view:  addressing  parameter identifiability by exploiting properties of both the sensitivity matrix and uncertainty quantification in the form of standard errors. In this series of papers, structural observability and identifiability were not explicitly addressed. For example in
\cite{MR2532016} the authors identify $(S_0,I_0,\beta/N, \gamma)$ based on incidence observations - akin to equation (\ref{cumul}) with $k=1$ - which we have proved to be structurally identifiable (see Theorem \ref{thm:cumulativeIncidence}). Similarly in \cite{Capaldi:2012aa} the authors seek to identify $(S_0, I_0, \beta, \gamma)$ which, with $N$ known, are structurally identifiable in principle. However the authors encounter practical identifiability issues. This is a typical example of strictly practical unidentifiability (as in our second example, the Plague in Bombay).

\m
\n
Although a structural observability and identifiability analysis should be done as a prerequisite to a practical identifiability analysis, it does not suffice. Moreover, when doing practical identifiability analyses, the error structure of the data should be considered. For instance, sensitivity analyses can be extended to non-constant error variance through Generalized Least Squares (GLS), which makes it possible to test different ways of weighting errors \cite{CRSC-tr08}. Appendix \ref{app:GLS} shows the principle of GLS and how to compute the Fisher Information Matrix in this case. We also provide an online example\footnote{https://github.com/nikcunniffe/Identifiability} based on the ``Influenza in a Boarding School'' data (since the ``Plague in Bombay'' example generates numerical optimization issues related to its practical non-identifiability).

\n
An additional issue may occur when the output signal is not sufficiently informative (i.e., not {\em persistently exciting} \cite{MR1261705}). For example when the data correspond to states near unobservability, e.g., near an equilibrium. In those cases, one has to wait to have data sufficiently far from equilibrium.

To conclude, the problem of observability and identifiability, either structural or practical, is far from being simple, even in relatively simple SIR models with seemingly good quality data \cite{CRSC-tr08}. Of course, the more complex the model, the more parameters there are to identify, the more serious the problem of identifiability.

\section{Observers in practice}
\label{obspratic}

In this section, we show how the various observers presented in Section \ref{secobservers} behave in practice, and the role of the tuning parameters. Up to know, we have assumed the measurements to be perfect i.e.~not tainted with any noise.
Since integration has good ``averaging'' properties, an observer is expected to filter noise or inaccuracies in the measurements. However, we will see that the filtering capacity of an observer is related to his convergence speed, which often leads to a ``precision-speed'' dilemma in the choice of the observer or his settings.

\smallskip

Let us underline that when identifiability/observability cannot be proved theoretically or is too difficult to be proven analytically, one can still look for an observer and study its asymptotic convergence, theoretically or numerically.

\subsection{Observers with linear assignable error dynamics}
\label{secobslinearpractic}

We illustrate the observer \eqref{obs3stades} of the age-structured model \eqref{model3stages} on simulations, for the following values of the parameters.

\begin{center}
	\begin{tabular}{|c|c|c|c|c|c|c|c|}
		\hline $\vphantom{\displaystyle \int}$
		$a_1$ & $a_2$ & $m_1$ & $m_2$ & $m_3$ & $k$ & $\bar r_{min}$ & $\bar r_{max}$\\
		\hline\hline
		$\vphantom{\displaystyle \int} \; $\verb"0.1" \; & \; \verb"0.1" \; & \; \verb"0.05" \; & \; \verb"0.07" \; & \; \verb"0.07" \; & \; \verb"1.0" \; & \; \verb"0.9" \; & \; \verb"1.1" \; \\
		\hline
	\end{tabular}
\end{center}\m
The following code is used to compute the gain vector \verb"G" for a set of desired eigenvalues.

\medskip

\begin{verbatim}
	a1=0.1;a2=0.1;m1=0.05;m2=0.07;m3=0.07;
	Sp=0.3*[-0.3,-0.33,-0.36];
	A=[-a1-m1,0,0;  
	a1,-a2-m2,0;
	0,a2,-m3];
	C=[0,0,1];
	B=[0;0;1];
	Obs=[C;C*A;C*A*A];
	L=inv(Obs)*B;
	P=[L,A*L,A*A*L];
	Abar=inv(P)*A*P;
	sigma=coeff(poly(Sp,'x'));
	G=P*(-sigma(1:3)'-Abar(:,3));
\end{verbatim}

\medskip

\begin{figure}[ht!] \centering 
	\includegraphics[scale=0.4]{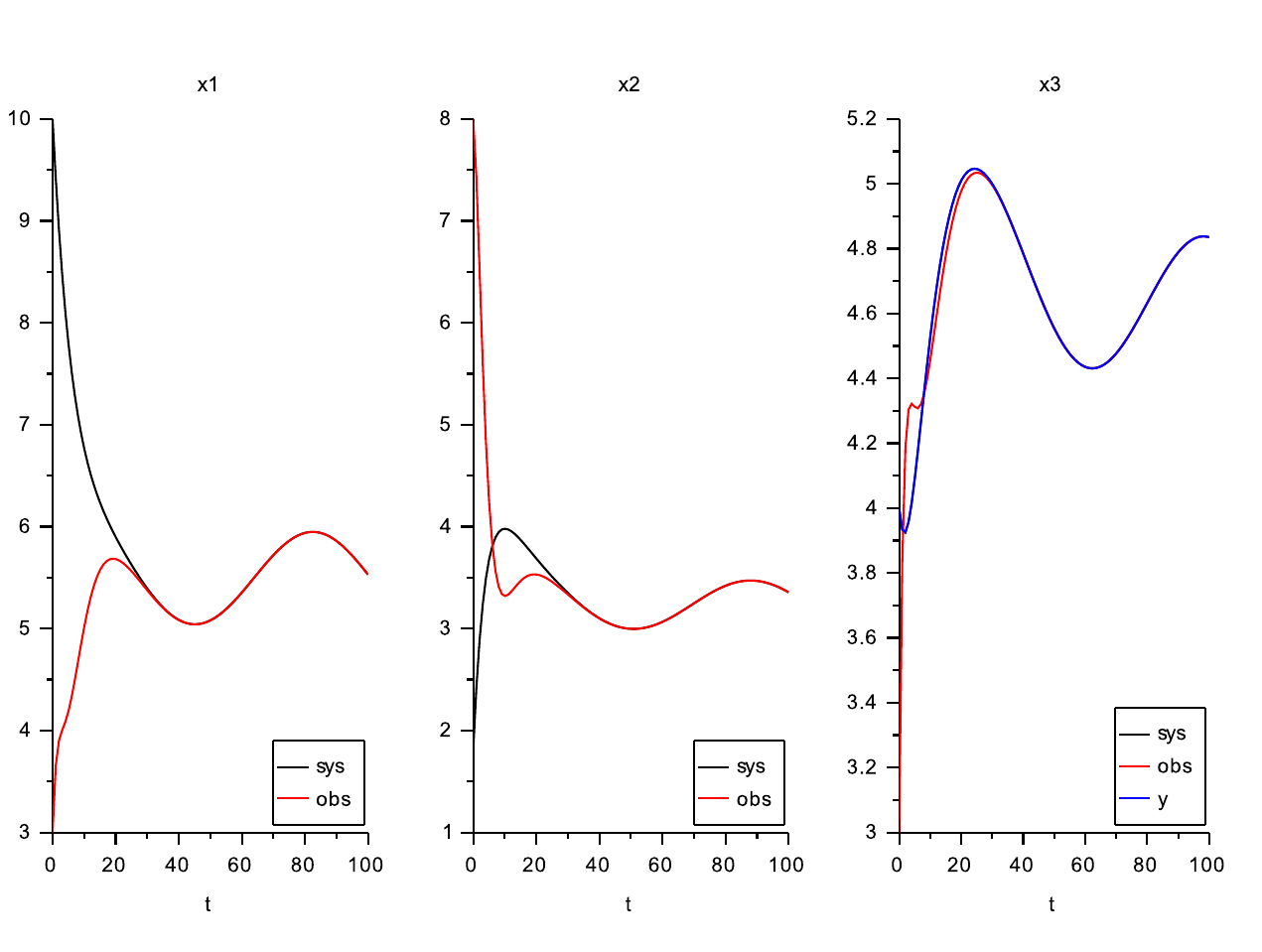}
	\caption{$Sp(A+GC)=\{-0.3,-0.33,-0.36\}$ without measurement noise}
	\label{simu1}
\end{figure}

\begin{figure}[!ht] \centering 
	\includegraphics[scale=0.4]{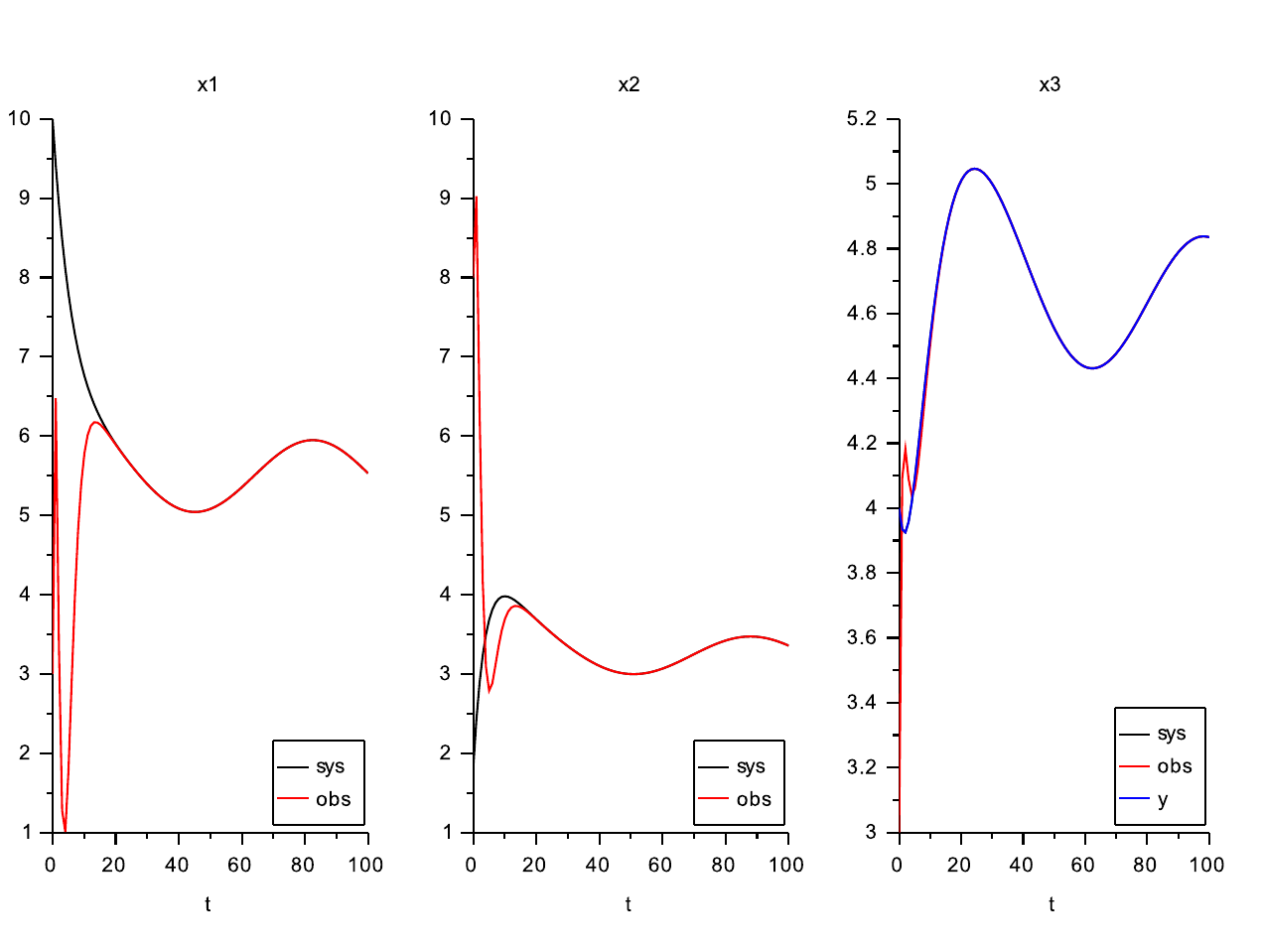} 
	\caption{$Sp(A+GC)=\{-0.6,-0.66,-0.72\}$ without measurement noise}
	\label{simu2}
\end{figure}

\begin{figure}[!ht] \centering 
	\includegraphics[scale=0.4]{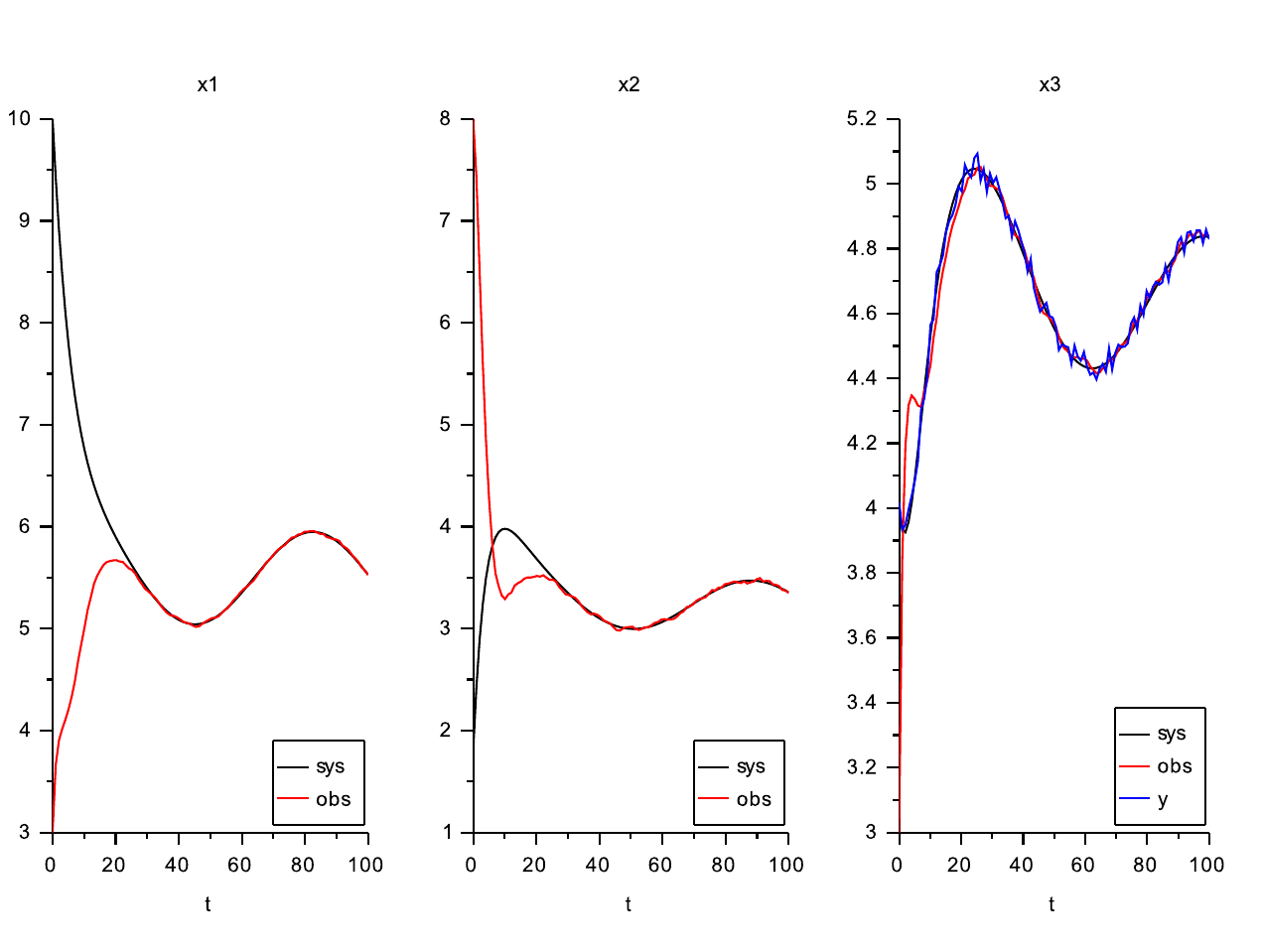} 
	\caption{$Sp(A+GC)=\{-0.3,-0.33,-0.36\}$  with measurement noise}
	\label{simu1noise}
\end{figure}

\begin{figure}[!ht] \centering 
	\includegraphics[scale=0.4]{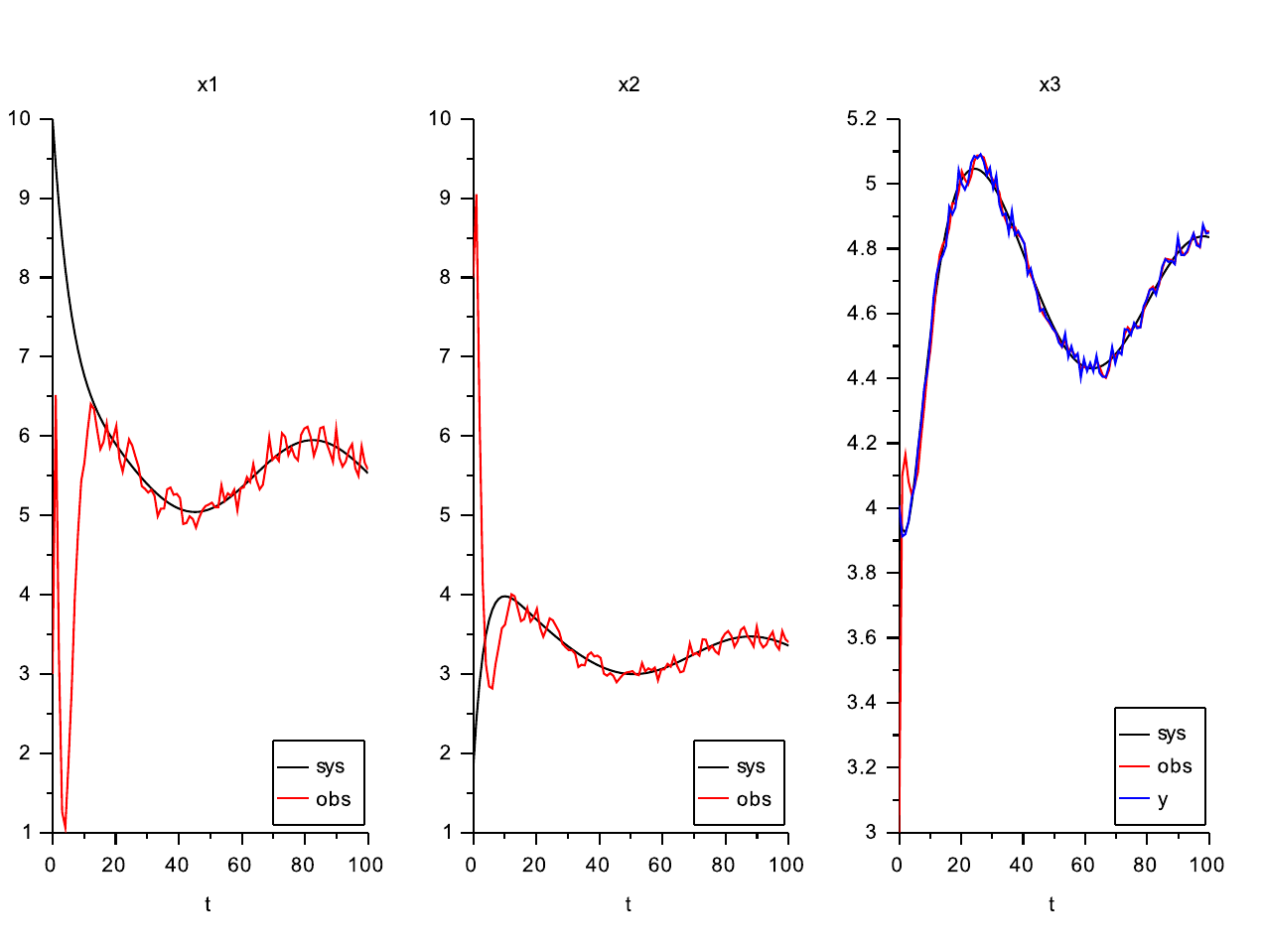} 
	\caption{$Sp(A+GC)=\{-0.6,-0.66,-0.72\}$ with measurement noise}
	\label{simu2noise}
\end{figure}

Figure \ref{simu1} shows convergence for a moderately negative spectrum, while Figure \ref{simu2} shows the acceleration of convergence obtained for a spectrum located further to the left in the complex plane. For the same choice of gains, Figures \ref{simu1noise} and \ref{simu2noise} show the effect of noise on the $y(\cdot)$ measurements. It can be seen that a faster convergence is more sensitive to noise and loses accuracy. In practice, one often has to make a compromise for the choice of the observer's setting.

\clearpage

\color{blue}
\subsection{About observers with Lipchitz non-linearity}
\label{obsLipnum}

Let us illustrate the result on Exemple \ref{exSIRS} with the following numerical values of the parameters
\begin{center}
    \begin{tabular}{|c|c|c|}
    \hline
    $\beta$ & $\gamma$ & $\mu$\\
    \hline\hline $\vphantom{\displaystyle \int}$
    $\vphantom{\displaystyle \int}$ \; \verb"0.13" \; & \; \verb"0.1" \; & \; \verb"0.05" \;\\
    \hline
    \end{tabular}
\end{center}
The gain vector $G$ such that $Sp(A+GC)=\{-1,-2\}$, and the corresponding matrix $P$ solution of the Lyapunov equation \eqref{Lyap_A+GC} can be numerically computed 
(using, for instance, scilab) as
\begin{verbatim}
    G  = | -23.889 |            P  = | 73.135667   83.710112 |
         | 20.909  |                 | 83.710112   96.199374 |
\end{verbatim}
with 
\begin{verbatim}
   norm(P) = 169.16821
\end{verbatim}
for which the condition \eqref{condobsLip} gives $\varepsilon \leq 0.001477$.
As an illustrattion, for a total population size of one billion, this gives $I \leq 1477$, which is a very small number... 

\medskip

This example shows that this technique is not well suited to epidemiological models such as the SIRS one, because it requires a too small Lipschitz constant of the non linear terms. However, we consider useful to have exposed this known approach and shown its drawback.
\color{black}

\subsection{Observers with asymptotic convergence}
\label{obssirrobustnum}

We illustrate on simulations the behavior of the asymptotic observer \eqref{obsSIRrobust} of the SIR model with fluctuating rates \eqref{SIRfluctuating}, for the following values of the parameters.

\smallskip

\begin{figure}[!h] \centering 
	\includegraphics[scale=0.4]{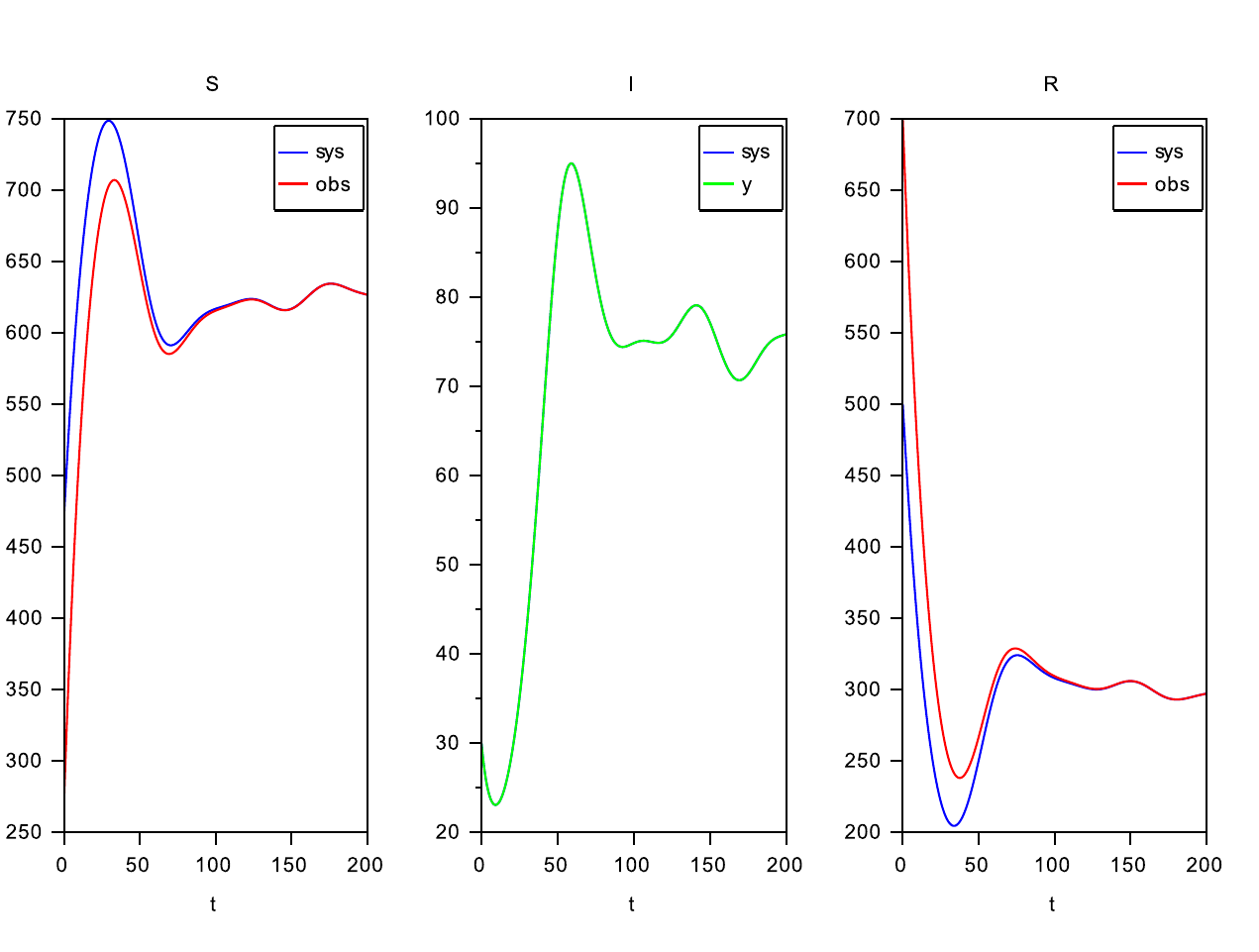} 
	\caption{The high gain observer without measurement noise}
	\label{simuSIRnat}
\end{figure}

\begin{figure}[!h] \centering 
	\includegraphics[scale=0.4]{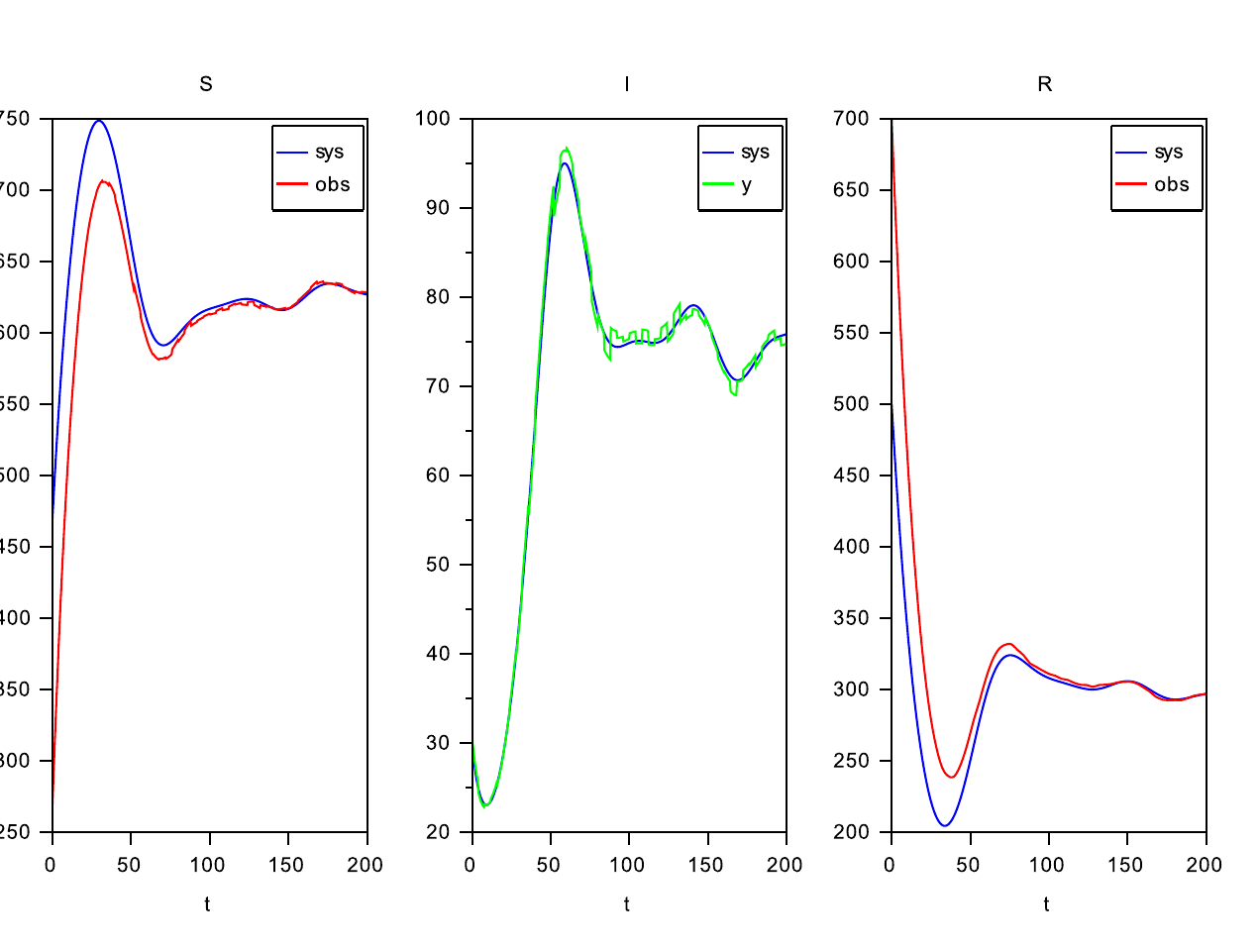} 
	\caption{The high gain with measurement noise}
	\label{simuSIRnatNoise}
\end{figure}

\begin{center}
	\begin{tabular}{|c|c|c|c|c|}
		\hline $\vphantom{\displaystyle \int}$
		$\beta$ & $\rho$ & $\nu$ & $\mu$ & $N$\\
		\hline\hline
	$\vphantom{\displaystyle \int}$	\; \verb"0.4"$\pm$\verb"0.08" \; & \; \verb"0.2"$\pm$\verb"0.04" \; & \; \verb"0.05" \; & \; \verb"0.05" \; & \; \verb"1000" \; \\
		\hline
	\end{tabular}
\end{center}
\smallskip
Here $\beta$ and $\rho$ are functions of time chosen randomly in between the bounds given in the table.
Figures \ref{simuSIRnat} and \ref{simuSIRnatNoise} show that the observer has a convergence relatively insensitive to measurement noise, but the speed of convergence is slow because the exponential decay of the error is equal to $\mu$, which is not adjustable.

Unlike the observers in previous sections, let us underline that the present observer is not based on a $\hat y-y$ innovation. Therefore, one is not informed of the quality of the estimate over time, which is a {\em price to pay} to have a observer insensitive to unknown variations of the epidemic parameters $\beta$, $\rho$.

\subsection{Observers with partially assignable error dynamics}
\label{obsmalarianum}

\begin{center}
	\begin{figure}[h!]
		\centering
		\includegraphics[width=0.9\textwidth]{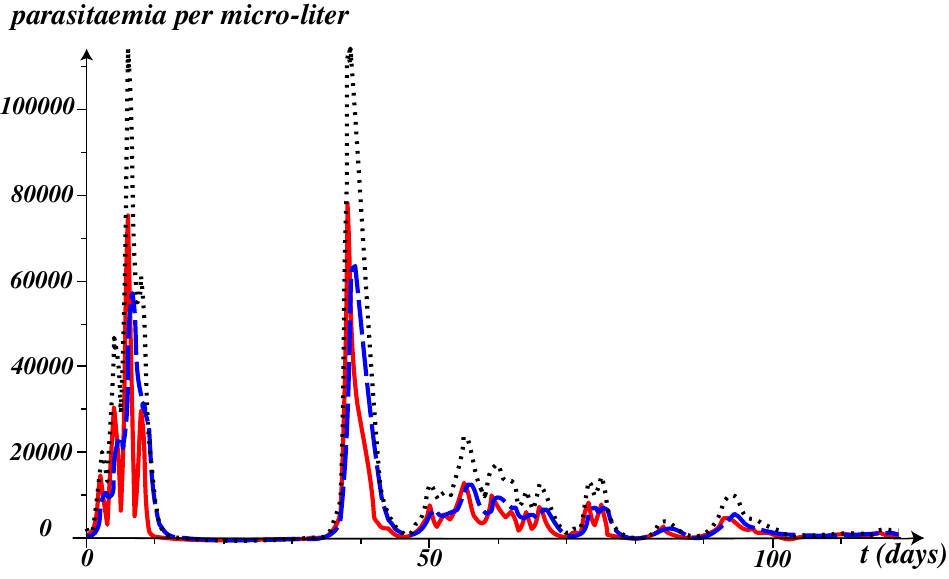}
		\caption{Example of patient S1204:  Measures (data) of peripheral parasitaemia are plotted with red solid line,  The estimations delivered by the observer~\eqref{AMGobs} are plotted with blue dashed line for the estimated  sequestered parasitaemia, and with black dotted line for the estimated total parasitaemia.
			The gain used is $L=(0, 5, 5, 0, 0, 0, 0)^T$.
			$Sp(\bar A-LC)=\{-10, -\mu_S, -\mu_M, -\mu_5-\gamma_5, -\mu_4-\gamma_4, -\mu_3-\gamma_3, -\mu_2-\gamma_2-\gamma_1\}$.}
		\label{fig:PatientS1204}
	\end{figure}
\end{center}

The observer given in  \eqref{AMGobs} for the intra-host malaria model \eqref{AMGreduit} is illustrated here on real data, as one can see for instance on Figure \ref{fig:PatientS1204}.

As already mentioned in Remark \ref{remobsnonassis}, the speed of convergence of this observer cannot be tuned as fast as desired. However, this is quite satisfactory in practice. Let us also underline that the observer does not require the reconstruction of the parameter $\beta$, although this parameter is identifiable (see Section \ref{secidentchgvar}). This is a strength of this observer, because the parameter $\beta$ could switch or fluctuate with time.

\subsection{High gain observer}

The non-linear observer \eqref{obsSIRhighgain}-\eqref{obsSIRhighgain} of the classical SIR model \eqref{SIRclassical} is illustrated on simulations for the following values

\smallskip

\begin{figure}[!h] \centering 
	\includegraphics [width=\linewidth]{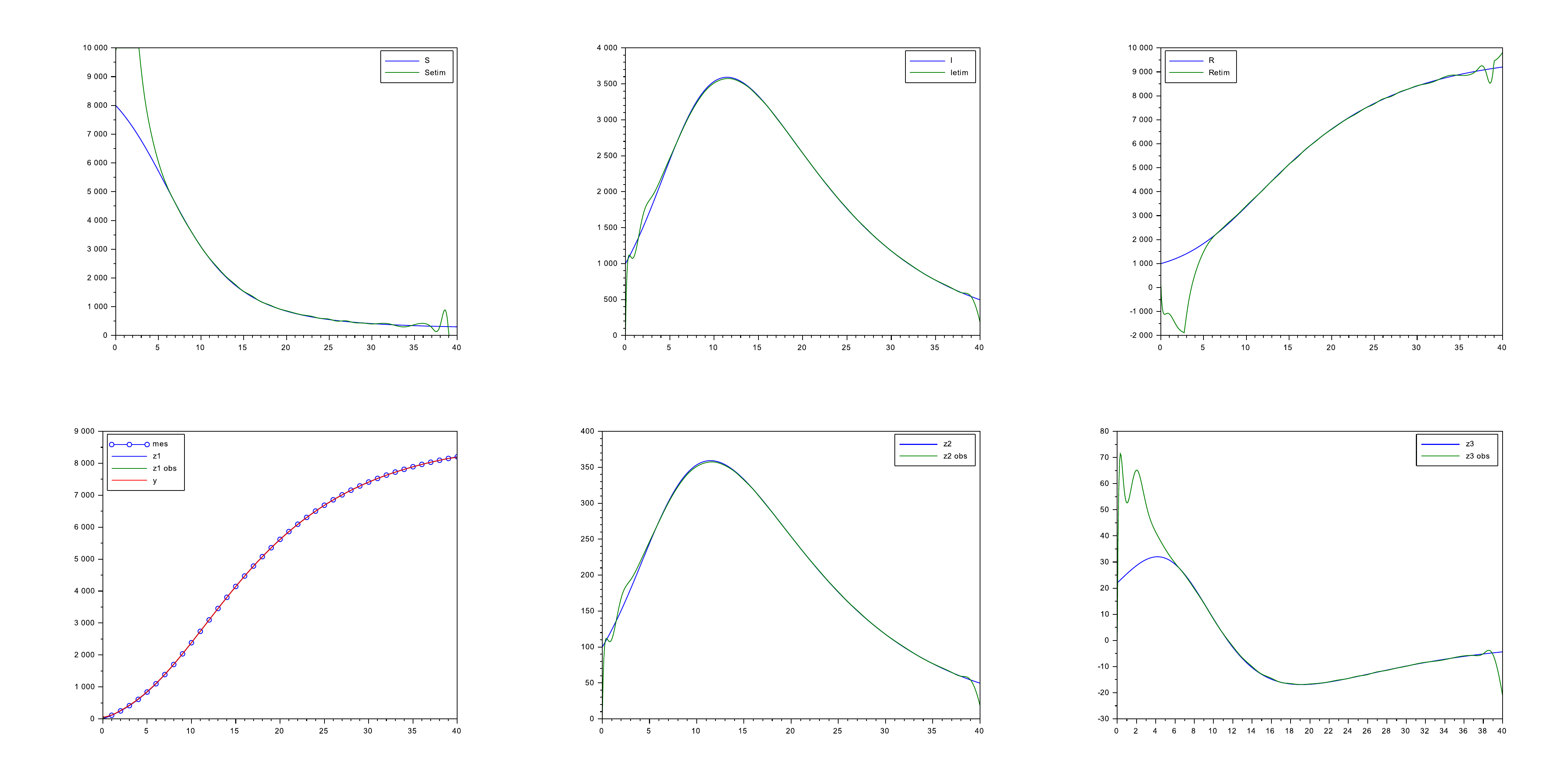} 
	\caption{Observer simulations (top: variables $S$, $I$, $R$ and their estimates; bottom:  coordinates $z_i$ with measurements points on the left)}.
	\label{simuSIR}
\end{figure}

\begin{figure}[!h] \centering 
	\includegraphics [width=\linewidth]{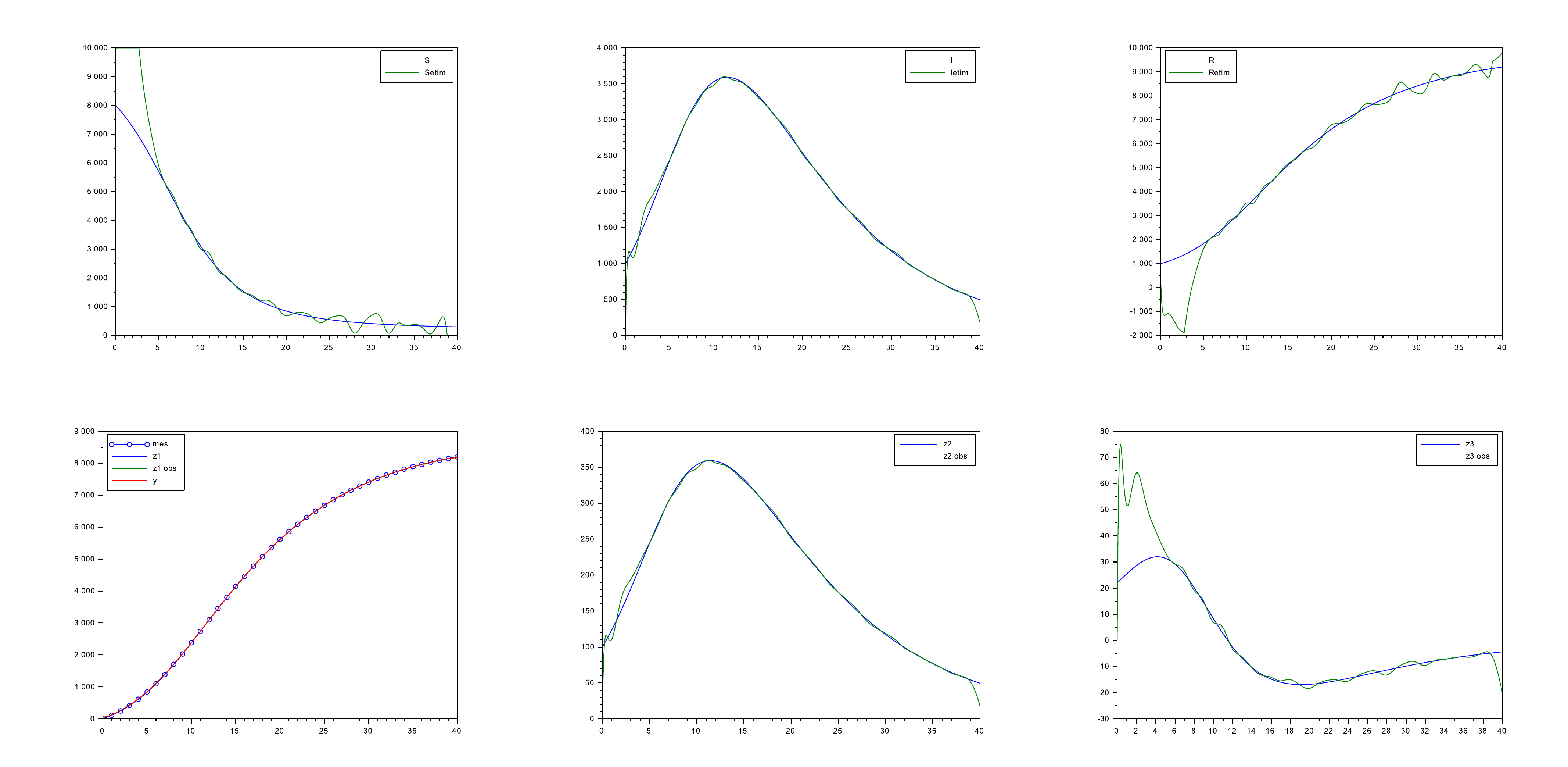} 
	\caption{Observer simulations with noisy measurements}
	\label{simuSIRnoise}
\end{figure}

\begin{center}
	\begin{tabular}{|c|c|c|c|}
		\hline $\vphantom{\displaystyle \int}$
		$\beta$ & $\rho$ & $N$\\
		\hline\hline
		 $\vphantom{\displaystyle \int}$ \; \verb"0.4" \;  & \; \verb"0.1" \; & \; \verb"10000" \;\\
		\hline
	\end{tabular}
\end{center}
\smallskip
where the $y$ cumulative measures were made discretely every day (rounded to the nearest integer). In order to obtain a time-continuous $y(\cdot)$ signal, we performed an interpolation by cubic splines. Figure \ref{simuSIR} shows the convergence of the observer for the eigenvalues $\{-2,-2.2,-2.4\}$.
We also simulated the observer when the measurements are corrupted by random counting errors up to $\pm 5$ individuals per day (see Figure \ref{simuSIRnoise}).
As for the adjustable observer in Section \ref{secobslinearpractic}, these simulations show the dilemma {\em accuracy versus speed} of the estimation in presence of measurement noise.

\clearpage

\bleu{
We end this section by showing an example for which an observer is used to reconstruct state and parameter simultaneously.
\section{\bleu{A case study : An observer to estimate state and parameter}}
Consider the following simple  bilharzia transmission model \cite{Macdonald78}:
\begin{equation}\label{schisto}
\left\{ 
\begin{array}{l}
\dot w=-\gamma w+a s\\
\dot s=b(1-s)w-\mu s
\end{array}\right.
\end{equation}
where $w$ is the number of female schistosomes (worms) \textit{per} single host and $s$ the proportion of infected snails. The female schistosomes \textit{per} single host decays at a per capita rate $\gamma$ and are replenished at a rate $a$. The latter process is proportional to the proportion of infected snails $s$. The proportion of susceptible snails $1-s$ are infected through an indirect contact with schistosomes that are excreted from hosts $w$ at rate $b$, and naturally die at a rate $\mu$. 

The number of female schistosomes (worms) \textit{per} single host can be measured using urine or faeces  samples. Therefore, we can assume that the measurable output is $y(t)=w(t)$. 

It is easy to show that the compact set 
\begin{equation}\label{schistocompact}
\Omega=\{(w,s)\in\R^2 \,\,\vert\;\; 0\leq w\leq \dfrac{a}{\gamma}; \,\, 0\leq s\leq 1\}
\end{equation}
is a positively invariant set for System~\eqref{schisto}.
The measurable output $y(t)$  satisfies the equation $\dot y= -\gamma y + a s$. Hence, $y(t) >0$ for all $t\geq 0$ if $y(0)>0$.

We assume that the ecosystem is in an endemic situation which means that the basic reproduction number  $\mathcal{R}_0=
\dfrac{a b}{\gamma \mu}>1$ and implies \cite{Macdonald78} that 
$y(t)$ converges to $ \bar y=\dfrac{a}{\gamma } \left(1-\dfrac{\gamma  \mu }{a b}\right)>0$ as time goes to infinity.

A crucial problem in epidemic models is the estimation of the transmission parameters.  For the schistosomiasis model, it is the parameter $a$ which represents the snail-host infection rate that is difficult to estimate 
\cite{Macdonald78}.

In order to study the observability and identifiability of model~\eqref{schisto},
we consider the augmented system by adding the unkown parameter $a$ to the augmented state $x=(w,s,a)^\top$:
\begin{equation}\label{schis-augm}
\left\{
\begin{array}{l}
\dot x = \left[\begin{array}{c}
\dot w \\
\dot s \\
\dot a 
\end{array}\right]
=
\left[\begin{array}{c}
-\gamma w+a s\\ 
b(1-s)w-\mu s \\
0
\end{array}\right]
=f(x) \\
\\
y= w =h(x)
\end{array}\right.
\end{equation}

One determines
\begin{align*}
    & \dot y= \mathcal L_f h(x) =   \langle \nabla h (x)  | f (x) \rangle
 = -\gamma w+a s
 \end{align*}
 \begin{align*}
 \ddot y=\mathcal L_{f}^{2} h(x) & =   
\langle \nabla \mathcal L_f h(x)  | f (x) \rangle
=\langle\;
\left[\begin{array}{c}
			-\gamma \\[3mm]
			a \\[3mm]			 			
			s
		\end{array}\right]
| \left[\begin{array}{c}
			-\gamma w+a s \\[3mm]
			b(1-s)w-\mu s \\[3mm]
			0
		\end{array}\right]
 \;\rangle\\
& =-\gamma (-\gamma w+a s) + a (b(1-s)w-\mu s).
\end{align*}
The map 
\[\begin{array}{ccl}
H_3 :\{0< w\leq \dfrac{a}{\gamma}; \,\, 0< s< 1; \,\, a>0\}&\rightarrow& \R^3 \\
x&\mapsto & H_3(x)=\Big(h(x), \mathcal L_{f} h(x), \mathcal L_{f}^{2} h(x)\Big)
\end{array}
\]
is injective.
Indeed 
\[H_3(\bar x)= H_3(x) \Rightarrow 
\left\{
\begin{array}{l}
\bar w = w,\\
-\gamma \bar w +\bar a \bar s=-\gamma w + a  s \\
 \bar a (b(1-\bar s)\bar w-\mu \bar s) =
 a (b(1-s)w-\mu s)
\end{array}
\right.
\]
\[
\Rightarrow 
\left\{
\begin{array}{l}
\bar w = w,\\
\bar a \bar s= a  s, \\
 \bar a b(1-\bar s)\bar w =
 a b(1-s)w
\end{array}
\right.
\Rightarrow 
\left\{
\begin{array}{l}
\bar w = w,\\
\bar a \bar s= a  s, \\
 \bar a  =
 a \;(\text{ since } b>0, w>0, \bar w >0, s<1, \bar s <1).
\end{array}
\right.
\]
\[
\Rightarrow 
\left\{
\begin{array}{l}
\bar w = w,\\
\bar s=  s \;(\text{ since } a>0),\\
\bar a  = a.
\end{array}
\right.
\Rightarrow \bar x =x.
\]
%
%
%
%
%
Therefore, thanks to Proposition~\ref{prop:H_k}, the augmented system~\eqref{schis-augm} is  observable. Using chap 1 Proposition~\ref{prop:augm}, we deduce that  Model~\eqref{schisto} is observable and identifiable.
Moreover, the state variables $w$ and $s$ as well as the parameter $a$ can be expressed as rational functions of $y$, $\dot y$ and $\ddot y$ as follows:
\[
\begin{cases}
w = y, \\
s = 
\dfrac{\left(\gamma  y +\dot y \right) b}
{(b  +  \mu)  \gamma  y +(b  +\gamma  +\mu ) \dot y +\ddot y}, 
\\[3mm]
a = 
\dfrac{(b  +  \mu)  \gamma  y +(b  +\gamma  +\mu ) \dot y +\ddot y}{b}
=\dfrac{\gamma  y +\dot y}{s}.
\end{cases}
\]
%
 %

Now, we shall built an observer that will allow to estimate the unmeasured state variable (here it is $s(t)$) as well as the unknown parameter $a$. To this end, we perform  the following change of coordinates:
\[
z_1=w, \quad z_2=-\gamma w+ x_2 a, \quad z_3=a.
\]
One has then
\[
\left\{\begin{array}{lll}
\dot z_1 & = & -\gamma z_1+z_3 s = -\gamma z_1+z_3 \dfrac{z_2 + \gamma z_1}{z_3}
=z_2 \\[2mm]
\dot z_2 & = & -\gamma (-\gamma z_1+z_3 s)+ (b(1-s) w-\mu s) z_3
= -\gamma z_2 + (b(1-s) w-\mu s) z_3 \\[2mm]
& = & -\gamma z_2 + (b(1-\dfrac{z_2 + \gamma z_1}{z_3}) z_1-\mu \dfrac{z_2 + \gamma z_1}{z_3}) z_3\\
& = & -\gamma z_2 + b(z_3-z_2 - \gamma z_1) z_1
-\mu z_2 -\mu \gamma z_1 \\[2mm]
\dot z_3 & = & 0
\end{array}\right.
\]
Note that, since $y=z_1$, the dynamics of $z_2$ can be written as
\[
\dot z_2 =-\gamma (\mu +b y) z_1 -(\mu+\gamma+b y)z_2 
+  b y z_3
\]
%
Therefore, the dynamics takes the form
\begin{equation}\label{adap0}
\left\{
\begin{array}{l}
\dot z(t) =A(y)\,z(t),\\[3mm]
y=C_0\, z, 
\end{array}
\right.
\end{equation}
where
\[
A(y)=\left[\begin{array}{ccccc}
0 && 1 && 0\\[3mm]
-\gamma (\mu +b y) && -(\mu+\gamma+b y) && b y \\[3mm]
0 && 0 && 0
\end{array}\right] , \quad
C_0=[1 \; 0 \; 0]
\]
For any fixed $y>0$, the corresponding observability matrix  is
\[O_{(C_0,A)}=\left[\begin{array}{c}
			C_0\\
			C_0 A\\\
			C_0 A^{2}
		\end{array}\right]
=\left[\begin{array}{ccc}
1 & 0 & 0 
\\
 0 & 1 & 0 
\\
 -\gamma  \left(b y +\mu \right) & -b y -\gamma -\mu  & b y  
\end{array}\right]
\]
that is of full rank if $y\neq 0$. Therefore, by the pole-shifting theorem (see \cite[page 61]{0424.93001}), it is possible to find a $y-$dependent gain $ K(y)$ such that $sp(A(y)-K(y) C_0)=\{-\lambda_1, -\lambda_2, -\lambda_3\}$, where $sp(M(y))$ denotes the spectrum of $M(y)$ and $\lambda_i$ are any positive real numbers.
This gain $ K(y)$  can be computed using for instance Ackermann's formula (see \cite{Antsaklis2007} page 382):
\[K(y)=\prod_{i=1}^{3}
(A(y)+\lambda_i I_3). 
O_{(C_0,A)}^{-1}. 
\left[\begin{array}{c}
0  \\
 0 \\
1\\
\end{array}\right],
\]
where $I_3$ is the $3\times 3$ identity matrix.
\[K(y)=\left[\begin{array}{c}
K_1(y)\\
K_2(y)\\
K_3(y)
\end{array}\right]=
\left[\begin{array}{l}
 -(\gamma +\mu) +\lambda_1 +\lambda_2 +\lambda_3  -b y
\\[3mm]
 \gamma^{2}+\mu^{2}+\mu \gamma
-(\mu +\gamma) \left( \lambda_1 +\lambda_2 +\lambda_3 \right)   +
\left(\lambda_2 +\lambda_3 \right) \lambda_1 +\lambda_2 \lambda_3
\\
+b^{2} y^{2}- \left(\lambda_1 +\lambda_2 +\lambda_3 -\gamma -2 \mu\right) b y  
\\[3mm]
 \dfrac{\lambda_1 \lambda_2 \lambda_3}{b y} 
\end{array}\right]
\]
The gain $K(y)$ is well defined since $y(t) >0$ for all $t\geq 0$.
An observer is then given by:
\begin{equation}\label{obs-adap0}
\dot {\hat z}(t) =A(y){\hat z}(t)-K(y(t)) (C_0\hat z(t) - y(t)).
\end{equation}
In coordinates:
\begin{equation}\label{obs-adap0c}
\left\{\begin{array}{l}
\dot {\hat z}_1 ={\hat z}_2 -({\hat z}_1-y) K_1(y)\\[2mm]
\dot {\hat z}_2 =-\mu \gamma {\hat z}_1 - (\mu + \gamma){\hat z}_2 +  b({\hat z}_3-{\hat z}_2 - \gamma {\hat z}_1) y
-({\hat z}_1-y) K_2(y)\\[2mm]
\dot {\hat z}_3 = -({\hat z}_1-y) K_3(y) 
\end{array}\right.
\end{equation}
The error equation
\begin{equation}\label{obs-adap0-error}
\dot e(t)= \Big(A(y) -K(y)\, C_0\Big)\, e(t)=M(y)\, e(t).
\end{equation}
The eigenvalues of the matrix $M(y)$ are $-\lambda_1$, $-\lambda_2$ and $-\lambda_3$.
It has been proved in \cite{bichara:hal-02189643} that the error $e$ converges exponentially fast to zero (the proof being quite long is omitted here), which gives the exponential convergence of the observer~\eqref{obs-adap0}.
%


%
%
%
%
}

\begin{figure}[!ht]
\begin{center}
\includegraphics[scale=0.3]{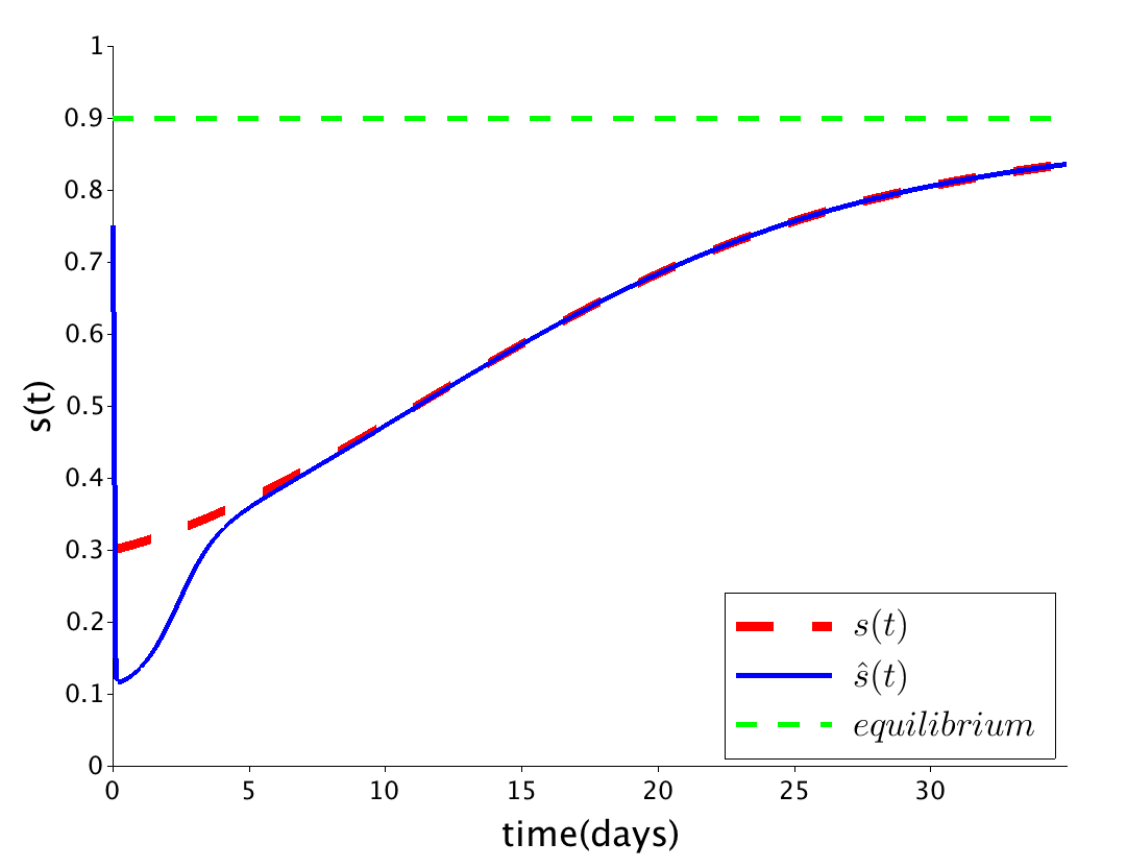}\includegraphics[scale=0.3]{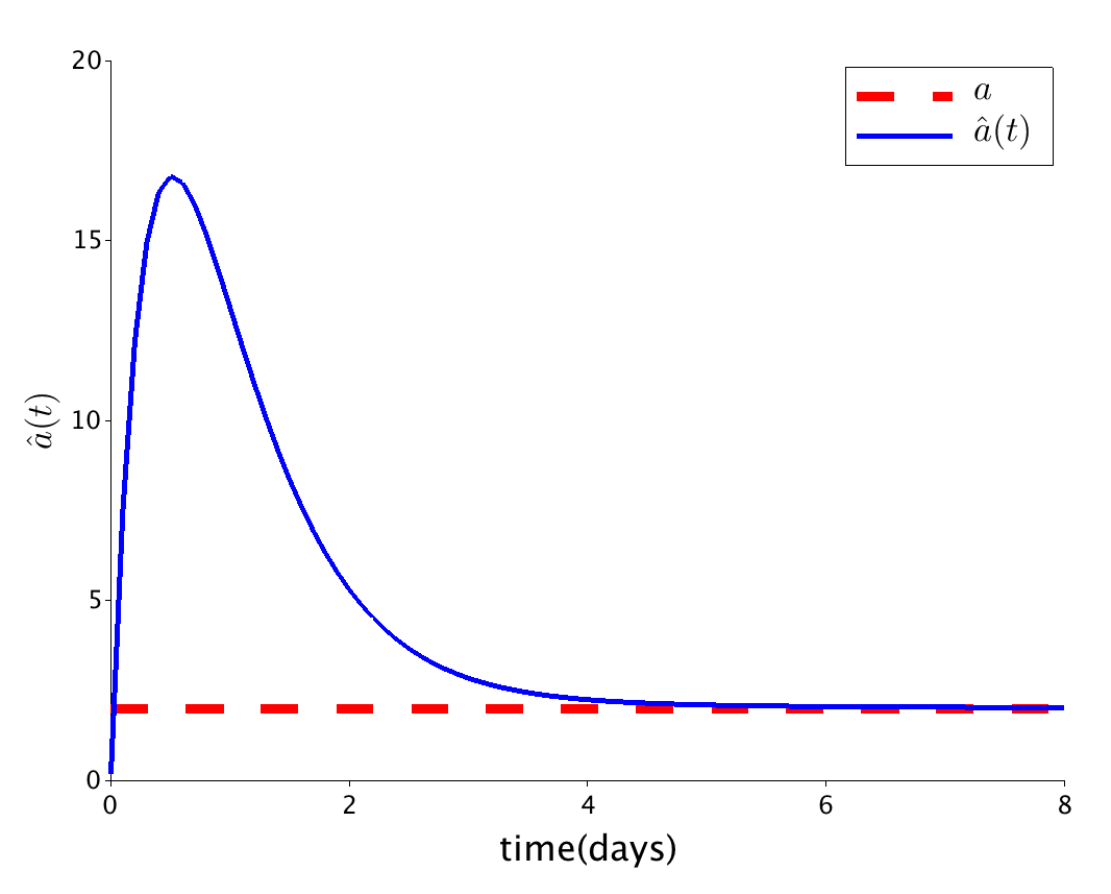}
\caption{Dynamics of proportion of infected snails $s(t)$ with its estimate $\hat s(t)$ (left) and of the estimation $\hat a(t)$ of parameter $a$ (right).  
 The green dashed line corresponds to the steady state value of $s$.}
\label{adapobs}
\end{center}
\end{figure}

\bigskip
\color{blue}
Finally, let us illustrate this observer on numerical simulations. We have taken the following values of the parameters
\begin{center}
	\begin{tabular}{|c|c|c|c|}
		\hline $\vphantom{\displaystyle \int}$
		$\beta$ & $\mu$ & $a$ & $b$\\
		\hline\hline
		 $\vphantom{\displaystyle \int}$ \; \verb"0.05" \;  & \; \verb"0.04" \; & \; \verb"2" \; & \; \verb"0.01" \;\\
		\hline
	\end{tabular}
\end{center}
The initial conditions are $(w(0),s(0))=(3,0.3)$ and $(\hat z_1(0),\hat z_2(0),\hat z_3(0))=(1,\;0.1,\; 0.2)$.
The set of eigenvalues of the matrix $M(y)$ have been chosen to be  $\{-0.4,-1.4,-2.4\}$. 
Figure~\ref{adapobs} shows the convergence of the estimation of the unmeasured proportion of infected snails $s(t)$
\[
\hat s(t)=\dfrac{\hat z_2(t) + \gamma \hat z_1(t)}{\hat z_3(t)}
\]
and of unknown parameter $a$
\[
\hat a(t)=\hat z_3(t)
\]
delivered by the observer~\eqref{obs-adap0c}.
\n

The same simulations have been conducted in {\tt Scilab} with measurement noise
\begin{verbatim}
y=x(5)+0.1*grand(1,1,"nor",1,1)
\end{verbatim}
\begin{figure}[!ht]
\begin{center}
\includegraphics[height=50mm]{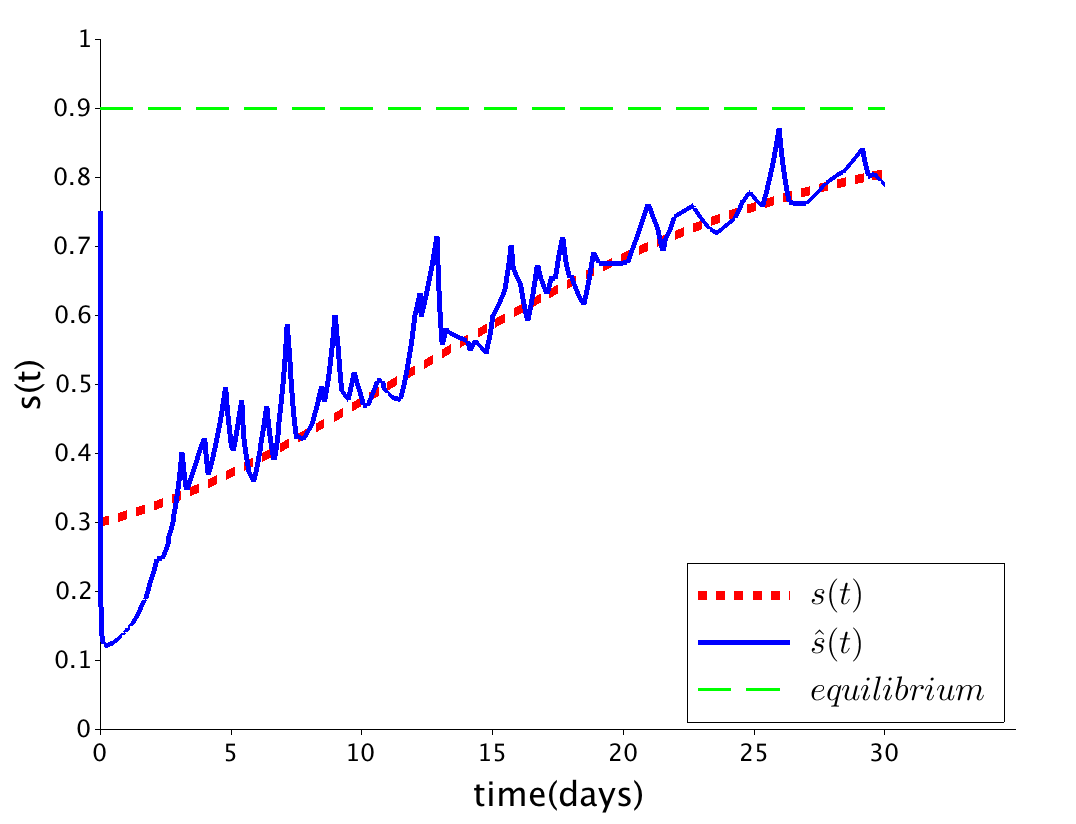}\includegraphics[height=50mm]{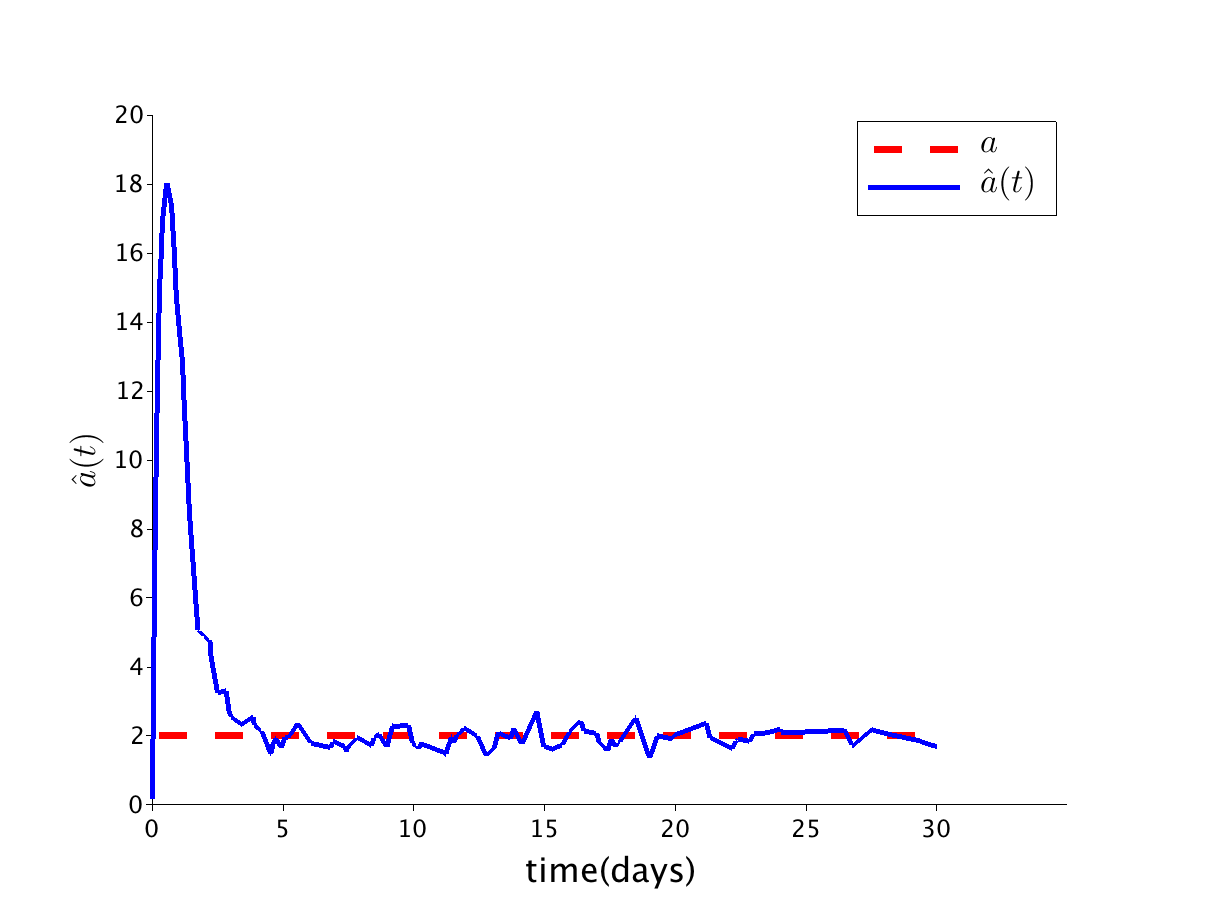}
\caption{Simulations with measurement noise.}
\label{adapobs-noise}
\end{center}
\end{figure}
Figure \ref{adapobs-noise} shows that for this choice of gains, the estimations are heckled but follow quite well the unknown state and parameter.
\color{black}

%
%

%
%
%


\appendix

\chapter{Proofs of some useful lemmas}

\section{Proof of Lemma \ref{lem_poles}}
\label{appendixLemmaPoles}\m

The proof is adapted from \cite{AndreaLara}.\m

Let us first consider pairs $(\bar A, \bar C)$ of the canonical form known as Brunovsky's form
\[
\bar A=\left[\begin{array}{cccccc}
	0 & \cdots & \cdots & \cdots & 0 & -a_{n}\\
	1 & 0 & \cdots & \cdots & 0 & -a_{n-1}\\
	\cdot & \ddots & \ddots & \cdot & \cdot & \cdot\\
	\cdot & \cdot & \ddots & \ddots & \cdot & \cdot\\
	\cdot & \cdot & \cdot & 1 & 0 & -a_{2}\\
	0 & \cdots & \cdots & 0  & 1 & -a_{1}
\end{array}\right], \quad \bar C=\left[\begin{array}{ccc}
	0 & \cdots 0 & 1\end{array}\right]
\]
where the $a_i$ are any numbers. Their observability matrices are lower triangular:
\[
\bar O = \left [\begin{array}{ccccc}
	& & & & 1\\
	\multicolumn{2}{c}{\text{\huge{$0$}}} & \reflectbox{$\ddots$} & & \\
	& \reflectbox{$\ddots$} & \multicolumn{3}{c}{\text{\Huge{$\star$}}}\\
	1 & & & &
\end{array}\right]
\]
therefore invertible.
It is easy to see that the characteristic polynomial of the matrix $\bar A$ is given by 
\[
\pi_{\bar A}(\xi)=\xi^n+a_{1}\xi^{n-1}+\cdots +a_{n-1}\xi+a_{n}
\]
Indeed, if $X$ is a left eigenvector of $\bar A$ for an eigenvalue $\lambda$ (possibly complex), 
$X\bar A =\lambda X$ gives
\[
\begin{array}{l}
	X_2=\lambda X_1 ,\\
	X_3=\lambda X_2=\lambda^2 X_1 ,\\
	\;\;\vdots\\
	X_n=\lambda X_{n-1}=\lambda^{n-1}X_1 ,\\
	-a_nX_1-a_{n-1}X_2- \cdots - a_1 X_n=\lambda X_n .
\end{array}
\]
Thus the line vector $X$ is of the form
\[
X=\left[\begin{array}{ccccc}
	1 & 
	\lambda & 
	\lambda^2 & 
	\cdots & 
	\lambda^{n-1}
\end{array}\right] X_1
\quad \mbox{with } X_1 \neq 0
\]
and $\lambda$ verifies
\[
\big(\lambda^n + a_1 \lambda^{n-1}+a_2 \lambda ^{n-2}+ \cdots + a_{n-1}\lambda+a_n\big)X_1=0 \ .
\]
Since $X_1$ is non-zero, we deduce that the eigenvalues are roots of the polynomial
\[
\lambda^n + a_1 \lambda^{n-1}+a_2 \lambda ^{n-2}+ \cdots + a_{n-1}\lambda+a_n=0
\]
which is of degree $n$ and whose coefficient of $\lambda^n$ is equal to $1$.

\medskip

The characteristic polynomial of the matrix $\bar A+\bar G\bar C$, where $\bar G$ is a vector of $\R^n$ with elements denoted $\bar g_i$, is written as follows
\[
\pi_{\bar A+\bar G\bar C}(\xi)=\xi^n+(a_{1}-\bar g_n)\xi^{n-1}+\cdots +(a_{n-1}-\bar g_2)\xi+(a_{n}-\bar g_1)
\]
Thus, one can arbitrarily choose the $n$ coefficients of this polynomial by choosing the $n$ elements of $\bar G$, and thus freely assign the spectrum of the matrix $\bar A+\bar G\bar C$.
For any set $\Lambda=\{\lambda_1,\cdots,\lambda_n\}$ of $n$ real or complex numbers two by two conjugates, one has just to identify the coefficients of the polynomial $\pi_{\bar A+\bar G\bar C}$ with those of
\[
\prod_{i=1}^n (\xi-\lambda_i)=\xi^n + \sum_{k=1}^{n} (-1)^k\sigma_{k}(\Lambda) \xi^{n-k}
\]
Thus, we obtain
\[
\bar g_i=a_{n+1-i}+(-1)^{n-i}\sigma_{n+1-i}(\Lambda), \quad i=1 \cdots n
\]

\medskip

Let us now show that for any pair $(A,C)$ such that $O$ is full rank, there is an invertible $P$ matrix such that 
$P^{-1}AP=\bar A$ and $CP=\bar C$, where the pair $(\bar A,\bar C)$ is in the Brunovsky's form. Consider the vector
\[
L=O^{-1}
\left[\begin{array}{c}
	0\\ \vdots\\ \vdots\\ 1
\end{array}\right] \; \Rightarrow \left\{\begin{array}{ll}
	CA^kL &=0, \quad k=0 \cdots n-2\\
	CA^{n-1}L &=1
\end{array}
\right.
\]
and the matrix consisting of the concatenation of the columns
\[
P=[L \; AL \; \cdots \; A^{n-1}L]
\]
We have
\[
OL=\left[\begin{array}{c}
	0\\ \vdots\\ \vdots\\ \vdots \\ 1
\end{array}\right] , \; 
OAL=\left[\begin{array}{c}
	CAL\\ CA^2L\\ \vdots\\ CA^{n-1}L\\ CA^nL
\end{array}\right]=\left[\begin{array}{c}
	0\\ \vdots\\ \vdots\\ 1\\ \star
\end{array}\right] , \;  
OA^2L=\left[\begin{array}{c}
	0\\ \vdots\\ 1\\ \star \\ \star\end{array}\right] , \cdots
\]
up to
\[
OA^{n-1}L=\left[\begin{array}{c}
	1\\ \star \\ \vdots \\ \vdots\\ \vdots
\end{array}\right] 
\]
Thus the $OP$ matrix is of the form
\[
OP=\left[\begin{array}{ccccc}
	& & & & 1\\
	\multicolumn{2}{c}{\text{\huge{$0$}}} & \reflectbox{$\ddots$} & & \\
	& \reflectbox{$\ddots$} & \multicolumn{3}{c}{\text{\Huge{$\star$}}}\\
	1 & & & &
\end{array}\right]
\]
which shows that $P$ is indeed an invertible matrix. Finally, the columns of the $AP$ matrix are
\[
AP=[AL \; A^2L \; \cdots \; A^{n-1}L \; A^n L]
\]
Its first $n-1$ columns are written as follows
\[
[AL \; A^2L \; \cdots \; A^{n-1}L]=P\left[\begin{array}{cccc}
	0 & \cdots &  &    \\
	1 & 0 & \cdots &    \\
	0 & \ddots & \ddots &   \\
	\vdots &    & 1 & 0 \\
	0  & \cdots  & 0  & 1 
\end{array}\right]
\]
By Cayley-Hamilton's Theorem, we have $\pi_A(A)=0$, which allows us to write the last column of $AP$ as 
\[
A^nL=-a_nL-a_{n-1}AL-\cdots -a_1A^{n-1}L=P\left[\begin{array}{c}
	-a_n\\
	-a_{n-1}\\
	\vdots\\
	-a_1
\end{array}\right]
\]
which shows that we have $PA=\bar AP$. Finally, one gets
\[
CP=[CL \; CAL \; \cdots \; CA^{n-}]=[0 \; \cdots \; 0 \; 1]=\bar C
\]
Then, for a vector $\bar G$ such that $Sp(\bar A+\bar G\bar C)=\Lambda$, we have $Sp(P^{-1}(\bar A+\bar G\bar C)P)=\Lambda$, where $P^{-1}(\bar A+\bar G\bar C)P=A+P\bar G C$. We conclude that for the gain vector $G=P\bar G$, we have $Sp(A+GC)=\Lambda$.

\bigskip

\section{Proof of Lemma \ref{lemVandermonde}}
\label{appendixLemmaVandermonde}\m

The proof is adapted from \cite{Ciccarella93}.\m

Let $X$ be a left eigenvector of $A+GC$ for the eigenvalue $\lambda_i$. By writing $X(A+GC)=\lambda_i X$, we obtain the $n-1$ inequalities.
\[
\begin{array}{lll}
	X_1 & = & \lambda_i X_2\\
	X_2 & = & \lambda_i X_3\\
	& \vdots & \\
	X_{n-1} & = & \lambda_i X_n
\end{array}
\]
Thus $X_n$ is necessarily non-zero and can be taken equal to $1$, which gives
\[
X = \left [\begin{array}{ccccc} \lambda_i^{n-1} & \lambda_i^{n-2} & \cdots & \lambda_i & 1 \end{array}\right]
\]
We then obtain the $n$ rows of the matrix $V_{\lambda_1,\cdots,\lambda_n}$, which defines a matrix of change of basis that diagonalizes the matrix $A+GC$.

\medskip

Now let's show how to determine the inverse of $V_{ \lambda_1,\cdots,\lambda_n }$. Let $w_{ij}$ be the coefficients of $V_{ \lambda_1,\cdots,\lambda_n }^{-1}$. The equality $(V_{ \lambda_1,\cdots,\lambda_n }^{-1})(V_{ \lambda_1,\cdots,\lambda_n })=Id$ gives
\begin{equation}
	\label{wkj}
	\sum_{k=1}^n w_{kj}\lambda_i ^{n_k}=\delta_{ij}:=\left|\begin{array}{ll}
		0 & \mbox{si } i \neq j\\
		1 & \mbox{si } i=j
	\end{array}\right.
\end{equation}
For each $j$ in $\{1,\cdots,n\}$, let's consider the polynomial
\begin{equation}
	\label{Pj}
	P_j(X)=\sum_{k=1}^n w_{kj}X ^{n-k}
\end{equation}
The conditions \eqref{wkj} amount to write $P_j(\lambda_i)=\delta_{ij}$, i.e. the polynomial $P_j$ has $n-1$ roots $\lambda_i$ for $i\neq j$ and $P_j(\lambda_j)$ is equal to $1$. So it has the following expression
\[
P_j(X) = \prod_{k \neq j} \frac{X-\lambda_k}{\lambda_j-\lambda_k}
\]
By identifying its coefficients with those of the expression \eqref{Pj}, we obtain
\begin{equation}
	\label{wij}
	w_{ij}=(-1)^{i-1}\frac{\sigma_{i-1}(\Lambda\setminus\{\lambda_j\})}{ \prod_{k\neq j}\lambda_j-\lambda_k}
\end{equation}
where the $\sigma_k$ are the symmetric functions defined in \eqref{sigma}.

\medskip

Let
\[
\varphi(\lambda_1,\cdots,\lambda_n)=\lambda_1+c||V_{\lambda_1,\cdots,\lambda_n}^{-1}||_\infty+\theta
\]
The expression \eqref{wij} shows that the norm $||V_{\lambda_1,\cdots,\lambda_n}^{-1}||_\infty$ becomes arbitrarily large when $\lambda_i-\lambda_j$ approaches $0$ (for $i\neq j$), which ensures the existence of numbers
$\lambda_n <\lambda_{n-1}< \cdots < \lambda_1<0$ such as
$\varphi(\lambda_1,\cdots,\lambda_n)>0$.
For $\lambda_i=-\alpha^i$ ($i=1,\cdots,n$), we obtain, for any $j$
\[
\lim_{\alpha \to +\infty} w_{ij} = \left|\begin{array}{ll}
	0 & i < n\\
	1 & i=n
\end{array}\right.
\]
and $||V_{-\alpha,-\alpha^2,\cdots,-\alpha^n }^{-1}|||_\infty$ thus tends towards $1$ when $\alpha$
tends towards $+\infty$, which shows the existence of numbers $\lambda_n <\lambda_{n-1}< \cdots < \lambda_1<0$ such that $\varphi(\lambda_1,\cdots,\lambda_n)<0$.
Finally, by continuity of $\varphi$, we deduce the existence of $\lambda_n <\lambda_{n-1}< \cdots < \lambda_1<0$ such that $\varphi(\lambda_1,\cdots,\lambda_n)=0$.

\bigskip

\color{blue}
\section{Proof of Theorem \ref{LyapunovTh}}
\label{appendixPropLyapunov}\m

The proof is adapted from \cite{Khalil}.\m

Let $P$ be a symmetric positive definite matrix satisfying $M^\top P+PM+Q=0$ where $Q$ is a symmetric positive definite matrix. We consider the Lyapunov function
\[
V(x)=||x||_P=x^\top Px
\]
whose time derivative along solutions of $\dot x=Mx$ is
\[
\frac{d}{dt}V(x(t))= -x(t)^\top Qx(t) \leq - \underbrace{\dfrac{\lambda_{min}(Q)}{\lambda_{max}(P)}}_\beta x(t)^\top P x(t)
\]
(where $\lambda_{min}$, $\lambda_{max}$ denote the smallest and largest real eigenvalues of a symmetric matrix).
Then, one has
\[
\frac{d}{dt}V(x(t)) \leq -\beta V(x(t))
\]
where $\beta>0$.
We deduce that $t \mapsto ||x(t)||_P$ converges exponentially to $0$.

\bigskip

Conversely, when $M$ is Hurwitz, one can consider the symmetric matrix
\[
P=\lim_{t \to +\infty}\int_0^{t} e^{\tau M^\top}Qe^{\tau M} d\tau
\]
(which exists as a sum of terms $t^{k}e^{\lambda_i t}$ where $Re \lambda _i<0$).
Then, one has
\begin{align*}
    M^\top P+PM & = \lim_{t \to +\infty}\int_0^{t} M^\top e^{\tau M^\top}Qe^{\tau M} 
    + e^{\tau M^\top}Qe^{\tau M}Md\tau\\
    & =  \lim_{t \to +\infty}\int_0^{t} \frac{d}{d\tau}\left(e^{\tau M^\top}Qe^{\tau M}\right)d\tau\\
    & = \lim_{t \to+\infty} e^{t M^\top}Qe^{t M} - Q \\
    & = -Q
\end{align*}
Let us show that $P$ is necessarily definite positive. If there exists $u \neq 0$ such that $P u=0$, then one has
\[
u^T P u =0 = u^T \left(\int_0^{\infty} e^{\tau M^\top}Qe^{\tau M} d\tau\right) u
=\int_0^{\infty} ||Q^{1/2}e^{\tau M} u ||^2 d\tau
\]
where $Q^{1/2}$ is the positive definite matrix such that $Q^{1/2}Q^{1/2}=Q$. Therefore, $u$ has to be null which shows that $P$ is non singular.

\bigskip

Finally, if there exists another symmetric definite positive matrix $\tilde P$ which satisfies $M^\top \tilde P +\tilde P M=-Q$, then one has
\[
M^\top(\tilde P-P)+(\tilde P -P)M =0
\]
But then
\[
e^{t M^\top}\left(M^\top(\tilde P-P)+(\tilde P -P)M\right)e^{t M}=\frac{d}{dt}\left(e^{t M^\top}(\tilde P-P)e^{t M}\right)=0
\]
i.e.~$t \mapsto e^{t M^\top}(\tilde P-P)e^{t M}$ is a constant function. In particular,
\[
e^{0 M^\top}(\tilde P-P)e^{0 M}=\tilde P-P=\lim_{t\to+\infty}e^{t M^\top}(\tilde P-P)e^{t M}=0
\]
which gives $\tilde P=P$.
\color{black}


\chapter{Implementation of the ``Boarding School'' example}
\label{appendixCode}

\section{Derivation of the Fisher Information Matrix}

In this example, we consider $x\in\mathbb{R}^n$ with $n=2$, $y\in\mathbb{R}^m$ with $m=1$, i.e.,
\[
x(t)=\left[\begin{array}{c}S(t)\\I(t)\end{array}\right]\,,\quad y(t)= I(t)\,,
\]
and $\theta\in\mathbb{R}^p$ with $p=2$, i.e.,
\[
\theta=\left[\begin{array}{c}\beta\\\gamma\end{array}\right]\,.
\]
We consider the following model, equivalent to model \eqref{KmcK} with $k=1$:
\[
\left\{\begin{array}{lll}
	\dot x&=&f(x,\theta)=\left[\begin{array}{c}\dot S\\\dot I\end{array}\right]=\left[\begin{array}{c}f_S(S,I,\theta)\\f_I(S,I,\theta)\end{array}\right]=\left[\begin{array}{c}-\beta SI/N\\\beta SI/N-\gamma I\end{array}\right]\,,\quad x(0)=x_0=\left[\begin{array}{c}S_0\\I_0\end{array}\right]\,,\\
 \\
	y &=& h(x,\theta) = I\,.
\end{array}\right.
\]
We disregard $\Theta=(\theta,x_0)$ since the initial conditions are assumed to be known in this example. The Jacobian of the observation with respect to the parameter $\theta$ is:
\[
\chi(t,\theta)=\frac{\partial y}{\partial\theta}(t)=\frac{\partial h}{\partial x}\frac{\partial x}{\partial \theta}(t)\,,
\]
since in this example,
\[
\frac{\partial h}{\partial \theta}(t)=0\,.
\]
This Jacobian has dimension $m\times p = 1 \times 2$. We have
\[
\frac{\partial h}{\partial x}=\left[\begin{array}{cc} \dfrac{\partial h}{\partial S} & \dfrac{\partial h}{\partial I}\end{array}\right]=\left[\begin{array}{cc}0 & 1\end{array}\right]
\]
and 
\[
z = \dfrac{\partial x}{\partial\theta}=\left[\begin{array}{cc}
	\dfrac{\partial S}{\partial \beta} & \dfrac{\partial S}{\partial \gamma}\vspace{.1cm}\\
	\dfrac{\partial I}{\partial \beta} & \dfrac{\partial I}{\partial \gamma} \end{array}\right]\,.
\]
This yields
\[
\chi=\dfrac{\partial h}{\partial x}\dfrac{\partial x}{\partial \theta}=   \left[\begin{array}{cc} \dfrac{\partial I}{\partial \beta} & \dfrac{\partial I}{\partial \gamma}\end{array}\right]\,.
\]
Let $\{t_i\}$, $i=0,1,2,\ldots,M$, be the sampling times. Fisher's Information Matrix is defined as:
\[
\text{FIM}(\theta,\sigma)=\frac{1}{\sigma^2}\sum_{i=1}^M \chi(t_i,\theta)^\top\chi(t_i,\theta)\,,
\]
where $\sigma^2$ is defined as the sum of the squared error (SSE) divided with $M-p$ instead of $M-(n+p)$ as in equation \eqref{estimbias}, since the initial conditions are assumed to be known in this example.

\paragraph{Computing Fisher's Information Matrix} Let 
\[
A(t)=\dfrac{\partial f}{\partial x}=\left[\begin{array}{cc}
	\dfrac{\partial f_S}{\partial S} & \dfrac{\partial f_S}{\partial I} \vspace{.1cm}\\
	\dfrac{\partial f_I}{\partial S} & \dfrac{\partial f_I}{\partial I} \end{array}\right]=\left[\begin{array}{ccc}
	-\dfrac{\beta I}{N} && -\dfrac{\beta S}{N} \\[2mm]
	\dfrac{\beta I}{N} && \dfrac{\beta S}{N} - \gamma \end{array}\right]
\]
and
\[
B(t)=\frac{\partial f}{\partial \theta}=\left[\begin{array}{cc}
	\dfrac{\partial f_S}{\partial \beta} & \dfrac{\partial f_S}{\partial \gamma} \vspace{.1cm}\\
	\dfrac{\partial f_I}{\partial \beta} & \dfrac{\partial f_I}{\partial \gamma} \end{array}\right]=\left[\begin{array}{ccc}
	-\dfrac{SI}{N} && 0 \\[2mm]
	\dfrac{SI}{N} && -I \end{array}\right]\,.
\]
The matrix $z$ can be computed by numerically solving the following system of ODE's:
\[
\left\{\begin{array}{lll}
	\dot x &=& f(x,\theta)\,,\quad x(0)=x_0\,,\\
	\dot z &=& Az+B\,,\quad z(0)=0_{n\times p}
\end{array}\right.
\]
which is a subsystem of (\ref{eq:fullSystem}) since the initial conditions are assumed to be known in this example (i.e., we disregard $w$). In the following code, the entries of $x$ and $z$ are indexed in this way:
\[
x=\left[\begin{array}{c} x_1 \\ x_2 \end{array}\right]\,,\quad z=\left[\begin{array}{cc} z_3 & z_5 \\ z_4 & z_6 \end{array}\right]\,,
\]
leading to 
\[
\chi=\left[\begin{array}{cc}z_4 & z_6\end{array}\right]\,.
\]
\color{black}
\section{Numerical implementation}

The code has been written with the {\tt Scilab} language and executed under {\tt SCILAB 6.0.0}\footnote{\tt https://www.scilab.org/}.  It consists in a function for identifying $\beta$ and $\gamma$ and using the {\tt lsqrsolve} function which implements the Levenberg-Marquard algorithm to perform ordinary least squares. We could have chosen the {\tt fminsearch} function which is an implementation of the Nelder-Mead algorithm, but this gives exactly the same results.  For solving ODE's, {\tt Scilab} uses the  {\tt lsoda} solver of {\tt ODEPACK}. It automatically selects between non-stiff predictor-corrector Adams method and stiff Backward Differentiation Formula (BDF) method. It uses non-stiff method initially and dynamically monitors data in order to decide which method to use.\m

\n
We define the following functions in the {\tt Scilab} environment:\m\n

\begin{verbatim}
	function [kguess_n,SSE]=identifKMK(OBS,T,kguess,N)
	// kguess_n = [BETA;GAMMA]
	t0=T(1);m=length(T);
	[x,SSE]=lsqrsolve(kguess,errorKmcK,m,[1.d-8,1.d-8,1.d-5,1d9,0,100]);
	nbparam=length(kguess)
	kguess_n=x;
	x0=[N-OBS(1);OBS(1)];
	sol=ode(x0,t0,T,list(KmcK,x(1),x(2)));
	kguess_n=kguess_n(:);
	SSE=sum(SSE.^2)
	xset("window",1)
	sol1=ode(x0,t0,T(1):0.01:T($),list(KmcK,x(1),x(2)));
	clf
	plot((T(1):0.01:T($))',sol1(2,:)')
	plot(T',OBS,'ro')
	endfunction
	
	function y=errorKmcK(k,m)
	x0=[N-OBS(1);OBS(1)];x0=x0(:);
	BETA=k(1);GAMMA=k(2);
	sol=ode(x0,t0,T,list(KmcK,BETA,GAMMA));
	predic=sol(2,:);
	predic=predic(:);
	OBS=OBS(:);
	y=OBS-predic;   
	endfunction
	
	function xdot=KmcK(t,x,BETA,GAMMA,N)   
	xdot=[-BETA/N*x(2),0;BETA/N*x(2),-GAMMA]*x
	endfunction
	
\end{verbatim}
\n\n

Then, the {\tt Scilab} session goes like this\m\n

\begin{verbatim}[commandchars=\\\{\}]
	--> load('databoarding')
	ans  =
	T
	--> OBS=dataBSFlu
	
	OBS  = 
	
	column 1 to 8
	
	1.   6.   26.   73.   222.   293.   258.   237.
	
	column 9 to 14
	
	191.   124.   68.   26.   10.   3.
	
	--> M=length(OBS);
	
	--> T=0:M-1;
	
	--> N=763;
	
	--> beta0=2;gamma0=0.5;param=[beta0,gamma0];
	
	--> p=length(param);
	
	--> [param,SSE]=identifKMK(OBS,T,param,N)
	
	SSE  = 
	
	4892.6472
	
	param  = 
	
	1.9605032
	0.4751562
	
	--> sigma2=SSE/(M-p)
	
	sigma2  = 
	
	407.72060
	
	--> BETA=param(1);GAMMA=param(2);
\end{verbatim} 

\m\n

\medskip

We then compute confidence intervals  using the formulas \eqref{FIM}, \eqref{estimtheta}, \eqref{estimbias}, \eqref{estimcovmat}, \eqref{SE}.\m\n

\begin{verbatim}[commandchars=\\\{\}]
	function FIM=fimKmcK(x0,T,BETA,GAMMA,sigma2)
	// Compute the sensitivity matrix 
	x0=x0(:);t0=T(1);
	X0=[x0;0;0;0;0]
	sol=ode(X0,t0,T,list(JKmcK,BETA,GAMMA));
	M=sol([4,6],:);
	FIM=M*M'./sigma2;
	endfunction
	
	function xdot=JKmcK(t,x,BETA,GAMMA)
	xdot(1)=-BETA*x(1)*x(2)/N
	xdot(2)=BETA*x(1)*x(2)/N-GAMMA*x(2) 
	xdot(3)=-BETA*x(2)*x(3)/N-BETA*x(1)*x(4)/N-x(1)*x(2)/N
	xdot(4)=BETA*x(2)*x(3)/N+(BETA*x(1)/N-GAMMA)*x(4)+x(1)*x(2)/N
	xdot(5)=-BETA*x(2)*x(5)/N-BETA*x(1)*x(6)/N
	xdot(6)=BETA*x(2)*x(5)/N+(BETA*x(1)/N-GAMMA)*x(6)-x(2)
	endfunction
\end{verbatim}

\medskip

Then the {\tt Scilab} session is 

\begin{verbatim}[commandchars=\\\{\}]
	
	--> x0=[N-OBS(1); OBS(1)];
	
	--> FIM=fimKmcK(x0,T,BETA,GAMMA)
	FIM  = 
	
	974.5073   -523.73985
	-523.73985   3132.2047
	
	--> cond(FIM)
	ans  =
	
	3.8082403
	
	--> CovMAT=inv(FIM)
	CovMAT  = 
	
	0.0011275   0.0001885
	0.0001885   0.0003508
	
	--> t=cdft("T",M-p,0.975,0.025)
	t  = 
	
	2.1788128
	
	--> confBETA=t*sqrt(CovMAT(1,1))
	confBETA  = 
	
	0.0731602
	
	--> confGAMMA=t*sqrt(CovMAT(2,2))
	confGAMMA  = 
	
	0.0408077
\end{verbatim}

\bigskip

\chapter{Implementation of the ``Plague in Bombay'' example}
\label{appendixCodeBombay}

\section{Derivation of the Fisher Information Matrix}

In this example, we consider $x\in\mathbb{R}^n$ with $n=2$, $y\in\mathbb{R}^m$ with $m=1$, i.e.,
\[
x(t)=\left[\begin{array}{c}S(t)\\I(t)\end{array}\right]\,,\quad y(t)= I(t)\,,
\]
and $\theta\in\mathbb{R}^p$ with $p=2$, i.e.,
\[
\tilde\beta=\frac{\beta}{N}\,,\quad \theta=\left[\begin{array}{c}\tilde\beta\\\gamma\end{array}\right]\,.
\]
We consider the following model, equivalent to model \eqref{KmcK} with $k=\gamma$:
\[
\left\{\begin{array}{lll}
	\dot x&=&f(x,\theta)=\left[\begin{array}{c}\dot S\\\dot I\end{array}\right]=\left[\begin{array}{c}f_S(S,I,\theta)\\f_I(S,I,\theta)\end{array}\right]=\left[\begin{array}{c}-\tilde\beta SI\\\tilde\beta SI-\gamma I\end{array}\right]\,,\quad x(0)=x_0=\left[\begin{array}{c}S_0\\I_0\end{array}\right]\,,\\
 \\
	y &=& h(x,\theta) = \gamma I\,.
\end{array}\right.
\]
We consider $\Theta=(\theta,x_0)$ since the initial conditions are assumed to be unknown in this example. The Jacobian of the observation with respect to the parameter $\Theta$ is:
\[
\chi(t,\Theta)=\frac{\partial y}{\partial\Theta}(t)=\frac{\partial h}{\partial x}\frac{\partial x}{\partial \Theta}(t)+\frac{\partial h}{\partial \Theta}(t)\,,
\]
This Jacobian has dimension $m\times (p+n) = 1 \times 4$. We have
\begin{eqnarray*}
	\frac{\partial h}{\partial x}&=&\left[\begin{array}{cc} \dfrac{\partial h}{\partial S} & \dfrac{\partial h}{\partial I}\end{array}\right]=\left[\begin{array}{cc}0 & \gamma\end{array}\right]\,,\\
	\frac{\partial h}{\partial \Theta}&=&\left[\begin{array}{cccc} \dfrac{\partial h}{\partial \tilde\beta} & \dfrac{\partial h}{\partial \gamma} & \dfrac{\partial h}{\partial S_0} & \dfrac{\partial h}{\partial I_0}\end{array}\right]=\left[\begin{array}{cccc} 0 & I & 0 & 0\end{array}\right]\,,
\end{eqnarray*}
and 
\[
\frac{\partial x}{\partial \Theta}=\left[\begin{array}{cccc}
	\dfrac{\partial S}{\partial \tilde\beta} & \dfrac{\partial S}{\partial \gamma} & \dfrac{\partial S}{\partial S_0} & \dfrac{\partial S}{\partial I_0}\vspace{.1cm}\\
	\dfrac{\partial I}{\partial \tilde\beta} & \dfrac{\partial I}{\partial \gamma} & \dfrac{\partial I}{\partial S_0} & \dfrac{\partial I}{\partial I_0}\end{array}\right]\,.
\]
This yields
\[
\frac{\partial h}{\partial x}\frac{\partial x}{\partial \tilde\theta}=\gamma\left[\begin{array}{cccc} \dfrac{\partial I}{\partial \tilde\beta} & \dfrac{\partial I}{\partial \gamma} & \dfrac{\partial I}{\partial S_0} & \dfrac{\partial I}{\partial I_0}\end{array}\right]\,.
\]
Therefore,
\[
\chi=\gamma\left[\begin{array}{ccccccc} \dfrac{\partial I}{\partial \tilde\beta} && \dfrac{\partial I}{\partial \gamma} + I && \dfrac{\partial I}{\partial S_0} && \dfrac{\partial I}{\partial I_0}\end{array}\right]\,.
\]
Let $\{t_i\}$, $i=0,1,2,\ldots,M$, be the sampling times. Fisher's Information Matrix is defined as:
\[
\text{FIM}(\Theta,\sigma)=\frac{1}{\sigma^2}\sum_{i=1}^M \chi(t_i,\Theta)^\top\chi(t_i,\Theta)\,,
\]
where $\sigma^2$ is defined as in equation \eqref{estimbias}, since the initial conditions are unknown in this example.

\paragraph{Computing Computing Fisher's Information Matrix} We make the following decomposition:
\[
\frac{\partial x}{\partial \Theta}=\left[\begin{array}{cc} \dfrac{\partial x}{\partial \theta} & \dfrac{\partial x}{\partial x_0}\end{array}\right]=\left[\begin{array}{cc}\displaystyle z & w\end{array}\right]\,,
\]
with
\[
z=\frac{\partial x}{\partial \theta}=\left[\begin{array}{cc}
	\dfrac{\partial S}{\partial \tilde\beta} & \dfrac{\partial S}{\partial \gamma} \vspace{.1cm}\\
	\dfrac{\partial I}{\partial \tilde\beta} & \dfrac{\partial I}{\partial \gamma} \end{array}\right]\quad\mbox{and}\quad w= \frac{\partial x}{\partial x_0}=\left[\begin{array}{cc}
	\dfrac{\partial S}{\partial S_0} & \dfrac{\partial S}{\partial I_0}\vspace{.1cm}\\
	\dfrac{\partial I}{\partial S_0} & \dfrac{\partial I}{\partial I_0}\end{array}\right]\,.
\]
Letting
\[
A(t)=\frac{\partial f}{\partial x}=\left[\begin{array}{cc}
	\dfrac{\partial f_S}{\partial S} & \dfrac{\partial f_S}{\partial I} \vspace{.1cm}\\
	\dfrac{\partial f_I}{\partial S} & \dfrac{\partial f_I}{\partial I} \end{array}\right]=\left[\begin{array}{cc}
	-\tilde\beta I & -\tilde\beta S \\
	\tilde\beta I & \tilde\beta S - \gamma \end{array}\right]\,,
\]
and
\[
B(t)=\frac{\partial f}{\partial \theta}=\left[\begin{array}{cc}
	\dfrac{\partial f_S}{\partial \tilde\beta} & \dfrac{\partial f_S}{\partial \gamma} \vspace{.1cm}\\
	\dfrac{\partial f_I}{\partial \tilde\beta} & \dfrac{\partial f_I}{\partial \gamma} \end{array}\right]=\left[\begin{array}{cc}
	-SI & 0 \\
	SI & -I \end{array}\right]\,.
\]
The FIM can be computed via numerically solving the following system of ODE's:
\[
\left\{\begin{array}{lll}
	\dot x &=& f(x,\theta)\,,\quad x(0)=x_0\,,\\
	\dot z &=& Az+B\,,\quad z(0)=0_{n\times p}\,,\\
	\dot w &=& Aw\,,\quad w(0)=\text{Id}_{n \times n}\,,
\end{array}\right. 
\]
which repeats equation (\ref{eq:fullSystem}). In the following code, the entries of $x$, $z$, and $w$ are indexed in this way:
\[
x=\left[\begin{array}{c} x_1 \\ x_2 \end{array}\right]\,,\quad z=\left(\begin{array}{cc} z_3 & z_5 \\ z_4 & z_6 \end{array}\right)\,,\quad w=\left[\begin{array}{cc} w_7 & w_9 \\ w_8 & w_{10} \end{array}\right]\,,
\]
leading to 
\[
\chi=\gamma\left[\begin{array}{ccccccc}z_4 && z_6 + x_2 && w_8 && w_{10}\end{array}\right]\,.
\]
\section{Numerical implementation} Although the code is very similar the one provided in the previous example (Appendix \ref{appendixCode}), we provide it for convenience, as it required a number of small changes.

\n
We define the following functions in the {\tt Scilab} environment:\m\n

\begin{verbatim}
	function [kguess_n ,SSE]=identifKMK(OBS,T,kguess)
	t0=T(1);m=length(T);
	[x,SSE]=lsqrsolve(kguess,errorKmcK,m,[1.d-8,1.d-8,1.d-5,1d9,0,100]);
	nbparam=length(kguess);
	kguess_n=x;
	x0=[kguess_n(3);kguess_n(4)];
	sol=ode(x0,t0,T,list(KmcK,x(1),x(2)));
	kguess_n=kguess_n(:);
	SSE=sum(SSE.^2);
	xset("window",1);
	sol1=ode(x0,t0,T(1):0.01:T($),list(KmcK ,x(1),x(2)));
	clf;
	plot((T(1):0.01:T($))',kguess_n(2)*sol1(2,:)','k');
	plot(T',OBS ,'ko')
	legend(["$\Large \gamma I(t)$", "$\Large\mbox{Data}$"])
	ylabel("$\Large \mbox{number of deaths per week}$","fontsize",3);
	xlabel("$\Large \mbox{time }t\mbox{ (in week)}$","fontsize",3);
	endfunction
	
	function y=errorKmcK(k,m)
	x0=[k(3);k(4)];
	x0=x0(:);
	B=k(1);
	GAMMA=k(2);
	sol=ode(x0,t0,T,list(KmcK,B,GAMMA));
	predic=GAMMA*sol(2,:);
	predic=predic(:);
	OBS=OBS(:);
	y=OBS-predic;
	endfunction
	
	function  xdot=KmcK(t,x,B,GAMMA)
	xdot=[-B*x(2),0; B*x(2),-GAMMA]*x
	endfunction
	
\end{verbatim}
\n\n

Then, the {\tt Scilab} session goes like this\m\n

\begin{verbatim}[commandchars=\\\{\}]
	--> load('databombay')
	ans  =
	T
	--> OBS=dataBSFlu
	
	OBS  = 
	
	column 1 to 9
	
	8.   10.   12.   16.   24.   48.   51.   92.   124.
	
	column 10 to 16
	
	178.   280.   387.   442.   644.   779.   702.
	
	column 17 to 23
	
	695.   870.   925.   802.   578.   404.   296.
	
	column 24 to 31
	
	162.   106.   64.   46.   35.   27.   28.   24.
	
	--> M=length(OBS);
	
	--> T=0:M-1;
	
	--> S0=15000;I0=7;gamma0=0.6;b0=8e-5;//Initial guesses
	
	--> param=[b0,gamma0,S0,I0];//With b=beta/N
	
	--> p=length(param);
	
	--> [param,SSE]=identifKMK(OBS,T,param,N)
	
	SSE  = 
	
	106336.49
	
	param  = 
	
	0.0000855
	3.7161743
	48113.13
	1.4213612
	
	--> n=2;sigma2=SSE/(M-(n+p))
	
	sigma2  = 
	
	4253.4597
	
	--> B=param(1);GAMMA=param(2);
\end{verbatim} 

\m\n

\medskip

We then compute confidence intervals using the formulas \eqref{FIM}, \eqref{estimtheta}, \eqref{estimbias}, \eqref{estimcovmat}, \eqref{SE}.\m\n

\begin{verbatim}[commandchars=\\\{\}]
	function  FIM=fimKmcK(x0,T,B,GAMMA,sigma2)
	//  Compute  the  sensitivity  matrix
	x0=x0(:);t0=T(1);
	X0=[x0;zeros(4,1);eye(2,2)(:)];
	sol=ode(X0,t0,T,list(JKmcK,B,GAMMA));
	M=GAMMA*sol([4,6,8,10],:);
	M(2,:)=M(2,:)+GAMMA*sol(2,:);
	FIM=M*M'./sigma2;
	endfunction
	
	function  xdot=JKmcK(t,x,B,GAMMA)
	xdot (1)=-B*x(1)*x(2);
	xdot (2)=B*x(1)*x(2)-GAMMA*x(2);
	
	xdot (3)=-B*x(2)*x(3)-B*x(1)*x(4)-x(1)*x(2);
	xdot (4)=B*x(2)*x(3)+(B*x(1)-GAMMA)*x(4)+x(1)*x(2);
	xdot (5)=-B*x(2)*x(5)-B*x(1)*x(6);
	xdot (6)=B*x(2)*x(5)+(B*x(1)-GAMMA)*x(6)-x(2);
	
	xdot (7)=-B*x(2)*x(7)-B*x(1)*x(8);
	xdot (8)=B*x(2)*x(7)+(B*x(1)-GAMMA)*x(8);
	xdot (9)=-B*x(2)*x(9)-B*x(1)*x(10);
	xdot (10)=B*x(2)*x(9)+(B*x(1)-GAMMA)*x(10);
	endfunction
\end{verbatim}

\medskip

Then the {\tt Scilab} session is 

\begin{verbatim}[commandchars=\\\{\}]
	
	--> x0=[param(3);param(4)];
	
	--> FIM=fimKmcK(x0,T,B,GAMMA)
	FIM  = 
	
	1.100D+14  -2.117D+09   196924.7    92108921.
	-2.117D+09   40885.251  -3.7835134  -1826.2845
	196924.7   -3.7835134   0.0003533   0.1603199
	92108921.  -1826.2845   0.1603199   118.22324
	
	--> cond(FIM)
	ans  =
	
	9.141D+24
	
	--> CovMAT=inv(FIM)
	Warning: Matrix is close to singular or badly scaled. 
	CovMAT  = 
	
	0.0000006   0.0093978  -220.96547  -0.0128142
	0.0093978   150.37071  -3535242.5  -204.997
	-220.96547  -3535242.5   8.313D+10   4820745.
	-0.0128142  -204.997     4820745.    279.63056
	
	--> t=cdft("T",M-(n+p),0.975,0.025)
	t  = 
	
	2.0595386
	
	--> confB=t*sqrt(CovMAT(1,1))
	confB  = 
	
	0.0015784
	
	--> confGAMMA=t*sqrt(CovMAT(2,2))
	confGAMMA  = 
	
	25.255243
	
	--> confS0=t*sqrt(CovMAT(3,3))
	confS0  = 
	
	593794.26
	
	--> confI0=t*sqrt(CovMAT(4,4))
	confI0  = 
	
	34.439929
	
\end{verbatim}

\chapter{Generalized Least Squares}\label{app:GLS}
With Ordinary Least Squares, constant variance has been assumed which may be not appropriate for some data. A relative error, i.e., when the error is assumed to be proportional to the size of the measurement, is an assumption that might be reasonable when counting individuals in a population.

\m
\n
In this case we assume that the observation are \cite{Banks2009,MR2532016,MR3203115,Capaldi:2012aa}:

\[Y_i=y(t_i,\Theta) + y(t_i,\Theta)^\rho \, \mathcal E_i\,,\]
with $\Theta=(\theta,x_0)$.

\n
The criterion to be minimized is 

\[\mathcal J(\Theta)= \sum_{i=1}^N \; w_i \;\left[Y_i - y(t_i,\Theta) \right]^2. \]

\n
The values of the weights ($w_i$) depend on the value of the model and are not known. The process is carried with an iterated re-weighted least squares:\m

\begin{enumerate}
\item Estimate $\hat \Theta_0$ with an OLS step ($\rho=0$): $w_i=1$ for all $i=1,\ldots,N$
\item Set $\rho=1$ (for instance) and $w_i=1/\left[y(t_i,\hat \Theta_0)\right]^{2\,\rho}$ for all $i=1,\ldots,N$
\item  Form $\mathcal J(\Theta)$ with these $w_i$ and estimate
\[\hat \Theta_1=\arg \min_{\Theta}  \mathcal J(\Theta)\]
\item Continue the procedure till  the estimates $\hat \Theta_k$ and $\hat \Theta_{k+1}$ are sufficiently close to each other, to obtain $\hat \Theta_{\text{GLS}}$.
\end{enumerate}
\m
With $\hat \Theta_{\text{GLS}}$, as in the OLS case, we can obtain the covariance matrix ($\Sigma$) of the estimated parameters, approximated as the inverse Fisher Information Matrix with weights:
\begin{equation*}
 \text{FIM}  (\hat  \Theta_{\text{GLS}},\hat\sigma_{\text{GLS}}^2) = \\\frac{1}{\hat\sigma^2_{\text{GLS}}}\sum_{i=1}^N  \;  \dfrac{1}{y(t_i,\hat \Theta_{\text{GLS}})^{2\,\rho}}\; \dfrac{\partial y}{\partial \Theta} (t_i, \hat \Theta_{\text{GLS}})^\top \;    \dfrac{\partial y}{\partial \Theta} (t_i, \hat \Theta_{GLS})\,,
\end{equation*}
with
\[
\hat \sigma^2_{\text{GLS}}=  \dfrac{1}{N-p} \; \sum_{i=1}^N \;   \dfrac{1}{y(t_i,\hat \Theta_{\text{GLS}})^{2\,\rho}}\;\left[ Y_i - y(t_i, \hat  \Theta_{\text{GLS}}) \right]^2. 
\]
We then obtain
\[
\hat\Sigma_{\text{GLS}}= \left[\text{FIM}(\hat\Theta_{\text{GLS}},\hat\sigma_{\text{GLS}}^2)\right]^{-1},
\]

The square roots of the diagonal elements of the approximation of the covariance matrix $\hat\Sigma_{\text{GLS}}$ give the standard errors.

\backmatter
\printindex


\end{document}